\documentclass[12pt,reqno,twoside]{amsart}
\usepackage[utf8]{inputenc}
\usepackage{amsmath}
\usepackage{tgtermes}

\usepackage{indentfirst,csquotes}

\usepackage{amssymb}
\usepackage[T1]{fontenc}
\usepackage[english]{babel}
\usepackage{fancyhdr}
\usepackage{indentfirst}
\usepackage{graphicx}
\usepackage{newlfont}
\usepackage{amssymb}
\usepackage{latexsym}
\usepackage{amsthm}
\usepackage{mathtools}
\usepackage{subfig}
\usepackage{floatflt}
\usepackage{floatrow}
\usepackage{color}
\usepackage{textcomp}
\usepackage{framed}
\usepackage{verbatim}
\usepackage{hyperref}
\usepackage{booktabs}
\usepackage{caption}
\usepackage[mathscr]{euscript}
\usepackage{braket}
\usepackage{bbm}
\usepackage{mathrsfs}
\newcommand*{\transp}[2][-3mu]{\ensuremath{\mskip1mu\prescript{\smash{\mathrm t\mkern#1}}{}{\mathstrut#2}}}%
\DeclarePairedDelimiter{\abs}{\lvert}{\rvert}
\DeclarePairedDelimiter{\norma}{\lVert}{\rVert}

\author[]{Davide Tramontana}
\title[]{Subelliptic Random Walks on Riemannian Manifolds and Their Convergence to Equilibrium}

\numberwithin{equation}{section}

\theoremstyle{plain}
\newtheorem{theorem}{Theorem}[section]
\newtheorem{lemma}[theorem]{Lemma}
\newtheorem{corollary}[theorem]{Corollary}
\newtheorem{proposition}[theorem]{Proposition}

\theoremstyle{definition}
\newtheorem{definition}[theorem]{Definition}

\newtheorem{exa}[theorem]{Example}

\newtheorem{remark}[theorem]{Remark}

\newtheorem*{theorem*}{Theorem}

\newtheorem{hp}{}

\address{Davide Tramontana
\newline \indent Dipartimento di Matematica, Universit\`a di Bologna
\newline \indent Piazza di Porta San Donato 5, 40126, Bologna, Italy}
\email{davide.tramontana4@unibo.it}
\thanks{The author is a member of the Research Group GNAMPA of INdAM}

\begin{document}
\begin{abstract}
The aim of this work is to study the convergence to equilibrium of an $(h,\rho)$-subelliptic random walk on a closed, connected Riemannian manifold $(M,g)$ associated with a subelliptic second-order  differential operator $A$ on $M$. 
In such a random walk, $h$ roughly represents the step size and $\rho$ the speed at which it is carried out. 
To construct the random walk and prove the convergence result, we employ a technique due to Fefferman and Phong, which reduces the problem to the study of a constant-coefficient operator $\tilde{A}$ that is locally equivalent to our second-order  subelliptic operator $A$, in the sense that the diffusion generated by $\tilde{A}$ induces a local diffusion for $A$. By using the compactness of $M$ this local diffusion can be lifted to a global diffusion, and the convergence result is then obtained via the spectral theory of the associated Markov operator. 
\end{abstract} 

\maketitle

\tableofcontents                        
\rhead[\fancyplain{}{\bfseries\thepage}]{\fancyplain{}{\bfseries\thepage}}
\lhead[\fancyplain{}{\bfseries\thepage}]{\fancyplain{}{\bfseries\thepage}}

\section{Introduction}		
The aim of this work is to construct an $(h,\rho)$-\textit{subelliptic diffusion} on a closed, connected smooth Riemannian manifold $
(M,g)$ and to establish convergence to equilibrium through the spectral theory of the associated Markov operator.
Specifically, starting from a subelliptic operator $A$ on $M$ we construct a Random Walk (see Definition \ref{def.ranwalk})  associated with it (cf. \cite{DLM1}, \cite{DLM2}, \cite{LM}, \cite{BNR}), which depends on two small parameters $h,\rho$ representing, respectively, the step size and the speed at which the step is performed.

This random walk can be interpreted as a \emph{subelliptic} version of the \textit{Gibbs–Metropolis algorithm}, a Markov chain widely used in statistical mechanics to reach low-energy configurations of many-body systems.
More precisely, the Gibbs-Metropolis algorithm proposes local moves that are always accepted if they reduce the system’s energy, and accepted with decreasing probability as the energy increases, according to the Boltzmann distribution.
Therefore, it is commonly used to study, for instance, liquid - solid phase
transitions in particle systems, as well as magnetic phase transitions in spin models (see for instance \cite{FS}, \cite{Lo} and the references therein).

On the other hand, from a mathematical perspective this work originated from the aim of generalizing the construction of a \textit{hypoelliptic random walk} developed by G. Lebeau and L. Michel in \cite{LM}. They considered the case of differential operators
on $M$ that can be written as $A=\sum_{k=1}^p X_k^2$,
where $X_1, \dots, X_p$ are divergence-free vector fields (with respect to a fixed measure $\mu$) 
and satisfy H\"ormander's condition (the term \textit{hypoelliptic} refers to this fact). 
In this context, the main novelty of the subelliptic construction we present
is that it applies 
to second-order  differential operator on $M$, which, in general, cannot be written as sums of squares of smooth vector fields. Namely, we consider a subelliptic operator $A$ of order $0<\varepsilon<1$ (see Definition \ref{def.subopM}) that
can be locally written as (with $(U,\varphi)$ local chart and $\tilde{U}=\varphi(U)$)
\[
A_U=-\sum_{i,j=1}^n a_{ij}(\tilde{x})\frac{\partial^2}{\partial {\tilde{x}_i}\partial {\tilde{x}_j}}+\sum_{k=1}^nb_{k}(\tilde{x})\frac{\partial}{\partial \tilde{x}_k}+d(\tilde{x}), \quad \tilde{x} \in \tilde{U},
\]
where $a_{ij}$, $b_k$, 
$d \in C^\infty(\tilde{U}_\alpha;\mathbb{\mathbb{R}})$
for all $i,j,k \in \lbrace 1,...,n \rbrace$ and $(a_{ij})_{i,j} \geq 0$.
Therefore, with our approach
we are able to include any
second-order  differential operators
(which satisfy the natural hypotheses \ref{hp1}, \ref{hp2}, \ref{hp3}) on a closed, connected Riemannian manifold.

This work is essentially divided into two main parts. The first is devoted to the construction of a global subelliptic diffusion, while the second focuses on establishing convergence results for this diffusion, obtained via an analysis of the spectral theory of the associated Markov operator. In what follows, we provide a brief overview of these two parts.

The diffusion we consider is based on the notion of subunit 
balls $B_A(x,\rho)$ (see Definition \ref{def.subball})
associated with our subelliptic operator $A$.
Loosely speaking, one moves from a point $x \in M$ to a point $x' \in B_A(x,\rho)$,
for sufficiently small $\rho$, along a (perturbation of) subunit curve for $A$, that is a curve $\gamma_A:[0,1] \rightarrow M$ from $x$ to $x'$ such that the velocity vector is controlled 
from above by the principal symbol $a_2(x,\xi)$ of $A$, almost everywhere. More precisely, we look for a curve $\gamma_A:[0,1] \rightarrow M$ that satisfies
\[
\begin{split}
&\gamma_A(0)=x, \quad \gamma_A(1)=x', \\
&\abs{\dot{\gamma}_A(t,\xi)}^2 \leq a_2(\gamma(t),\xi), \quad \text{for} \ \text{a.e.} \ t \in [0,1], \ \forall \xi \in T^\ast_{\gamma(t)}M,
\end{split}
\]
where $a(x,\xi)$ is the principal symbol of $A$ (see Remark \ref{rmk.princdiff}) 
and $\dot{\gamma}_A(t,\xi)$ is the principal symbol of the vector fields 
$\dot{\gamma}_A(t)$ along the curve (see Remark \ref{princvecfield}).
The first basic construction that we develop is a local construction based on the Oleinik–Radkevich vector fields. However, as discussed in Remark \ref{rmk.notopt}, this approach is not optimal. This motivates the introduction of a new more refined construction based on an iterative procedure, which we next explain.
The idea is to reduce locally our operator $A$ to an equivalent (the term equivalent should be more clear soon) constant coefficients operator $\tilde{A}$, which will be our universal model. Therefore, we reduce matters to constructing a diffusion for such an operator and then we "pullback" this diffusion to a local diffusion on our manifold $M$. Finally, using the compactness of $M$ we are able to extend this local construction to a global one. 

To employ this procedure we use a technique of Fefferman and Phong 
based on the Calder\' on-Zygmund localization 
(see \cite{F}, \cite{FP1}, \cite{FP2}, \cite{FSC}).
To briefly explain this technique, since
it is of local nature, we will reduce to the case
of the standard cube $Q \subseteq \mathbb{R}^n$ 
with dilations $Q^\ast$ and $Q^{\ast\ast}$ 
(see Section \ref{sec.precomm} for the notation).
Essentially,
by choosing the smallest $\delta_1=2^{-N_1}$ for which 
\[
\max_{1 \leq i \leq n}\max_{x \in \overline{Q(0,\delta_1)}}\abs*{a_{ii}(x)} \geq 10R_1\delta_1^2, 
\]
with $R_1>0$ a sufficiently large universal constant, it
is possible to construct a diffeomorphism 
\[
\begin{split}
\phi_1:]-r_1,r_1[ \times Q_{n-1}(0,1) & \rightarrow \phi_1 (]-r_1,r_1[ \times Q_{n-1}(0,1)) \subseteq Q(0,\delta_1) \\
(w_1,z) \mapsto x,
\end{split}
\]
for some universal $r_1>0$, such that,
in these new coordinates, the operator $A$ looks
like
\[
A_{\phi_1}=-\partial_{w_1}^2 +\sum_{i,j=2}^n b_{i,j}(w_1,z) \partial_{z_i} \partial_{z_j}.
\]
At this point the key tool is that, for $\rho>0$ sufficiently small, one has
\[
\lbrace \abs*{w_1}\leq \rho \rbrace \times B_{\bar{A}_{\phi_1}}(0,\rho) \subseteq B_{A_{\phi_1}}(0,C\rho), \quad C>0 \ \ \text{universal},
\]
where, 
\[
\bar{A}_{\phi_1}=-\sum_{i,j=2}^n\bar{b}_{ij}(z)\frac{\partial^2}{\partial z_i z_j}, \quad \bar{b}_{ij}(z)=\frac{1}{2\rho}\int_{-\rho}^{\rho}b_{ij}(w_1,z)dw_1.
\]
Therefore, by iterating this argument we obtain that locally
\[
A  \sim \tilde{A}=-\sum_{i=1}^nc_j^2 \rho^{2\kappa_j}\frac{\partial^2}{\partial w_j^2}, 
\]
where $1=\kappa_1 \leq \kappa_2\leq \dots \leq \kappa_n=1/\varepsilon$ (recall that $\varepsilon$ is the order of subellipticity of $A$) and $c_1,\dots,c_n$ are universal constants. This equivalence should be understood in the sense that
subunit curves for $\tilde{A}$ induce subunit curves for $A$ and vice versa.
For our purposes, it is important to emphasize the fact that subunit curves for $\tilde{A}$
induce subunit curves for 
$A$. We highlight this in the following theorem, which is one of the main results of the paper.
\begin{theorem}\label{fundtheorem}
Let $A$ be a
second-order differential operator of the form 
\[
A=-\sum_{i,j=1}^n a_{ij}(x)\frac{\partial^2}{\partial {x_i}\partial {x_j}}+\sum_{k=1}^nb_k(x)\frac{\partial}{\partial x_k}+d(x),
\]
where $a_{ij},b_k,d \in C^\infty(Q^{\ast\ast},\mathbb{R})$
for every $i,j,k=1,...,n$ and $\mathsf{A}_2(x)=(a_{ij}(x))_{i,j}$
is symmetric and positive semidefinite for any $x \in Q^{\ast\ast}$.

\noindent Then, if $A$ is subelliptic of order $\varepsilon$
$(0<\varepsilon<1)$, there exists a universal block
\[
Q_\rho=]-c_1\rho,c_1\rho[ \times ]-c_2\rho^{\kappa_2},c_2\rho^{\kappa_2}[\times \dots \times ]-c_n\rho^{\kappa_n},c_n\rho^{\kappa_n}[,
\]
where $\rho>0$ is sufficiently small, $c_1=1,\dots,c_n>0$ are universal constants and 
$\kappa_1=1\leq \kappa_2\leq \dots \leq \kappa_n=1/\varepsilon$,
such that for any 
$x_0 \in Q$ there exists a diffeomorphism,
\[
\Phi:Q_\rho \rightarrow \Phi(Q_\rho), \quad \Phi(0)=x_0,
\]
with uniformly bounded derivatives.

\noindent In addition, there exists a second-order differential
operator with constant coefficients on this
universal block defined as
\[
\tilde{A}:=-\sum_{j=1}^n c_j^2\rho^{2\kappa_j} \frac{\partial^2}{\partial w_j^2},  
\]
with the following property.

\noindent If we define $\gamma_j:[0,1]\rightarrow \mathbb{R}^n$,
for $1\leq j \leq n$ as
\[
\gamma_j(t)=e^{tc_j\rho^{\kappa_j}\frac{\partial}{\partial w_j}}(0),
\]
we have that this curve is a subunit curve for $\tilde{A}$,
starting at $0$ with speed less than $1$ and it determines
a curve $\gamma:[0,1]\rightarrow Q^\ast$,
subunit for $A$, with speed less than $C\rho$,
for some universal constant $C>0$, such that 
\[
\gamma(0)=x_0, \quad \quad \gamma(1)=:x \in B_{A}(x_0,C\rho).
\]
\end{theorem}
To make this construction global,  we use the uniform boundedness of the coefficients along with some geometric arguments, which lead to the following universal inclusions. 
\begin{theorem}\label{theosubM}
There exist $\nu,\rho,C>0$ universal constants,
such that, for every $x \in M$, 
\begin{equation}
\label{geomrelationM}
\Phi(Q_\rho)\subseteq \Phi(Q_\rho^\ast)\subseteq B_{A_\ell}(\tilde{x},C\rho)\subseteq Q(\tilde{x},\nu),
\end{equation}
where, $\tilde{x}=\varphi_{\ell}(x)\in \mathbb{R}^n$
for some $\ell=1,\dots,m$,
$A_\ell$ is the operator corresponding to $A$
on the open set $\tilde{V}_\ell \subseteq \mathbb{R}^n$,
$Q_\rho$ is the universal block of radius $\rho$ (see Definition \ref{univcube}),
$Q_\rho^\ast$ is a fixed dilation of $Q_\rho$ (i.e. $Q^\ast_\rho=Q_{c^\ast\rho}$, with $c^\ast>1$)
and $\Phi$ is a diffeomorphism with uniform bounds.
\end{theorem}
Therefore, we fix $\rho_0>0$ sufficiently small and for $\rho \in ]0,\rho_0]$ 
we construct an atlas $\lbrace (Q_{i},\Psi_i) \rbrace_{i=1}^N$ of $M$ where the local open sets $Q_i$ are subunit cubes - i.e. image (through $\Psi_i^{-1}$) of an universal block $Q_\rho$ of $\mathbb{R}^n$ - (see Definition \ref{def.subM}). Hence for $h \in ]0,h_0]$, with $h_0$ sufficiently small,
we define the subunit $(h,\rho)$-diffusion as follows.\\
$\bullet$ \textit{Step 0}. Fix $x_0 \in M$.\\
$\bullet$ \textit{Step k}. At step $k \in \mathbb{N}$ we choose at random $i \in \lbrace 1,...,N \rbrace$, $j \in \lbrace 1,...,n \rbrace$ and $t \in ]-h,h[$ (with respect the uniform probability) and set 
\[
x_{k+1}=
\begin{cases}
\Phi_{i,j,h,\rho}(t,x_k), & \text{if $x_k\in Q_{i}$ and $\Phi_{i,j,h,\rho}(t,x_k)\in Q_{i}$} \\
x_k, & \text{otherwise},
\end{cases}
\]
where $\Phi_{i,j,h,\rho}(t,x_k)$ is the flow of $X_{i,j}$, the divergence-free 
version of the vector field $\partial^{(i)}_{x_j}$ with respect to the coordinates induced by $\Psi_i$.
In order to obtain the convergence to equilibrium for this diffusion 
we study the spectral theory of the self-adjoint $(h,\rho)$-Markov operator $S_{h,\rho}$ associated with it (see \eqref{defShrho} and the discussion above for the construction).

We express this operator
in term of its Markov-Riesz kernel 
$s_{h,\rho}(x,dy)$ (see \eqref{eq.MR}) as
\[
S_{h,\rho}f(x)=\int_M f(y)s_{h,\rho}(x,dy), \quad f \in C^\infty(M),
\]
and analogously, for $k \in \mathbb{N}$, we denote by $s_{h,\rho}^k(x,dy)$ the Markov-Riesz kernel of the iterated operator $S_{h,\rho}^k=\underbrace{S_{h,\rho} \circ \dots \circ S_{h,\rho}}_{\text{k-times}}$.

The key tools to reach the spectral results and the convergence to equilibrium are Theorem \ref{lowestimatekernel} 
and Theorem \ref{theohighlow}.
In the first we obtain 
a lower bound for our $(h,\rho)$-Markov kernel, which yields an $L^\infty$ upper bound for the eigenfunctions of $S_{h,\rho}$, associated with energy levels close to $1$. In the second 
we decompose the $L^2$ functions with bounded Dirichlet energy into low and high frequency parts, $u=u^L+u^H$, where $u^H \in H^1(M)$ (the Sobolev space on $M$ of $L^2$ distribution with $L^2$ derivatives) and $u^L$
has small $L^2$ norm.
This decomposition is the key tool to obtain a Weyl-type estimate (see Subsection \ref{sec.weylest}) which yields 
the following results for the spectrum of the operators $S_{h,\rho}$ (Theorem \ref{Theorem1}). The latter, in turn, leads to convergence to equilibrium (Theorem \ref{Theorem2})
i.e. the convergence of the $(h,\rho)$-Markov-Riesz kernel of the iterated $(h,\rho)$-Markov operator to a
fixed stationary 1-density on $M$, with respect to the so-called total variation distance, when the number of the step of the walk goes to infinity. Finally, it is interesting to note that the 
rate of convergence is given by the spectral gap $g_{h,\rho}$ of $I-S_{h,\rho}$, which corresponds distance between the eigenvalue $1$ of $S_{h,\rho}$ and the rest of the spectrum, and satisfies $g_{h,\rho} \approx h^2$. Recall that here $\rho \in ]0,\rho_0]$, where $\rho_0>0$ is chosen sufficiently small to ensure the validity of the previous construction.

\begin{theorem}\label{Theorem1}
Let $\rho \in ]0,\rho_0]$ be fixed.
There exists $h_0>0$ such that, for any $h \in ]0,h_0]$ the  following facts hold. \\
$(i)$ $S_{h,\rho}1=1$ and $1$ is a simple eigenvalue of the operator $S_{h,\rho};$ \\
$(ii)$ There exists a fixed constant $\delta_1>0$ such that $\mathrm{Spec}(S_{h,\rho}
)\subseteq [-1+\delta_1,1]$; \\
$(iii)$ There exists a fixed constant $\delta_2>0$ such that $\mathrm{Spec}(S_{h,\rho}) \cap [1-\delta_2,1]$ is discrete. Moreover, for any $0 \leq \zeta \leq \delta_2h^{-2}$,
\[
\sharp\Bigl(\mathrm{Spec}(S_{h,\rho})\cap [1-h^2\zeta,1]\Bigr)\leq C_1(1+\zeta)^m,
\]
for some fixed constants $C_1,m>0$;\\
$(iv)$ There exist two universal constants $C_2,C_3>0$ such
that the spectral gap $g_{h,\rho}$ of the operator $I-S_{h,\rho}$ 
(see Definition \ref{def.specgap}) satisfies \\
\begin{equation*}
C_2h^2 \leq g_{h,\rho} \leq C_3h^2. 
\end{equation*} 
\end{theorem}

\begin{theorem}\label{Theorem2}
Let $\rho \in ]0,\rho_0]$ be fixed. There exists $h_0>0$ such that for any $h \in ]0,h_0]$ one has 
\begin{equation}\label{convequi}
\sup_{x \in M} \norma*{s_{h,\rho}^k(x,dy)-\mu}_{TV}\leq Ce^{-kg_{h,\rho}},
\end{equation} 
for any $k \in \mathbb{N}$, where $C>0$ is a universal constant.
\end{theorem}

The plan of the paper is the following. In Section \ref{sec.subopM}, we recall basic facts about subelliptic differential operators on smooth manifolds, along with the notion of subunit geometry.
Section \ref{sec.subwalk} is devoted to the construction of the 
$h$-subelliptic random walks, first via the Oleinik–Radkevich vector fields, and then using the Fefferman–Phong localization argument.
In Section \ref{sec.aux}, we prove several auxiliary results that play a key role in establishing the spectral properties of the associated operators.
Section \ref{sec.maintheo} is dedicated to the proofs of Theorem \ref{Theorem1} and Theorem \ref{Theorem2}.
An appendix, in which we collect some basic results from differential geometry and spectral theory and prove a few technical lemmas, concludes the paper.

\textbf{Notation}
 Let $(M,g)$ be a connected, closed Riemannian manifold endowed with a $1$-density $d\mu$ to be chosen later. In some situations - as we will specify - the basic definitions apply to a general smooth manifold $M$.
 By $(\cdot, \cdot)_0$ we mean the $L^2$-product on $M$ (with respect the fixed $1$-density $d\mu$) and by $ (\cdot, \cdot)_s$ the Sobolev product of order $s \in \mathbb{R}$ on $M$.
 By $\langle \cdot, \cdot \rangle$ we mean the standard Euclidean product. 
We will use the notation $\langle \cdot \rangle=(1+\abs{\cdot}^2)^{1/2}$.
We write $a\approx b$ if there exist two universal constants $c,c'>0$ such that $a\leq cb$ and $b\leq c'a$.
We denote by $\mathbb{N}_0=\mathbb{N} \cup \lbrace 0 \rbrace$.
By $\Gamma(M,E)$ we mean the set of sections of the fiber bundle $E$ over $M$.
We denote by $\mathscr{B}$ the $\sigma$-algebra of Borel on $M$, by $\mathscr{P}(M)$ the set of the probability measures on $(M,\mathscr{B})$, by $\mathscr{M}(M,\mathbb{R})$
the set of real-valued (with respect the $\sigma$-algebra of Borel) measurable functions on $M$ and finally by $\mathscr{M}_b(M,\mathbb{R})$ the set of real-valued, bounded measurable functions on $M$.   
If $f,g:I\subseteq \mathbb{R} \rightarrow \mathbb{C}$ we write 
$f=O_{\mathrm{unif}}(g)$ if there exists a constant $C>0$ such that $\abs{f(x)} \leq C \abs{g(x)}$ for all $x \in I$. 
We denote by $\mathscr{S}(\mathbb{R}^n)$ the space of Schwartz functions and by $\mathscr{S}'(\mathbb{R}^n)$ the space of tempered distributions.
If $U,V \subseteq \mathbb{R}^n$, we write $U \subseteq V$ if there exists a universal constant $C>0$ such that $C^{-1}U\subseteq V \subseteq CU$, where $cU=\lbrace y \in \mathbb{R}^n; \quad y=cx, \quad x \in U \rbrace$, for $c>0$.

\section{Subelliptic Calculus on Manifolds}\label{sec.subopM}
The purpose of this chapter is to introduce the basics
of the \textit{subelliptic} calculus for a second-order
differential operator on a smooth closed
Riemannian manifold of dimension $n$.
In the first section, we give the notion of subellipticity
and establish our setting of work.
Then, in the second section, we review some fundamental
geometric tools related to this notion,
which will play a crucial role in studying our diffusion problems. 

\subsection{Subelliptic Operators on Manifolds and Main Hypotheses}

Let $(M,g)$ be a connected, closed Riemannian manifold of dimension $n$,
with $\mathscr{U}=\lbrace (U_\alpha,\varphi_\alpha); \ \alpha \in \mathscr{A} \rbrace$ the
differentiable structure
on $M$, where $\varphi_\alpha : U_\alpha \rightarrow \tilde{U}_\alpha\subseteq \mathbb{R}^n$
are the local maps. 
 
The aim of the work is to study a diffusion associated
with a second-order differential operator on $M$
(by differential operator on $M$
we mean an operator defined as in Definition \ref{diffM})  
that satisfy a subellipticity condition stated as follows. 

\begin{definition}\label{def.subopM}
Let $A$ be a second-order differential operator on $M$. The operator $A$ is said to be \textit{subelliptic of order $\varepsilon$}, $0<\varepsilon<1$, if there exist some constants $c,C>0$, such that, 
\begin{equation}\label{sub1}
\norma*{Au}_0+C \norma*{u}_0 \geq c \norma*{u}_{2\varepsilon}, \quad \forall u \in C^\infty (M).
\end{equation} 
\end{definition}

Throughout the work we will deal with a second-order differential operator
defined on a closed manifold, but,
since we usually work locally,
it is convenient to clarify early what we mean by a \textit{local expression} of our operator. 

Let $A$ be a second-order differential operator on $M$.
By definition (see again \ref{diffM}), for every $\alpha \in \mathscr{A}$, we have that
the operator 
\[
A_\alpha f:=A(f \circ \varphi_\alpha^{-1}) \circ \varphi_\alpha, \quad f \in C^\infty(\tilde{U}_\alpha)
\]
(that we call \textit{local espression of $A$}
on $\tilde{U}_\alpha$) 
has the form 
\[
A_\alpha=-\sum_{i,j=1}^n a^{(\alpha)}_{ij}(\tilde{x})\frac{\partial^2}{\partial {\tilde{x}_i}\partial {\tilde{x}_j}}+\sum_{k=1}^nb^{(\alpha)}_{k}(\tilde{x})\frac{\partial}{\partial \tilde{x}_k}+d^{(\alpha)}(\tilde{x}), \quad \tilde{x} \in \tilde{U}_\alpha,
\]
where $a^{(\alpha)}_{ij}$, $b^{(\alpha)}_k$, 
$d^{(\alpha)} \in C^\infty(\tilde{U}_\alpha,\mathbb{C})$
for all $i,j,k \in \lbrace 1,...,n \rbrace$.
In addition, we denote by 
$\mathsf{A}_2^{(\alpha)}=(a_{ij}^{(\alpha)})_{i,j=1,\dots,n}$
the matrix of the coefficients
of the principal part of the local expression $A_\alpha$.

We will work throughout with a subelliptic
second-order differential operator $A$
on $M$ under the following hypotheses.
\begin{hp}\label{hp1}
For every local expression $A_\alpha$ we assume that
the associated coefficients $a^{(\alpha)}_{ij},b^{(\alpha)}_k,d^{(\alpha)}$ are
real valued functions, for all $i,j,k \in \lbrace 1,...,n \rbrace$, 
and the matrix of coefficients $\mathsf{A}^{(\alpha)}_2(\tilde{x})=(a^{(\alpha)}_{ij}(\tilde{x}))_{i,j=1,...,n}$ is
symmetric and positive semidefinite for every $\tilde{x} \in \tilde{U}_\alpha$.
\end{hp}  
\begin{hp}\label{hp2}
$A$ is a formally self-adjoint operator.
\end{hp}
\begin{hp}\label{hp3}
For every local expression $A_\alpha$ we assume that
the coefficients $a^{(\alpha)}_{ij}$ does not vanish on open sets,
that is,
for all $i,j \in \lbrace 1,\dots n \rbrace$ and 
for all given $W\subseteq \mathbb{R}^n$
one has 
\[
{a_{i,j}^{(\alpha)}}_{\bigl|_{W}}=0 \Rightarrow \mathrm{Int(W)=\emptyset}.
\]
\end{hp}
Finally, we wish to emphasize that the subellipticity condition is weaker than the ellipticity one in the following sense. 
\begin{remark}
By Proposition \ref{prop.tec.1} we have that \eqref{sub1} holds if and only if  
\begin{equation}
\label{sub2}
\mathrm{Re}(Au,u)_0+C' \norma*{u}_0^2 \geq  c' \norma*{u}^2_{\varepsilon}, \quad \forall u \in C^\infty(M).
\end{equation}
Therefore, since the latter with $\varepsilon=1$ is the celebrated
G\r{a}rding inequality (see for instance Taylor \cite{Tay}, p. 55),
we have that if $A$ is an elliptic operator of order two (see \cite{Tay}, p. 60), 
it satisfies \eqref{sub2},
and consequently \eqref{sub1}, with $\varepsilon=1$. 
\end{remark}

\subsection{Subunit Vector Fields}
The aim of this subsection is to give some geometric tools
related to the subellipticity notion. 
To do that we refer to \cite{FP2}, \cite{FSC}, \cite{JSC}. 
The key tool will be the so-called \textit{subunit ball} associated
with a subelliptic operator $A$ and to introduce it we first need
to recall the notion of \textit{subunit vector field} associated with $A$,
satisfying
the hypotheses \ref{hp1}, \ref{hp2} and \ref{hp3}.

In what follows, we denote by
$a_2(x,\xi) \in C^\infty(T^\ast M\setminus 0)$, the principal symbol of $A$
(see Remark \ref{rmk.princdiff}), as the function that corresponds locally to 
\[
a^{(\alpha)}_2(\tilde{x},\tilde{\xi}):=\langle \mathsf{A}_2^{(\alpha)}(\tilde{x})\tilde{\xi},\tilde{\xi} \rangle, \quad (\tilde{x},\tilde{\xi}) \in \tilde{U}_\alpha \times \mathbb{R}^n.
\] 
Moreover, given a smooth vector field $X \in \Gamma(M,TM)$
we denote by $X(x,\xi)$ its principal symbol (see Remark \ref{princvecfield}). 
\begin{definition}
Let $X \in \Gamma(M,TM)$ and let $A$ be a subelliptic 
operator on $M$ with the aforementioned properties.
We say that \textit{$X$ is subunit for $A$ at $x \in M$} if 
\begin{equation}
\label{SubVe}
\abs*{X(x,\xi)}^2 \leq a_2(x,\xi), \quad \forall \ \xi \in T_x^\ast M.
\end{equation}
\end{definition}
The subunit condition could be equivalently
stated using the fixed Riemannian metric $g$ on $M$, or in coordinates,
as shown below. 
\begin{remark}
\label{rmkmetricprincsymb}
If we fix $x \in M$ and we take $\xi \in T_x^\ast M$, there exists $v_\xi \in T_xM$ such that 
\[
g_x(v_\xi,v)=\xi(v), \quad \forall v \in T_xM.
\]
Thus, we have
\[
X(x,\xi)=g_x(v_\xi,X(x))=\xi(X(x)), \quad \forall \xi \in T^\ast_xM.
\]
Hence, in particular, condition \eqref{SubVe} is equivalent to 
\[
g_x(X(x),v_\xi)^2 \leq a_2(x,\xi), \quad \forall \ \xi \in T_x^\ast M.
\]
\end{remark}
The notion of subunit vector field allows us to give
a notion of \textit{subunit curve}.
Roughly speaking, a subunit curve for $A$,
is a curve for which the velocity vector is controlled
from above by the principal symbol of $A$, for almost every time.

\begin{definition}
Let $\gamma:[0,1]\rightarrow M$ be a Lipschitz curve.
We say that \textit{$\gamma$ is a subunit curve}
for $A$ if $\dot{\gamma}(t)$ is subunit for $A$ for $\mathrm{a.e.t} \in [0,1]$. 
\end{definition}

Finally, using the notion of subunit curve,
it is possible to define a new distance on $M$,
the so-called \textit{subunit distance} that yields
to the notion of \textit{subunit ball}. 

\begin{definition}\label{def.subball}
Let $A$ be a subelliptic operator defined on $M$.
\textit{The subunit distance associated with $A$}
between two points $x,x' \in M$ is defined as 
\begin{center}
$d_A(x,x')=\inf \lbrace r>0; \ \exists \ \gamma:[0,1]\rightarrow M$ subunit curve for $r^2A$ and $\gamma(0)=x, \gamma(1)=x' \rbrace$
\end{center}
and \textit{the subunit ball associated with $A$ centered at $x \in M$ of radius $\rho>0$}
is defined as 
\[
B_A(x,\rho)=\lbrace x' \in M; \ d_A(x,x')<\rho \rbrace.
\]
\end{definition}
\begin{remark}\label{subballChangeofcoordinates2}
Note that if we consider two open sets $U,V \subseteq \mathbb{R}^n$,
a subelliptic operator $A$ defined on $U$ and a diffeomorphism $\Phi: V\rightarrow U$,
we have that setting $A_{\Phi}f:=A(f\circ \Phi^{-1})\circ \Phi$, with $f \in C^\infty(V)$,
the pullback of $A$ by this diffeomorphism, we get
\[
B_{A_{\Phi}}(\Phi^{-1}(0),\rho)=\Phi^{-1}(B_A(0,\rho)), \quad \rho>0.
\]
Equivalently, if $\gamma=\gamma(t)$ is a subunit curve for $A$
such that $\gamma(0)=0$ and $\gamma(1)=x$, then $\eta(t):=\Phi^{-1}(\gamma(t))$
is a subunit curve for $A_\Phi$ and $\eta(0)=\Phi^{-1}(0)$, $\eta(1)=\Phi^{-1}(x)$.
Hence, essentially, subunit balls are sent to subunit balls through a diffeomorphism.
This remark, in particular, states that the notion of subunit balls on manifolds can also be equivalently given in coordinates.
\end{remark}
\begin{exa}
Let us consider $(\mathbb{R}^n,g_E)$, where $g_E$ is the standard Euclidean metric and $A=-\Delta$. Then, the subelliptic ball associated with $A$ is just the Euclidean ball, i.e.  
\[
B_A(x,\rho)=B_E(x,\rho)=\lbrace y \in \mathbb{R}^n; \ \abs{x-y}<\rho \rbrace,
\]
(or equivalently $d_\Delta(x,y)$ is just the ordinary Euclidean distance). \\
\end{exa}

\begin{exa}
If we consider on $(\mathbb{R}^2, g_E)$ the Grushin operator 
\[
A=\frac{\partial^2}{\partial x_1^2}+x_1^{2}\frac{\partial^2}{\partial x_2^2}
\]
by proceeding similarly it is easy to see that for $\rho>0$
\[
B_A(0,\rho)\approx [-\rho,\rho] \times [-\rho^2,\rho^2].
\]
Indeed, by definition $x=(x_1,x_2) \in B_A(0,\rho)$ if and only if there exists $\gamma:[0,1] \rightarrow \mathbb{R}^2$ subunit for $\rho^2A$, $\gamma(0)=0$ and $\gamma(1)=x$ such that (with $\gamma(t)=(\gamma_1(t),\gamma_2(t))$)
\begin{equation}\label{eq.Gru}
|\gamma_1(t)\xi_1+\gamma_2(t)\xi_2|^2 \leq \rho^2(\xi_1^2+\gamma_1(t)^2\xi_2^2), \quad \text{for a.e.} \ t \in [0,1]
\end{equation}
Then, if $x \in [-c\rho,c\rho] \times [-c\rho^2,c\rho^2]$, for $c$ sufficiently small we have that the curve $\gamma(t)=(t+(1-t)x_1,t+(1-t)x_2)$ satisfies \eqref{eq.Gru}. On the other hand, if $x \in B_A(0,\rho)$ there exists a curve $\gamma=\gamma(t)$ that satisfies \eqref{eq.Gru} and $\gamma(1)=x$. Hence $x \in [-C\rho,C\rho] \times [-C\rho^2,C\rho^2]$, for $C>0$ sufficiently large. 
\end{exa}
\section{\texorpdfstring{$h$-Subelliptic Random Walks}{h-Subelliptic Random Walks}}\label{sec.subwalk}
The aim of this section is to investigate different ways
of constructing an \textit{$h$-subelliptic random walk} on our manifold $M$, 
where $h \in ]0,h_0]$, $h_0 \ll 1$, can 
be thought of as the semiclassical parameter (or as the Planck constant).
\subsection{Preliminary Comments}\label{sec.precomm}
In the first place it is useful to specify what we mean by a Random Walk on a manifold, which we usually call a diffusion. To do so we refer to \cite{JJ} (Section 8.2).
\begin{definition}\label{def.ranwalk}
    Let $(M,g)$ be a closed manifold. A \textit{Random Walk} on $M$ is a family $P=(p_x)_{x \in M}$ of probability measures
    $p_x$ on $(M,\mathscr{B})$, depending measurably on $x$.
\end{definition}
Roughly speaking, for $B \in \mathscr{B}$ the value $p_x(B)$ is interpreted as the probability of the random walk to land on $B$
after one step starting from $x \in M$. If we consider a discrete stochastic process $\lbrace X_n \rbrace_{n \in \mathbb{N}_0}$ that satisfy the Markov property (see for instance \cite{JJ}, p. 233 or \cite{T}), by the Chapman-Kolmogorov equation, the first transition probability functions $p_1(x,dy)$ determine all the others. Hence, for a discrete Markov process we have that it is determined by the random walk $(p_1(x,\cdot))_{x \in M}$.
Therefore, when speaking about random walks one can refer interchangeably to a stochastic process $(X_n)_{n \in \mathbb{N}_0}$ or to the family
$(p_1(x,\cdot))_{x \in M}$.

The main novelty of the constructions that we will carry out is that they apply 
to general second-order  differential operators on $M$, i.e. operators
that in coordinates look like (with $(U,\varphi)$ local chart) 
\begin{equation}\label{opcoordA}
A_U=-\sum_{i,j=1}^n a_{ij}(\tilde{x})\frac{\partial^2}{\partial {\tilde{x}_i}\partial {\tilde{x}_j}}+\sum_{k=1}^nb(\tilde{x})\frac{\partial}{\partial \tilde{x}_k}+d(\tilde{x}), \quad \tilde{x} \in \varphi(U)\subseteq \mathbb{R}^n,
\end{equation}
where $a_{ij}, \, b_k, \, d$ are $C^\infty$ real-valued functions
for all $i,j,k \in \lbrace 1,...,n \rbrace$ and $\mathsf{A}_2=(a_{ij})_{i,j} \geq 0$.
\begin{remark}
In particular, we wish to point out that our operators
include second-order  differential
operators that cannot be written 
as sums of squares of smooth vector fields. Indeed, if $A$ is locally of the form  \eqref{opcoordA}, denoting as usual with $\mathsf{A}_2(\tilde{x})=(a_{ij}(\tilde{x}))_{i,j=1,\dots,n} \geq 0$ the matrix of the coefficients of the principal part, it is well known (see for instance \cite{JSC}, p. 48) that the matrix $\mathsf{A}_2$ can always be factorized - though not uniquely - as $\mathsf{A}_2=B\transp{B}$ where $B(\tilde{x})=(b_{ij}(\tilde{x}))_{i,j}$ is an $n \times m$-matrix. However, in general, the coefficients $b_{ij}$ are not smooth functions. Therefore, although we may write 
\[
A_U=-\sum_{j=1}^mX_j^\ast X_j=-\sum_{j=1}^mX_j^2+X_0,
\]
where $X_j=\sum_{i=1}^n b_{ij}\partial/\partial \tilde{x}_i$, the vector fields $X_j$ are \emph{not smooth}. 
Other interesting examples of second-order differential operators that cannot be written as sums of squares of vector fields can be obtained from local operators of the form $\sum_{j=1}^n Y_j^\ast p(x)Y_j$,
where $p(x)=\sum_{\abs{\alpha}=k}c_\alpha x^{2\alpha}$, with $k \in \mathbb{N}$ and $c_\alpha>0$ for each $\abs{\alpha}=k$ (cf. \cite{PK}).
\end{remark}

In this section, we begin by dealing with cubes in $\mathbb{R}^n$,
then, to avoid any confusion, we start by fixing the notation.
\begin{definition}
Let $x_0 \in \mathbb{R}^n$. \textit{The cube of radius $r>0$ centered 
at $x_0$} in $\mathbb{R}^n$ is defined as the open set
\[
Q(x_0,r)=\bigl\lbrace x \in \mathbb{R}^n; \ \abs{x-x_0}_\infty <r\bigr \rbrace,
\] 
where 
\[
\abs{x-x_0}_\infty:=\max_{j=1,\dots,n} \abs{x_j-x_{0,j}}.
\]
With this notation, the cube $Q(0,1)$ of radius $1$ centered at $0$
in $\mathbb{R}^n$
is called \textit{the unit cube}.
In what follows we denote it by $Q=Q(0,1)$. Moreover, for $c^{\ast\ast}>c^\ast>1$ 
we define the \textit{dilations} of $Q$, with constants $c^\ast$ and $c^{\ast\ast}$, respectively, as
\[
Q^\ast=\bigl\lbrace y \in \mathbb{R}^n; \ \exists \, x \in Q \; \text{s.t.} \; y=c^\ast x \bigr\rbrace, 
\]
\[
Q^{\ast\ast}=\bigl\lbrace y \in \mathbb{R}^n; \ \exists \, x \in Q \; \text{s.t.} \; y=c^{\ast \ast} x \bigr\rbrace.
\]
\end{definition}

Our first step is to construct a local diffusion using
the so-called Oleinik-Radkevich's vector fields 
(see for instance  \cite{OR} or
 \cite{JSC}),
where the term "local" comes from the fact
that we will work on $Q^\ast$, fixed dilation of the unit cube $Q$.

Nevertheless, as we will see in Remark \ref{rmk.notopt}
this construction is not optimal.
Then, we will construct a global "$(h,\rho)$-subelliptic random walk" using
some local results of Fefferman
and Phong, based on the Calder\' on-Zygmund localization
(see again \cite{FP1}, \cite{FP2}, \cite{JSC}),
which will lead to a global construction
due to the compactness of $M$.
\subsection{The Oleinik and Radkevich Local Diffusion}
This subsection is devoted to constructing
the \textit{Oleinik-Radkevich local diffusion}.
The term local, once more, comes from the fact that
we will work in an open cube of $\mathbb{R}^n$.
In this context we we work with $A$ subelliptic formally self-adjoint
operator of the form stated above, i.e.
\[
A=-\sum_{i,j=1}^n a_{ij}(x)\frac{\partial^2}{\partial x_i\partial x_j}+\sum_{k=1}^nb_k(x)\frac{\partial}{\partial x_k}+d(x),
\]
where $a_{ij}, \, b_k, \, d \in C^\infty(Q^\ast,\mathbb{R})$
for all $i,j,k \in \lbrace 1,...,n \rbrace$,
such that the matrix $(a_{ij})_{i,j}$ is positive semidefinite.
Note that the coefficients of the operator are bounded on $Q$ 
since they are defined on $Q^\ast\supseteq \bar{Q}$. 

The Oleinik-Radkevich local diffusion
is based on the
\textit{Oleinik-Radkevich
vector fields} defined as follows.

\begin{definition}\label{ORvec}
\textit{The $j$-th Oleinik-Radkevich's vector field on $Q^\ast$}, where $j\in \lbrace 1,\dots,n \rbrace$, is the vector field $Y_j \in \Gamma(Q^\ast,TQ^\ast)$ defined by
\[
Y_j(x):=\sum_{i=1}^na_{ij}(x) \frac{\partial}{\partial x_i}, \quad x \in Q^\ast.
\]
\end{definition}
For these vector fields we have the following theorem (see \cite{JSC}, Proposition 2.1$'$, p. 66).
\begin{theorem}
If $A$ is subelliptic of order $\varepsilon$, there exist $c=c(\varepsilon)>0$ and $N=N(\varepsilon) \in \mathbb{N}$ universal constants, such that 
\begin{equation}\label{cond21'}
\sum_{\abs*{J}\leq N}\abs*{Y_{J}(0,\xi)} \geq c \abs*{\xi}, \quad \forall \xi \in \mathbb{R}^n,
\end{equation}
where $Y_J=[Y_{j_1},...[Y_{j_{s-1}},Y_{j_s}]...]$, for $J=(j_1,...,j_s)$, s $\in \mathbb{N}$.
\end{theorem}
\begin{remark}\label{quanthor}
This theorem is a quantitative version of H\"ormander's 
condition at the origin. Indeed, if \eqref{cond21'} holds
the family of vector fields $\lbrace Y_J \rbrace_{\abs{J}\leq N}$
contains $n$ vector fields linearly indipendents at the origin.
However, in general, there is no known simple relationship between such $N$
and the order of subellipticity $\varepsilon$.   
\end{remark}
\begin{exa}
Let $k \in \mathbb{N}$ and let $A$ be the $k$-Grushin operator
on $\mathbb{R}^2$ defined by
\[
A=\frac{\partial^2}{\partial x_1^2}+x_1^{2k}\frac{\partial^2}{\partial x_2^2}.
\]
We have that 
\[
\mathsf{A}_2(x)=(a_{ij}(x))_{i,j}=
\begin{pmatrix}
1 & 0 \\
0 & x_1^{2k}
\end{pmatrix},
\]
and so, in this case,
\[
Y_1=\frac{\partial}{\partial x_1}, \quad \quad Y_2=x_1^{2k}\frac{\partial}{\partial x_2}.
\]
\end{exa}

We next show that the Oleinik-Radkevich vector fields 
are subunit and satisfy H\"ormander's condition 
(possibly after a rescaling). Therefore, 
we will use such a property to construct
a local subelliptic diffusion.

\begin{theorem}
Let $Q, Q^\ast$ be as before,
let $A$ be a subelliptic second-order formally self-adjoint
differential operator defined 
on $Q^\ast$ with the aforementioned properties
and let $Y_1,\dots,Y_n$ be the associated Oleinik-Radkevich vector fields.
Then, there exists $C>0$ such that, the vector fields $(\tilde{Y}_j)_j$ defined as
\[
\tilde{Y}_j:=\frac{1}{\sqrt{C}}Y_j
\]
satisfy H\"ormander's condition and they are subunit
vector fields on this cube, i.e. for all $j \in \lbrace 1,\dots n \rbrace$, one has
\[
\abs*{\langle \tilde{Y}_j(x), \xi \rangle}^2 \leq \sum_{i,j=1}^n a_{ij}(x)\xi_i \xi_j, \quad  \forall (x,\xi) \in Q \times \mathbb{R}^n.
\]
\end{theorem}
\begin{proof}
We first note that, since $A$ is subelliptic, by Remark \ref{quanthor}, there exists $r \in \mathbb{N}$, locally independent of $x$, such that $Y_1,...,Y_n$ and its commutators
of length at most $r$, span the tangent space for every $x \in Q$.
To prove the subunit property it is sufficient to observe that
if we consider $A=(a_{ij}(x))_{i,j} \geq 0$,
the matrix of coefficients of $A$, there exists a constant $C>0$ such that 
\[
2a_{jj}(x) \leq C, \quad \forall x \in Q,
\]
for all $j=1,\dots n$, since they are defined on $Q^\ast$
that contains $Q$.
Hence, in view of Proposition \ref{propmatrixsub}, setting
\[
\tilde{Y}_j:=\frac{1}{\sqrt{C}}Y_j,
\]
for all $j \in \lbrace 1,..,n \rbrace$, we have that they satisfy 
\[
\begin{split}
\abs*{\langle \tilde{Y}_j(x), \xi \rangle}^2&=\frac{1}{C}\Bigl(\sum_{i=1}^n a_{ij}(x)\xi_i\Bigr)^2 \\
&\leq 2\frac{a_{jj}(x)}{C} \sum_{i,j=1}^n a_{ij}(x)\xi_i \xi_j \\
&\leq \sum_{i,j=1}^n a_{ij}(x)\xi_i\xi_j,
\end{split}
\]
for all $x \in Q$ and for all $\xi \in \mathbb{R}^n$.
Thus they are subunit vector fields.
\end{proof}
\noindent Now, for $h \in ]0,h_0]$, with $h_0>0$ sufficiently small
we construct a local $h$-Markov process through this theorem as follows.\\
$\bullet$ Fix $x_0 \in Q$;\\
$\bullet$ at step $k$ we choose at random $j \in \lbrace 1,...,n \rbrace$ 
and $t \in ]-h,h[$ (with respect the uniform probability) and set 
\[
x_{k+1}:=e^{t\tilde{Y}_j}x_k.
\]
If $ h \in ]0,h_0]$, with $h_0>0$ sufficiently small,
and $x \in Q$, we define 
\[
\mathsf{Y}_{j,h}f(x):=\frac{1}{2h}\int_{-h}^h f(e^{t\tilde{Y}_j}x)dt, \quad f \in C^\infty(Q^\ast;\mathbb{R})
\]
and the $h$-Markov operator associated with our local $h$-Markov process 
\[
\mathsf{Y}_hf:= \frac{1}{n}\sum_{j=1}^n\mathsf{Y}_{j,h}f, \quad f \in C^\infty(Q^\ast;\mathbb{R}).
\]
\subsection{CZ-localization}
In the previous subsection using Oleinik-Radkevich's
vector fields we have constructed a local subelliptic random walk.
However, as we will see in the next example,
such a construction is not optimal. 
\begin{remark}\label{rmk.notopt}
Let $A$ be the $k$-Grushin operator on $\mathbb{R}^2$ defined by 
\[
A=\frac{\partial^2}{\partial x_1^2}+x_1^{2k}\frac{\partial^2}{\partial x_2^2}, \quad k \in \mathbb{N}.
\]
Note that $A$ can be written as 
\[
A=X_1^2+X_2^2,
\]
with $X_1=\frac{\partial}{\partial x_1}$,
$X_2=x_1^k\frac{\partial}{\partial {x_2}}$
and that they satisfy H\"ormander condition at order $k+1$ 
(i.e. $X_1,X_2$ and their commutators 
up to order $k+1$ span the tangent space
at every poinx $x \in \mathbb{R}^n$).
On the other hand, in this case,
the Grushin vector fields are defined by
$Y_1=\frac{\partial}{\partial x_1}$,
$Y_2=x_1^{2k}\frac{\partial}{\partial x_2}$
and they satisfy H\"ormander condition at order 
$2k+1$.
Roughly speaking, this example shows that the speed of convergence
is not optimal (since we need more commutators).
Then our purpose is to improve it
using the tools introduced by Fefferman and Phong in \cite{FP2}
based on the so-called Calder\' on-Zygmund localization.
\end{remark}
In order to give the fundamental tool,
let us first recall the Calder\' on-Zygmund localization in details.
Let us consider again $Q$ be the unit cube and let $Q^\ast, Q^{\ast \ast}$ be
two fixed dilations of constants $c^\ast$ and $c^{\ast \ast}$ respectively,
with $1<c^\ast <c^{\ast \ast}$.
We also always suppose that $A$ is defined 
on $Q^{\ast\ast}$,
in order to ensure that the coefficients are bounded in $Q^\ast$. 
Finally, recall that, for $U,V\subset \mathbb{R}^n$ open sets,
$\Phi:U\rightarrow V$ diffeomorphism and $A$ differential operator defined on $V$,
we denote the pullback of $A$ between $\Phi$ defined on $U$ 
by $A_\Phi f=A(f\circ \Phi^{-1})\circ \Phi$, $f \in C^\infty(U)$.

The main result of this subsection is the following theorem,
that essentially allow us to straighten the first vector field 
by using the Calder\' on-Zygmund localization (cf. \cite{F} 
or \cite{JSC}).
\begin{theorem}\label{CZtheo}
Let $A$ be a
second-order differential operator defined on $Q^{\ast\ast}$ of the form 
\[
A=-\sum_{i,j=1}^n a_{ij}(x)\frac{\partial^2}{\partial {x_i}\partial {x_j}}+\sum_{k=1}^nb_k(x)\frac{\partial}{\partial x_k}+d(x),
\]
where $a_{ij},b_k,d \in C^\infty(Q^{\ast\ast},\mathbb{R})$
for every $i,j,k=1,...,n$ and $\mathsf{A}_2(x)=(a_{ij}(x))_{i,j}$
is symmetric and positive semidefinite for any $x \in Q^{\ast\ast}$.\\
Then, there exists a change of coordinates, diffeomorphism on the image
\[
\begin{split}
\phi:]-r_1,r_1[ \times  Q_{n-1}&(0,1) \rightarrow \phi(]-r_1,r_1[ \times  Q_{n-1})\subseteq Q(0,\delta_1), \quad \phi(0)=0,\\
&(w_1,z) \mapsto x,
\end{split} 
\]
for some $\delta_1=2^{-k}$, $k \in \mathbb{N}$, and $r_1>0$ constant
that depends only on the bounds of the $a_{ij}$'s on $Q^\ast$.\\
Moreover, the pullback of $A$ through $\phi$ has the form 
\[
A_{\phi}=-\frac{\partial^2}{\partial w_1^2}-\sum_{i,j=2}^nb_{ij}(w_1,z) \frac{\partial^2}{\partial z_i \partial z_j} \ + \ \text{l.o.t}
\]
where $b_{ij} \in C^\infty(]-r_1,r_1[ \times  Q_{n-1};\mathbb{R})$ for all $i,j=2,\dots,n$
(and l.o.t stands for 'lower order terms'). 
\end{theorem}
\begin{proof}
To prove the theorem we need the following Calder\' on-Zygmund
localization argument (see \cite{FSC} or \cite{JSC}). 
Choose $\delta_1=2^{-N_1}$ for the
smallest $N_1\in \mathbb{N}$, such that 
\begin{equation}
\label{Cz1}
\max_{1 \leq i \leq n}\max_{x \in \overline{Q(0,\delta_1)}}\abs*{a_{ii}(x)} \geq 10R_1\delta_1^2,
\end{equation}
where $R_1>0$ is a fixed constant so that
\[
\abs{\partial_x^\alpha a_{ij}(x)}\leq R_1
\]
for all $x\in Q^\ast$, for all $\alpha \in \mathbb{N}^n$
such that $\abs{\alpha}\leq 2$ and for all $i,j=1,\dots,n$.
This choice of $R_1$ will be more clear soon.
We can assume, without loss of generality,
that the index for which the maximum is attained is $i=1$.

\noindent The key point of the argument is the following lemma (for the proof see \cite{F}, p. 183),
which provides estimates that will lead to the
universal bounds (see Lemma \ref{universalbound} below) 
for the coefficients of the principal part.

\begin{lemma}\label{lemR1}
Let $\mathsf{A}_2(x)=(a_{ij}(x))_{ij}$ and $R_1>0$ be as before,
then for all $x \in Q(0,\delta_1)$
\[
\begin{split}
&(i) \ \abs*{a_{ij}(x)} \leq 40R_1\delta_1^2, \quad \forall i,j \in \lbrace 1, \dots, n\rbrace, \\
&(ii)\ a_{11}(x) \geq R_1\delta_1^2.
\end{split}
\]
\end{lemma}
We now introduce a first change of coordinates 
\[
\begin{split}
\psi:Q(0,10)& \rightarrow Q(0,\delta_1), \\
y & \mapsto x=\frac{\delta_1 y}{10}.
\end{split}
\]
Hence, if we consider the pullback of $A$ by this diffeomorphism,
we have in these new coordinates
\[
A_{\psi}=-\sum_{i,j=1}^n \tilde{a}_{ij}(y)\frac{\partial^2}{\partial y_i \partial y_j}+ \ \text{l.o.t.},
\]
where $\tilde{a}_{ij}(y)=100\delta_1^{-2}a_{ij}(\frac{\delta_1 y}{10})$. 
Note now that these coefficients 
have bounds that do not depend on $\delta_1$.
More precisely, we have the next lemma (see \cite{F}, p. 184).

\begin{lemma}\label{universalbound}
Let $\psi:Q(0,10)\rightarrow Q(0,\delta_1)$ be the
previous diffeomorphism and let $A_\psi f:=A(f\circ \psi^{-1}) \circ \psi $,
for $f \in C^\infty(Q(0,10))$,
be the pullback of $A$ by it.
Then, the coefficients of the operator 
\[
A_{\psi}=-\sum_{i,j=1}^n \tilde{a}_{ij}(y)\frac{\partial^2}{\partial y_i \partial y_j}+ \ \text{l.o.t.},
\]
have bounds that do not depend on $\delta_1>0$,
i.e. there exists a universal constant $C>0$,
independent of $\delta_1$, such that 
\[
\abs{\partial_y^{\alpha}\tilde{a}_{ij}(y)}\leq C,
\]
for all $y \in Q(0,10)$, for all $\alpha \in \mathbb{N}^n$ and for all $i,j=1,\dots,n$. 
\end{lemma}
We now come back to the proof of our theorem and write  
\[
\begin{split}
A_{\psi}&=-\sum_{i,j=1}^n \tilde{a}_{ij}(y)\frac{\partial^2}{\partial y_i \partial y_j}+ \ \text{l.o.t.} \\
&=-\tilde{a}_{11}(y)\Bigl( \Bigl(\frac{\partial}{\partial y_1}+\sum_{j=2}^n\frac{\tilde{a}_{1j}(y)}{\tilde{a}_{11}(y)}\frac{\partial}{\partial y_j}\Bigr)^2+\sum_{i,j=2}^n \tilde{a}'_{ij}(y)\frac{\partial^2}{\partial y_i \partial y_j}\Bigr)+ \ \text{l.o.t.},
\end{split}
\]
where, for all $i,j=2,\dots,n$, and for $y \in Q(0,10)$, 
\[
\tilde{a}'_{ij}(y):=\Bigl(\frac{\tilde{a}_{ij}(y)}{\tilde{a}_{11}(y)}-\frac{\tilde{a}_{1i}(y)\tilde{a}_{1j}(y)}{(\tilde{a}_{11}(y))^2}\Bigr).
\]
Hence, setting $Y:= \Bigl(\frac{\partial}{\partial y_1}+\sum_{j=2}^n\frac{\tilde{a}_{1j}(y)}{\tilde{a}_{11}(y)}\frac{\partial}{\partial y_j}\Bigr)$, we have
\[
A_\psi=-\tilde{a}_{11}(y)\Bigl(Y^2+\sum_{i,j=2}^n \tilde{a}'_{ij}(y)\frac{\partial^2}{\partial y_i \partial y_j}\Bigr)+ \ \text{l.o.t.}.
\]
Next, we change variables $y=\tilde{\varphi}(\tilde{w}_1,\tilde{z})$,
where $\tilde{\varphi}(\cdot,\tilde{z})$ 
is the integral curve of the vector field $Y$ starting at $(0,\tilde{z})$. 
In these new coordinates, if we set $\tilde{\phi}:=\psi \circ \tilde{\varphi}$,
we have 
\[
A_{\tilde{\phi}_1}=-\tilde{a}_{11}(\tilde{\varphi}(\tilde{w}_1,\tilde{z}))\Bigl(\frac{\partial^2}{\partial \tilde{w}_1^2}+\sum_{i,j=2}^{n}\tilde{b}_{ij}(\tilde{w}_1,\tilde{z})\frac{\partial^2}{\partial \tilde{z}_i\partial \tilde{z}_j}\Bigr)+ \ \text{l.o.t.},
\]
for some smooth coefficients $\tilde{b}_{ij}$, for $i,j=2,\dots,n$.\\
Now, if $\abs*{\tilde{w}_1}<1$, $\abs*{\tilde{z}}<1$ we have $\abs*{\tilde{\varphi}_1(\tilde{w}_1,\tilde{z})}=\abs*{\tilde{w}_1}<1$. 
Moreover, by Taylor's expansion, we may write  
\[
\tilde{\varphi}_j(\tilde{w}_1,\tilde{z})=\tilde{\varphi}_j(0,\tilde{z})+\tilde{w}_1\dot{\tilde{\varphi}}_j(\tilde{w}_1',z), \quad \tilde{w}_1' \in (0,\tilde{w}_1).
\]
Hence, since $\abs*{\tilde{a}_{1j}(y)}\leq  \tilde{a}_{11}(y)$ (by $\abs*{\tilde{a}_{1j}(y)} \leq \tilde{a}_{jj}(y)\tilde{a}_{11}(y)\leq \tilde{a}_{11}(y)^2$), one has
\[
\abs*{\tilde{\varphi}_j(\tilde{w}_1,\tilde{z})}< \abs*{\tilde{z}}+\abs*{\tilde{w}_1\dot{\tilde{\varphi}}_j(\tilde{w}_1',\tilde{z})}<1+8\abs*{\tilde{w}_1}<10, \quad \forall j \in \lbrace 2,...,n \rbrace.
\]
Thus, we have obtained that $\tilde{\varphi}(\tilde{w}_1,\tilde{z}) \in Q(0,10)$
if $\abs*{\tilde{w_1}}<1$, $\abs*{\tilde{z}}<1$, and
\[
\begin{split}
\tilde{\phi}=\psi \circ \tilde{\varphi} : Q(0,1)& \xrightarrow{\tilde{\varphi}} \tilde{\varphi}(Q(0,1))\subseteq Q(0,10) \xrightarrow{\psi} \tilde{\phi}(Q(0,1))\subseteq Q(0,\delta_1), \\
(\tilde{w}_1,\tilde{z})&\mapsto y=\tilde{\varphi}(\tilde{w}_1,\tilde{z}) \mapsto x=\frac{\delta_1 y}{10}.
\end{split}
\]
Finally, since $\tilde{a}_{11}$ is bounded
away from $0$ and from above by
constants independent of $\delta_1$,
if we consider the change of variables $(w_1,z)=\kappa(\tilde{w}_1,\tilde{z}_1)$,
using as before the flow $\kappa(\cdot,\tilde{z})$ starting form $(0,\tilde{z})$
associated with the vector field $\sqrt{\tilde{a}_{11}}\frac{\partial}{\partial \tilde{w}_1}$,
we obtain the last diffeomorphism with range into $Q(0,1)$ 
\[
\begin{split}
\kappa: ]-r_1,r_1[ \times  Q&_{n-1}(0,1)  \rightarrow \kappa(]-r_1,r_1[ \times Q_{n-1}(0,1)) \subseteq Q(0,1), \\
& (w_1,z) \mapsto (\tilde{w}_1,\tilde{z}),
\end{split}
\]
where $r_1>0$ depends only on the bounds of the
coefficients and not on $\delta_1$. 
Then, if we call $\varphi:=\tilde{\varphi} \circ \kappa$ and $\phi:=\psi \circ \varphi$,
we have that
\[
\begin{split}
\phi : ]-r_1,r_1[ \times Q_{n-1}(0,1)& \rightarrow \phi(]-r_1,r_1[ \times Q_{n-1}(0,1))\subseteq  Q(0,\delta_1), \\
(w_1,z)& \mapsto x=\frac{\delta_1 \varphi(w_1,z)}{10}
\end{split}
\]
is a diffeomorphism onto its image, $\phi(0)=0$ and $A_{\phi}$, 
the pullback of $A$ through $\phi$, satisfies  
\[
A_{\phi}=-\frac{\partial^2}{\partial w_1^2}+\sum_{i,j=2}^{n}b_{ij}(w_1,z)\frac{\partial^2}{\partial z_i\partial z_j}+ \ \text{l.o.t.},
\]
for some $b_{ij} \in C^\infty(]-r_1,r_1[ \times Q_{n-1}(0,1);\mathbb{R})$,
for all $i,j=2,\dots,n$. 
\end{proof}
Now, in order to iterate this argument on the dimension we make
the following observation.
\begin{remark}\label{rmkbij}
Note that the matrix 
\[
\tilde{\mathsf{A}}'_2(y)=(\tilde{a}'_{ij}(y))_{i,j=2,..,n}:=\Bigl(\frac{\tilde{a}_{ij}(y)}{\tilde{a}_{11}(y)}-\frac{\tilde{a}_{1i}(y)\tilde{a}_{1j}(y)}{(\tilde{a}_{11}(y))^2}\Bigr)_{i,j=2,...,n}, \quad y \in Q(0,10),
\]
satisfies $\tilde{a}'_{ii}(y) \geq 0$ for all $i=2,\dots n$ 
(since $\abs{\tilde{a}_{ij}(y)}^2\leq \tilde{a}_{ii}(y)\tilde{a}_{jj}(y)$). 
Then, setting
\[
\bar{A}_{\phi}=\sum_{i,j=2}^n\bar{b}_{ij}(z)\frac{\partial^2}{\partial z_j \partial z_k}, \quad \bar{b}_{ij}(z)= \frac{1}{2}\int_{-1}^{1}b_{ij}(w_1,z)dw_1,
\]
(with $b_{ij}'s$ as in the proof of the theorem) 
by induction we may apply the Calder\' on-Zygmund
localization for the operator $\bar{A}_{\phi}$ in $(n-1)$-variables.
\end{remark}
\subsection{The Subunit Geometry}
We are now ready to approach the problem of constructing
a global subelliptic diffusion on $M$.
To do that we need to understand the geometry
of the subunit balls (see \cite{FP2}). 
The goal of this subsection is to prove Theorem \ref{fundtheorem}
in which we will see some geometric consequences of the subellipticity.
This theorem, in some sense, is an alternative version of Theorem 1 in \cite{FP2}.

The goal is to locally reduce matter to a universal cube
where the operator $A$ is "equivalent" to a new second-order
differential operator with constant coefficients $\tilde{A}$.
This will be the most important geometric property that we will use
to construct a global diffusion.
The term "equivalent" will be more clear in the following.
We first do this locally on a dilation of the unit cube 
$Q^\ast \subseteq \mathbb{R}^n$
where the dilation constant is chosen such that 
for all $x_0 \in Q$, we have $Q(x_0,1)\subseteq Q^\ast$,
where $Q(x_0,1)$ is the unit cube centered at $x_0$. 
We also assume as before that $A$ is defined on
$Q^{\ast\ast}$ fixed dilation of $Q^\ast$.

The idea is to consider for every point $x_0 \in Q$ 
the Calder\' on-Zygmund cube centered at $x_0$
and reduce the diffusion on a neighborhood 
of this point to the diffusion on \textit{the universal cube} 
(see Definition \ref{univcube}), which is our universal model
and where we consider the diffusion associated
with the operator $\tilde{A}$.

\begin{proof}[Proof of Theorem \ref{fundtheorem}]
Let us fix $x_0\in Q$. The first aim is to construct a diffeomorphism $\Phi$ between a neighborhood of $0$ and a neighborhood of $x_0$. Possibly by composing with a translation,
we may assume that $x_0=0$ and work on a neighborhood of $0$.

Then we use the Calder\' on-Zygmund localization centered at $0$. Choose $\delta_1=2^{-N_1}$ for the first $N_1 \in \mathbb{N}$ such that 
\[
\max_{1 \leq i \leq n}\max_{x \in \overline{Q(0,\delta_1)}}\abs*{a_{ii}(x)} \geq 10 R_1\delta_1^2,
\]
where $R_1>0$ is a constant much larger than
a bound for the $a_{ij}$'s and their derivatives
on $Q^\ast$ determined on Lemma \ref{lemR1}.
Without loss of generality, as before, we can also
assume that the index for which the maximum is attained is $i=1$.
By Theorem \ref{CZtheo}, there exists a change of coordinates,
diffeomorphism onto its image
\[
\begin{split}
\phi_1:]-r_1,r_1[ \times  Q_{n-1}&(0,1) \rightarrow \phi_1(]-r_1,r_1[ \times  Q_{n-1}(0,1)) \subseteq Q(0,\delta_1), \\
&(w_1,z) \mapsto x,
\end{split}
\]
for some $\delta_1=2^{-N_1}$, $N_1 \in \mathbb{N}$, $r_1>0$ 
constants that depends only on the bounds
of the $a_{ij}$'s on $Q^\ast$, such that $\phi_1(0)=0$ 
and, moreover, the pullback of $A$ by $\phi_1$ has
the form 
\begin{equation}
\label{CZA1}
A_{\phi_1}=-\frac{\partial^2}{\partial w_1^2}-\sum_{i,j=2}^nb_{ij}(w_1,z)\frac{\partial^2}{\partial z_i \partial z_j}+ \ \text{l.o.t.},
\end{equation}
with the $b_{ij}$'s smooth with bounds independent of $\delta_1$. 

We now want to apply the following fundamental Lemma
due to Fefferman and Phong (cf. \cite{FP2}, Lemma 2 and the proof of Lemma 1) 
to this operator.
\begin{lemma}\label{FPlemma}
Let $L$ be a subelliptic operator of order $0<\varepsilon<1$
defined on $\mathbb{R}^{n}$ with the properties discussed above
(see \ref{hp1}, \ref{hp2} and \ref{hp3}) of the form 
\[
L=-\frac{\partial^2}{\partial t^2}-\sum_{i,j=2}^{n}a_{ij}(t,x)\frac{\partial^2}{\partial x_i\partial x_j}
\]
and let $\lambda>0, \, N>\varepsilon^{-50}$. Define
\[
L_\lambda=-\frac{\partial^2}{\partial t^2}-\sum_{i,j=2}^{n}a_{ij}(t,x)\frac{\partial^2}{\partial x_i\partial x_j}-\lambda^{2N}\Delta_x,
\]
and
\[
\bar{L}_\lambda=-\sum_{i,j=2}^n\bar{a}_{ij}(x)\frac{\partial^2}{\partial x_i\partial x_j}-\lambda^{2N}\Delta_x,
\]
where
\[
\bar{a}_{ij}(x)=\frac{1}{2\lambda}\int_{\abs*{t}\leq \lambda}a_{ij}(t,x)dt.
\]
Then, there exist universal constants $c',C'>0$ such that 
\[
B_{L_\lambda}(0,c'\lambda) \subseteq \lbrace \abs*{t}\leq \lambda \rbrace \times B_{\bar{L}_\lambda}(0,\lambda) \subseteq B_{L_\lambda}(0,C'\lambda). 
\]
Moreover, the subelliptic balls of the operators
$L$ and $L_{\lambda}$ are equivalent in the sense 
\[
B_{L_\lambda}(0,c''\lambda) \subseteq B_{L}(0,\lambda) \subseteq B_{L_\lambda}(0,C''\lambda), 
\]
for some universal constants $c'',C''>0$.
\end{lemma}
In order to use this lemma in our case we consider for $\rho>0$ the operator  
\[ 
\rho^2A_{\phi_1}=-\rho^2\frac{\partial^2}{\partial w_1^2}-\rho^2\sum_{i,j=2}^nb_{ij}(w_1,z)\frac{\partial^2}{\partial z_i \partial z_j}+ \ \text{l.o.t.}.
\]
If we now choose a universal $\rho>0$ and sufficiently small 
(recall that $w_1 \in ]r_1,r_1[$) we may consider another change 
of variables $\psi^{(1)}_\rho$ defined by   
\begin{equation}\label{wt}
Q(0,2) \ni (t_1,t') \xmapsto{\psi^{(1)}_\rho}  (\rho t_1,\rho t')=(w_1,z) \in ]-r_1,r_1[ \times Q_{n-1}(0,1).
\end{equation}
Thus, the pullback of the operator $\rho^2A_{\phi_1}$ through $\psi^{(1)}_\rho$ is $L=(\rho^2A_{\phi_1})_{\psi^{(1)}_\rho}$ defined by
\[
L=-\frac{\partial^2}{\partial t_1^2}-\sum_{i,j=2}^nb_{ij}(\rho t_1,\rho t')\frac{\partial^2}{\partial t_i' \partial t_j'}+  \ \text{l.o.t.}.
\] 
Therefore, by applying Lemma \ref{FPlemma} to the operator $L$ (with $\lambda=1$) we have
\begin{equation}\label{inclusiont}
B_{L}(0,c) \subseteq \lbrace \abs*{t_1}\leq 1  \rbrace \times B_{\bar{L}}(0,\tilde{c}) \subseteq B_{L}(0,C),
\end{equation}
that in our coordinates $(w_1,z)$ becomes
\begin{equation}\label{inclusionw}
B_{\rho^2A_{\phi_1}}(0,c) \subseteq \lbrace \abs*{w_1}\leq \rho  \rbrace \times B_{\rho^2\bar{A}_{\phi_1}}(0,\tilde{c}) \subseteq B_{\rho^2A_{\phi_1}}(0,C),
\end{equation}
for some universal constants $c,C,\tilde{c}>0$ and some universal $\rho=\rho(r_1)>0$
(recall that $r_1$ is universal as well) sufficiently small and where 
\begin{equation}\label{equivsubcubeballpreliminary}
\bar{A}_{\phi_1}=-\sum_{i,j=2}^n \bar{b}_{ij}(z)\frac{\partial^2}{\partial z_i \partial z_j}, \quad \bar{b}_{ij}(z)=\frac{1}{2\rho}\int_{-\rho}^{\rho}b_{ij}(w_1,z)dw_1.
\end{equation}
Finally, note that, since $B_{\rho^2A}(0,c')=B_A(0,c'\rho)$ for all $c'>0$, we may rewrite
\eqref{equivsubcubeballpreliminary} as
\begin{equation}\label{equivsubcubeball}
B_{A_{\phi_1}}(0,c\rho) \subseteq \lbrace \abs*{w_1}\leq \rho  \rbrace \times B_{\bar{A}_{\phi_1}}(0,\tilde{c}\rho) \subseteq B_{A_{\phi_1}}(0,C\rho).
\end{equation}
We now look at the operator 
\[
\bar{A}_{\phi_1}=-\sum_{i,j=2}^n \bar{b}_{ij}(z)\frac{\partial^2}{\partial z_i \partial z_j}, \quad \bar{b}_{ij}(z)=\frac{1}{2\rho}\int_{-\rho}^{\rho}b_{ij}(w_1,z)dw_1.
\]
By construction it is defined on $Q_{n-1}(0,1)$ 
we can repeat Calderón-Zygmund argument in this $(n-1)$-variables (see Remark \ref{rmkbij}),
i.e. we choose the first strictly positive integer $N_2$,
such that $\delta_2=2^{-N_2}$ satisfies
\[
\max_{2 \leq i \leq n}\max_{z \in \overline{Q(0,\delta_2)}}\abs*{\bar{b}_{ii}(z)} \geq 10 R_2\delta_2^2,
\]
with $R_2$ positive universal constant sufficiently large.
Thus, following the previous discussion also in this case,
there exists a change of coordinates 
\[
\begin{split}
\phi_2:]-r_2,& r_2[ \times  Q_{n-2}(0,1) \rightarrow \phi_2(]-r_2,r_2[ \times  Q_{n-2}(0,1)) \subseteq Q(0,\delta_2), \\
&u'=(w_2,u'') \mapsto z,
\end{split}
\]
with $r_2>0$ universal. 
Hence, if we look at the operator $\bar{A}_{\phi_1}$
in these new variables, one has
\[
\bar{A}_{\phi_1,\phi_2}=-\frac{\partial^2}{\partial w_2^2}+\sum_{i,j=3}^nc_{ij}(w_2,u')\frac{\partial^2}{\partial u'_i \partial u'_j},
\]
for some smooth real coefficients $c_{ij}$, with $i,j=3,\dots,n$.
Furthermore, also in this case,
we may introduce the last change of variables
\[
\psi^{(2)}_\rho:(t_2,t'')\mapsto (\rho w_2,\rho u''),
\]
for $\rho=\rho(r_2)$ universal sufficiently small 
and use again Lemma \ref{FPlemma} for the operator
in these new variables to obtain the corresponding
inclusion \eqref{inclusiont} in the "$t$-variables".
Moreover, since the operator $\bar{A}_{\phi_1}$ 
satisfies the same microlocal estimates of the operator $A_{\phi_1}$ 
(cf. \cite{FP2})
the inclusion \eqref{inclusionw} in "$w$-variables"
(the meaning of $t$-variables and $w$-variables will be more clear in the following)
reads as
\begin{equation}
\label{FP1}
B_{\bar{A}_{\phi_1,\phi_2}}(0,c'\rho) \subseteq \lbrace \abs*{w_2}\leq c_2 \rho^{\kappa_2} \rbrace \times B_{{\bar{\bar{A}}}_{\phi_1,\phi_2}}(0,\tilde{c}'\rho^{\kappa_2}) \subseteq B_{\bar{A}_{\phi_1,\phi_2}}(0,\tilde{c}\rho ),
\end{equation}
for some $c_2,c',\tilde{c}'>0$
universal constant and $1\leq \kappa_2 \leq 1/\varepsilon$.

\noindent Hence, we obtained a new diffeomorphism
onto its image (recall $u'=(w_2,u'')$)
\[
(w_1,w_2,u'')\xmapsto{(I_{w_1},\phi_2)}(w_1,\phi_2(u'))\xmapsto{\phi_1} x \in Q(0,\delta_1).
\]
Iterating $n$ times and composing
all the changes of variables
we obtain a diffeomorphism $\Phi^0$ (corresponding to the "$t$-variables")
between the standard unit cube $Q=Q(0,1)$ and its image,
i.e setting $Q_A(0,\rho):=\Phi^0(Q)$,
\[
Q \ni (t_1,\dots,t_n) \xmapsto{\Phi^0} \Phi^0(t_1,\dots,t_n) \in Q_A(0,\rho)\subseteq Q(0,\delta_1),
\]
such that $\Phi^0(0)=0$. 
Analogously, one has a diffeomorphism 
(corresponding to the "$w$-variables")
\[
Q_\rho \ni w \xmapsto{\Phi} \Phi(w) \in Q_A(0,\rho)\subseteq Q(0,\delta_1),
\] 
such that $\Phi(0)=0$, where 
\[
Q_\rho:=]-\rho,\rho[ \times ]-c_2\rho^{\kappa_2},c_2\rho^{\kappa_2}[\times \dots \times ]-c_n\rho^{\kappa_n},c_n\rho^{\kappa_n}[,
\]
with $1=\kappa_1 \leq \kappa_2 \leq \kappa_3 \leq\dots \leq \kappa_n=1/\varepsilon$ (cf. \cite{FP2}) and
$c_1=1,c_2, \dots,c_n$ universal constants. Note also that, 
by construction, $\Phi^0(Q)=\Phi(Q_\rho)=Q_A(0,\rho)$. 
In such a way, we have constructed the universal cube $Q_\rho$
centered at $0$ and the diffeomorphism $\Phi$ requested 
(note that we may consider also the diffeomorhpsim $\Phi^0$
and work on the standard cube $Q$ on $t$-variables,
but we prefer to work in the $w$-variables in order
to emphasize the dependence on $\rho$).

Furthermore, iterating the inclusion \eqref{FP1} 
(after $n$ integrations),
choosing $\rho>0$ sufficiently small and universal 
(such that all the changes of variables are done inside the resulting cubes),
one has
\[
B_{A_\Phi}(0,c\rho) \subseteq \lbrace \abs*{w_1} \leq \rho \rbrace \times \lbrace \abs*{w_2}\leq c_2\rho^{\kappa_2} \rbrace \times \cdot \cdot \cdot \times \lbrace \abs*{w_n}\leq c_n\rho^{\kappa_n} \rbrace \subseteq B_{A_\Phi}(0,C\rho),
\]
for some universal constants $c,C>0$, 
where $A_\Phi$ is the pullback of $A$
through the diffeomorphism $\Phi$ obtained
by composing all the change of variables,
and then 
\[
\Phi(Q_\rho)\subseteq B_A(0,C\rho).
\]

\noindent Finally, on the cube $Q_\rho$, we define
\[
\tilde{A}:=-\sum_{j=1}^nc_j^2\rho^{2\kappa_j}\frac{\partial^2}{\partial w_j^2}.
\]
If we now choose $j \in \lbrace 1,...,n \rbrace$
and consider the curve $\gamma_j:[0,1]\rightarrow \mathbb{R}^n$
defined by 
\[
\gamma_j(t):=e^{tc_j\rho^{\kappa_j}\frac{\partial}{\partial w_j}}(0),
\]
one has that this curve is subunit for $\tilde{A}$
(with speed less than $1$) and coming back to the appropriated step as explained before,
we obtain that the point $\gamma_j(1)=c_j\rho^{\kappa_j}$
corresponds to some point $x\in B_A(0,C\rho)$. 
Thus, by definition of subunit ball,
there exists $\gamma:[0,1]\rightarrow Q(0,\delta_1)\subseteq Q^\ast$
subunit curve for $A$ with speed less than $C\rho$, 
for some universal constant $C>0$, such that $\gamma(0)=0$ and $\gamma(1)=x$.
\end{proof}
Before commenting this theorem it is useful to give the following definition.
\begin{definition}\label{univcube}
We define \textit{the universal block of radius $\rho>0$} ($\rho>0$ universal and sufficiently small) centered at $0$ as 
\[
Q_\rho:=]-\rho,\rho[ \times ]-c_2\rho^{\kappa_2},c_2\rho^{\kappa_2}[\times \dots \times ]-c_n\rho^{\kappa_n},c_n\rho^{\kappa_n}[.
\]
\textit{The subunit cube of radius $\rho$ centered at $0$ associated with $A$} is instead defined as 
\[
Q_A(0,\rho):=\Phi(Q_\rho),
\]
i.e. as the image of the universal block through the diffeomorphism
$\Phi$ constructed in the proof of the previous theorem. 
\end{definition}
\begin{remark}\label{rmkcubeball}
Summing up, once again, in the proof of the theorem
we have constructed a universal block $Q_\rho$ in $w$-variables
(once again we prefer to use it instead of the
unit cube $Q$ in the $t$-variables because the dependence on $\rho$ is more evident) 
and an universal operator $\tilde{A}$ defined 
on it that is equivalent to our subelliptic operator $A$.\\
In other terms, we could say that the geometry
of the subunit cubes, and the subunit balls are equivalents, i.e.
\begin{equation}
\label{cubesandballs}
B_A(0,c\rho)\subseteq Q_A(0,\rho) \subseteq B_A(0,C\rho),
\end{equation}
for some $c,C>0$ universal constants, and $\rho>0$ universal
and sufficiently small. This equivalence allows us to conclude
that our operator $A$ is locally equivalent to the operator $\tilde{A}$, i.e. locally 
\[
A \sim \tilde{A}:=-\sum_{j=1}^n c_j^2\rho^{2\kappa_j} \frac{\partial^2}{\partial w_j^2},  
\]  
where the equivalence should be understood in the sense
that \emph{subunit curves for $\tilde{A}$ induce subunit curves
for $A$ and viceversa} (in the theorem we emphasized 
one direction of the equivalences but both are true thanks
to the inclusion \eqref{cubesandballs}).
\end{remark}

\begin{remark}
\label{rmkdilation}
Note also that, rescaling by a fixed constant, we can
construct another universal cube $Q_{\rho'}$ centered at $0$, 
where $\rho'=c\rho$, with $0<c<1$ positive constant, such that
\[
\Phi(Q_{\rho'})\subseteq \Phi(Q_\rho)\subseteq B_A(0,C\rho).
\] 
\end{remark}

To become familiar with this construction, it is very useful
to see it applied, for instance, to the case of the Grushin operator.
\begin{exa}
For $k \in \mathbb{N}$ the $k$-Grushin operator is defined as the operator
\[
A=-\frac{\partial^2}{\partial x_1^2}-x_1^{2k}\frac{\partial^2}{\partial x_2^2}, \quad x=(x_1,x_2)\in \mathbb{R}^2.
\]
Recall that in this case
the matrix of coefficients of the principal part is 
\[
\mathsf{A}_2(x)=(a_{ij}(x))_{i,j}=
\begin{pmatrix}
1 & 0 \\
0 & x_1^{2k}
\end{pmatrix}, \quad x \in \mathbb{R}^2.
\]
Our aim in this example is to calculate explicitly
the subunit cubes for the operator $A$. 
In the first place, note that the operator is already written
in "$(w_1,z)$-variables", i.e. has the form for which 
we can apply Lemma \ref{FPlemma}. 

\noindent Thus, we consider, for $\rho>0$, 
\[
\rho^2A=-\rho^2\frac{\partial^2}{\partial x_1^2}-\rho^2x_1^{2k}\frac{\partial^2}{\partial x_2^2}
\]
and the change of variables 
\[
\mathbb{R}^2\ni(t_1,t') \xmapsto{\psi^{(1)}_{\rho}} (\rho t_1,\rho t')=(x_1,x_2)\in \mathbb{R}^n
\]
(note that since in this case the operator 
is defined on $\mathbb{R}^2$
and not on a small cube it is not necessary
to choose $\rho$ small in this case).
Then, the pullback of $A$ through this diffeomorphism
is the operator $L=(\rho^2A)_{\psi_{1,\rho}}$ defined by 
\[
L=-\frac{\partial^2}{\partial t_1^2}-(t_1\rho)^{2k}\frac{\partial^2}{\partial {t'}^2}.
\]
By applying now Lemma \ref{FPlemma} to the operator $L$ one has 
\[
B_{L}(0,c) \subseteq \lbrace \abs*{t_1}\leq 1  \rbrace \times B_{\bar{L}}(0,\tilde{c}) \subseteq B_{L}(0,C),
\]
where the constants $c,\tilde{c},C$ are given 
by Lemma \ref{FPlemma}. Thus, coming back to the 
$(x_1,x_2)$ coordinates (i.e. in $(w_1,z)$ coordinates) we have
\begin{equation}\label{eqGrushin}
B_A(0,c\rho)\subseteq \lbrace \abs*{x_1}\leq \rho  \rbrace \times B_{\bar{A}}(0,\tilde{c}\rho) \subseteq B_{A}(0,C\rho).
\end{equation}
We next examine, the operator $\bar{A}$ (in the proof of Theorem \ref{fundtheorem})
that in this case has the form 
\[
\bar{A}=b(x_2)\frac{\partial^2}{\partial {x_2}^2},
\]
where 
\[
b(x_2):= \frac{1}{2\rho}\int_{-\rho}^{\rho}x_1^{2k}dx_1=\frac{\rho^{2k}}{2(2k+1)}.
\]
Thus 
\[
\bar{A}=\frac{\rho^{2k}}{2(2k+1)}\frac{\partial^2}{\partial {x_2}^2},
\]
and 
\[
\lbrace \abs*{x_2}\leq c_2 \rho^{k+1} \rbrace \subseteq B_{\bar{A}}(0,\tilde{c}\rho),
\]
for some universal $c_2>0$. 
In conclusion, by \eqref{eqGrushin},
\[
\lbrace \abs*{x_1} \leq \rho \rbrace \times \lbrace \abs*{x_2}\leq c_2 \rho^{k+1} \rbrace \subseteq B_{A}(0,C\rho).
\]
Note also that in this case 
\[
Q_A(0,\rho)=]-\rho,\rho[ \times ]-c_2\rho^{k+1},c_2\rho^{k+1}[
\]
and hence
\[
Q_A(0,\rho)\subseteq B_A(0,C\rho)\approx [-C\rho,C\rho]\times [-(C\rho)^{k+1},(C\rho)^{k+1}]
\]
for a universal constant $C>0$.
\end{exa}
\subsection{Subunit Cubes on Manifolds}
The first aim of this subsection is to prove that
we may choose an open cover of $M$ for which the coefficients
of the principal part of the local expressions of $A$ 
are uniformly bounded.
This property, as we saw in the previous subsection, is fundamental 
to construct a universal block (and a universal diffeomorphism). 
To reach this property, we may consider the finite cover constructed 
trough the exponential map (see Appendix \ref{app.exp}). 
For every $x \in M$, given $r:=r_g(M)>0$ we denote
by $B(0_x,r)$ the open ball in 
$T_xM$ with respect the metric $g_x$ and by
$B(x,r):=\mathrm{exp}_x(B(0_x,r))\subseteq M$.
Note also that by Theorem \ref{BandB_g} we have $B(x,r)=B_g(x,r)$. 

Since $M$ is compact, we may consider the
finite cover given by $M=\bigcup_{\ell=1}^m V_\ell$, for $m \in \mathbb{N}$,
where $V_\ell:=B(x_\ell,r)$ 
and the atlas (that is contained in the maximal atlas by definition)
\[
\mathscr{O}=\lbrace (U_\ell,\varphi_\ell); \ \ell=1,\dots L \rbrace,
\]
where $U_\ell:=B(x_\ell,2r)$ and $\varphi_\ell:=B^{-1}\circ \mathrm{exp}_{\bigl|_{B(0_{x_\ell},r)}}^{-1}$
(see Appendix \ref{app.exp}).  

By definition we have that the operator $A_\ell$
induced by $A$ on $\tilde{U}_\ell:=\varphi_\ell(\tilde{U}_\ell) \subseteq \mathbb{R}^n$,
that is \textit{the local expression} of $A$ on $\tilde{U}_\ell$, 
is a second-order differential operator of the form
\[
A_\ell=-\sum_{i,j=1}^n a^{(\ell)}_{ij}(\tilde{x})\frac{\partial^2}{\partial {x_i}\partial {x_j}}+\sum_{k=1}^nb^{(\ell)}_{k}(\tilde{x})\frac{\partial}{\partial x_k}+d^{(\ell)}(\tilde{x}), \quad \tilde{x} \in \tilde{U}_\ell,
\]
where $a^{(\ell)}_{ij}$, $b^{(\ell)}_k$, $d^{(\ell)} \in C^\infty(\tilde{U}_\ell,\mathbb{C})$
for all $i,j,k \in \lbrace 1,...,n \rbrace$.
It is immediate to note that the operator $A_\ell$
has the coefficients of the principal part that satisfy 
(with $\tilde{V}_\ell:=\varphi_\ell(V_\ell)$)
\[
\begin{split}
\forall \alpha \in \mathbb{N}^n & \ \exists \ C_{\alpha,\ell}>0  \ \mathrm{s.t.} \\
&\sup_{\tilde{x} \in \tilde{V}_\ell} \abs*{\partial_{\tilde{x}}^{\alpha}a^{(\ell)}_{ij}(\tilde{x})}\leq C_{\alpha,\ell}, \quad \forall i,j \in \lbrace 1,...,n \rbrace. 
\end{split}
\]
Therefore, for each $\alpha \in \mathbb{N}^n$
there exists a universal constant $C_\alpha>0$
(with $C_\alpha:=\max \lbrace C_{\alpha,1}, \dots C_{\alpha,m}\rbrace$), such that
\begin{equation}
\label{ABound}
\sup_{\ell=1,\dots,m}\sup_{\tilde{x} \in \tilde{V}_{\ell}} \abs*{\partial_{\tilde{x}}^{\alpha}a^{(\ell)}_{ij}(\tilde{x})}\leq C_{\alpha}, \quad \forall i,j=1,...,n. 
\end{equation}

The next goal is to further investigate the 
properties of the subunit metric. To do that we first define
a Riemannian metric on $B(x_\ell,r)$
as the Euclidean metric induced by $\varphi_{\ell}^{-1}$,
i.e. if $ x \in B(x_\ell,r)$, $v \in T_xM$ we put 
\[
\abs*{v}_x^E:=\abs*{\varphi'_{\ell}(x) v}_{E},
\]
where $\varphi'_{\ell}(x)$ is the tangent map 
at $x \in M$ (see \cite{L1}, p. 63).
This norm is locally equivalent to the Riemannian metric
defined on $M$ in the following sense. 
\begin{lemma}\label{lem.eg}
There exist constants $c_{x_\ell},C_{x_\ell}>0$
such that for every $x \in \overline{B(x_\ell,r)}$
and for every $v \in T_xM$, 
\[
c_{x_\ell}\abs*{v}_x^E\leq \abs*{v}^g_x \leq C_{x_\ell}\abs*{v}_x^E.
\]
In particular, if ${\gamma}:[0,1]\rightarrow M$
is a piecewise regular curve with ${\gamma}([0,1])\subseteq \overline{B(x_\ell,r)}$
one has
\begin{equation}
\label{Agammaequiv}
c_{x_\ell}L_E(\gamma)\leq L(\gamma)\leq C_{x_\ell} L_E(\gamma),
\end{equation}
where $L_E$ is the length of curve with respect the Euclidean metric 
and $L$ is the length with respect the Riemannian metric $g$ defined on $M$, respectively. 
\end{lemma}
\begin{proof}
Let us define
\[
K_\ell:=\lbrace (x,v); \ x \in \overline{B(x_\ell,r)}, \ v \in T_xM, \ \abs*{v}_x^E=1 \rbrace
\]
and a continuous
function $F:TM\rightarrow \mathbb{R}$ by $F(x,v):=\abs*{v}^g_x$.

\noindent Since $K_\ell$ is compact and $F$
is strictly positive on $K_\ell$,
there exist positive constants $c_{x_\ell}, C_{x_\ell}>0$ such that 
\[
c_{x_\ell}\leq \abs*{v}_x^g \leq C_{x_\ell}, \quad \forall (x,v) \in K_\ell.
\]
Finally, for each $x \in \overline{B(x_\ell,r)}$ and
for each $0\neq v \in T_xM$ we may write $v:=\lambda \bar{v}$,
where $\lambda=\abs*{v}_x^E$ and $\bar{v}:=\frac{v}{\abs*{v}_x^E}$. 
Then $(x,\bar{v}) \in K_\ell$ and 
\[
c_{x_\ell}\abs*{v}_x^E\leq \lambda c_{x_\ell}\leq \lambda \abs*{\bar{v}}^g_x=\abs*{v}^g_x=\lambda \abs*{\bar{v}}^g_x\leq \lambda C_{x_\ell}=C_{x_\ell}\abs*{v}_x^E.
\]
From this it immediately follows that for
every curve $\gamma$ with image contained 
in $\overline{B(x_\ell,r)}$ one has
\[
c_{x_\ell}L_E(\gamma)\leq L(\gamma)\leq C_{x_\ell} L_E(\gamma).
\]
\end{proof}
Now by Lebesgue's Lemma (see Lemma \ref{Lebesguelemma}) 
we have that there exists $\delta>0$ such that for
every $x \in M$ $\exists \ \ell=\ell(x) \in \lbrace 1,...,m \rbrace$ such that 
\[
B_g(x,\delta)\subseteq B(x_\ell,r),
\]
thanks to which we reach the next result. 
\begin{proposition}\label{geodballandcube}
There exists a universal constant $\nu>0$ such that, 
for every $x \in M$, if $B_g(x,\delta) \subseteq B(x_\ell,r)$ one has
\[
\varphi_{x_\ell}^{-1}(Q(\varphi_{x_\ell}(x),\nu)) \subseteq B_g(x,\delta).
\]
\end{proposition}
\begin{proof}
Since the Euclidean ball is equivalent to the
Euclidean cube (i.e. there exist $c,C>0$ universal constants
such that $Q(x,cr)\subseteq B_E(x,r) \subseteq Q(x,Cr)$
for every $x \in \mathbb{R}^n$ and $r>0$) we may reduce to proving 
that there exists $\nu>0$ such that, if $y, y' \in \tilde{V}_{\ell}$
satisfy $\abs*{y-y'}<\nu$ then $d_g(x,x')< \delta$, 
for $x=\varphi_{\ell}^{-1}(y)$ and $x'=\varphi_{\ell}^{-1}(y')$.
Let $\gamma_E:[0,1]\rightarrow \mathbb{R}^n$ be the curve
in $\mathbb{R}^n$ defined as $\gamma_E(t)=ty+(1-t)y'$ and 
$\gamma(t)=\varphi_{\ell}^{-1}(\gamma_E(t))$ the curve 
from $[0,1]$ to $M$ induced by $\gamma_E$. 
By \eqref{Agammaequiv}, setting
$\bar{c}:=\min \lbrace c_{x_1},...,c_{x_m} \rbrace$
and $\bar{C}:=\max \lbrace C_{x_1},...,C_{x_m} \rbrace$ we obtain
\[
\bar{c}L_E(\gamma)\leq L(\gamma)\leq \bar{C} L_E(\gamma).
\]
Thus 
\[
L(\gamma)\leq \bar{C}L_E(\gamma)=\bar{C}L_E(\varphi_{\ell} \circ \gamma)=\bar{C}L_E(\gamma_E)=\bar{C}\abs*{y-y'},
\]
(on the second and third inequality,
by an abuse of notation,
we have used $L_E$ for the length of a curve in $\mathbb{R}^n$
with respect the Euclidean metric) and so,
choosing $d_g(x,x')<\delta$, we get
$\nu<\frac{\delta}{\bar{C}}$.
\end{proof}

We finally show the next universal inclusion. 
\begin{proposition}\label{subballsconteuclballs}
If $C>0$ is a universal constant, there exists $\rho>0$ universal, such that, 
for every $x \in M$, when $ B_g(x,\delta)\subseteq B(x_\ell,r)$, one has
\[
B_{A_\ell}(\varphi_{\ell}(x),C\rho) \subseteq Q(\varphi_{\ell}(x),\nu),
\]
where $A_\ell$ is the operator corresponding to $A$ on the open set $\tilde{U}_{\ell}\subseteq \mathbb{R}^n$.
\end{proposition}
\begin{proof}
As before we prove that if $y \in B_{A_\ell}(\varphi_{\ell}(x),C\rho)$ 
for $\rho>0$ sufficiently small one has $\abs*{y-\varphi_{\ell}(x)}<\nu$.

\noindent If $y \in B_{A_\ell}(\varphi_{\ell}(x),C\rho)$ by definition
there exists $\gamma_A:[0,1]\rightarrow \mathbb{R}^n$ Lipschitz 
such that $\gamma_A(0)=\varphi_{\ell}(x)$, $\gamma_A(1)=y$ and 
\[
\abs*{\langle \dot{\gamma}_A(t), \xi \rangle}^2 \leq C^2\rho^2 \sum_{i,j=1}^n a_{ij}^{(\ell)}(\gamma(t))\xi_i\xi_j=C^2\rho^2\langle \mathsf{A}^{(\ell)}_2(\gamma(t))\xi,\xi \rangle, \quad \text{for a.e.\,t} \in [0,1],
\]
where the $a_{ij}^{(\ell)} \, 's$ are the coefficients of 
the operator on the open set $\tilde{U}_{\ell}\subseteq \mathbb{R}^n$. 
If for a.e. \, t $\in [0,1]$ we choose $\bar{\xi}=\bar{\xi}(t)$ 
linearly dependent on $\dot{\gamma}_A(t)$ we obtain
\begin{equation}
\label{CS1}
\abs*{\dot{\gamma}_A(t)}^2 \abs*{\bar{\xi}}^2 \leq C^2\rho^2\langle \mathsf{A}^{(\ell)}_2(\gamma_A(t)) \bar{\xi}, \bar{\xi} \rangle.
\end{equation}
Thus, since $A_\ell$ is defined on a dilation of $\tilde{V}_{\ell}$,
we may construct a continuous function $F$ on a compact set as in the proof 
of Lemma \ref{lem.eg}.
By repeating this argument for every $\ell=1,...,m$, we are able to find a
fixed constant $C'>0$ such that for all $\ell=1,...,m,$  $\tilde{x} \in \tilde{V}_{\ell}$, 
one has 
\[
\langle \mathsf{A}^{(\ell)}_{2}(x) \xi, \xi \rangle \leq C', \quad \forall \xi \in \mathbb{R}^n, \ \abs*{\xi}=1. 
\]
Hence, by \eqref{CS1},
\[
\abs*{\dot{\gamma}_A(t)}^2\leq C^2 \rho^2 \langle \mathsf{A}_2^{(\ell)}(\gamma(t))\frac{\bar{\xi}}{\abs*{\bar{\xi}}},\frac{\bar{\xi}}{\abs*{\bar{\xi}}} \rangle\leq C^2\rho^2C',
\]
so that $\abs*{\dot{\gamma_A}(t)}\leq C\rho \sqrt{C'}$
for a.e.\,t $\in [0,1]$. Therefore, 
choosing $\rho < \nu/(C\sqrt{C'})$, gives
\[
\abs*{y-\varphi_{\ell}(x)}=d_E(y,\varphi_{\ell}(x))\leq L_E(\gamma_A)=\int_0^1 \abs*{\dot{\gamma}_A(t)}dt<\nu
\]
and hence the proof is concluded.
\end{proof}

In conclusion, putting together these two last
propositions and the uniform bound of the coefficients 
in \eqref{ABound}, we prove Theorem \ref{theosubM}.

\begin{proof}[Proof of Theorem \ref{theosubM}]
The result follows directly from the previous propositions.
Indeed, for every $x \in M$, by Proposition \ref{geodballandcube} 
we can find a universal cube $Q(\tilde{x},\nu)$ on some open
coordinates set $\tilde{V}_{\ell}$ on $\mathbb{R}^n$, for which the image 
through $\varphi_{\ell}^{-1}$ is contained in the Lebesgue ball centered
at $x$ of geodesic radius $\delta$. Moreover, possibly by composing by a fixed
(universal) dilation and a translation we may assume that $\tilde{x}=\varphi_{\ell}(x)=0$ and $\nu=1$. 
\noindent At this point,
since the coefficients of the operator in coordinates
are uniformly bounded, choosing $\rho>0$ sufficiently small, 
we have that, by applying Theorem \ref{fundtheorem}, 
we may construct two universal blocks $Q_\rho$, $Q_\rho^{\ast}$, 
where the latter is a dilation by a fixed constant $c^\ast>1$ (that is $Q^\ast_\rho=Q_{c^\ast\rho}$).
Thus, the first two inclusions follow from
Theorem \ref{fundtheorem} (see also Remark \ref{rmkcubeball}
and Remark \ref{rmkdilation}).
The last inclusion is due to Proposition \ref{subballsconteuclballs},
by possibly shrinking $\rho$. 
Finally, the fact that the diffeomorphism $\Phi$ 
has uniform bounds follows from
Theorem \ref{fundtheorem} (cf. \cite{FSC}).
\end{proof}
Following this argument it is also possible to construct 
a covering of $M$ made by \textit{subunit cubes} on $M$ as follows.
 
Let $(M,g)$ be our closed Riemannian manifold with a
finite atlas 
\[
\mathscr{O}=\lbrace (B(x_\ell,r),\varphi_\ell); \ \ell=1,...,m \rbrace.
\]
By Lebesgue's Lemma there exists $\delta>0$ such that
if we take $x \in M$ one has $B_g(x,\delta) \subseteq B(x_\ell,r)$
for some $\ell=\ell(x)\in \lbrace 1,...,m \rbrace$. Now, by \eqref{geomrelationM}, for
every $y \in Q_\rho$, if we consider the integral curve
associated with the vector fields $c_j\rho^{\kappa_j}\partial/\partial w_j$
starting from $y$, we have that 
\begin{equation}\label{intecurvetilde}
\tilde{y}=e^{tc_j\rho^{\kappa_j}\frac{\partial}{\partial w_j}}y \in Q_\rho^\ast
\end{equation}
for $t \in ]-h,h[$, $h \in ]0,h_0]$, $h_0 \ll 1$
and 
\[
\varphi_{\ell}^{-1} \circ \Phi(\tilde{y}) \in B_{A}(x,C\rho).
\]
Then, for $x\in M$ we define
\[
Q_A(x,\rho):=\Psi_x^{-1}(Q_\rho),
\]
where $\Psi_x:=\Phi^{-1} \circ \varphi_{\ell}$,
and analogously $Q_{A}^\ast(x,\rho):=\Psi_x^{-1}(Q_\rho^\ast)$.
We are now in a position to define 
the subunit cubes of $M$.
\begin{definition}\label{def.subM}
\textit{The subunit cube on $M$ centered at $x$ of radius $\rho$} is defined 
as the open set $Q_{A}(x,\rho)\subseteq M$.
The open set $Q_{A}^\ast(x,\rho)$ is said to be \textit{a dilation} of the subunit cube. 
\end{definition}
By repeating this argument for every $x \in M$
we obtain a cover of $M$ of such subunit cubes
\[
M=\bigcup_{x \in M} Q_{A}(x,\rho).
\]
Therefore, by using the compactness of $M$ we can extract
a finite subcover that will allow us to construct a global diffusion.
\begin{remark}
Note that the inclusion 
\[
\lbrace \abs*{w_1}\leq \rho  \rbrace \times B_{\bar{A}}(0,\tilde{c}\rho) \subseteq B_A(0,C\rho),
\]
says that the end-point $\varphi_\ell \circ \Phi (\tilde{y})$ belongs to $B_A(x,C\rho)$,
so it can be reached by subunit curve of $A$, 
but this curve could be different from the curve 
$\eta(t)=\Psi_x^{-1}\circ e^{tc_j\rho^{\kappa_j}\partial/\partial w_j}(y)$.
Then, in the following, we will use the term "perturbed subunit curves"
referring to the curves $\eta(t)=\Psi_x^{-1}\circ e^{tc_j\rho^{\kappa_j}\partial/\partial w_j}(y)$, with $j=1,\dots,n$.
\end{remark}
\subsection{\texorpdfstring{Construction of a Global $(h,\rho)$-subelliptic Diffusion}{Construction of a Global (h,rho)-subelliptic Diffusion}}
The goal of this subsection is to construct,
with the tools introduced before, a global $(h,\rho)$-diffusion, 
with $h,\rho>0$ sufficiently small.
First of all, by compactness, using the construction above,
we may write, for some $N=N(\rho) \in \mathbb{N}$,
\[
M=\bigcup_{i=1}^N Q_{A}(x_i,\rho)
\]
and may consider the finite atlas 
\[
\lbrace (Q_{i}^\ast,\Psi_i); \ i=1,...,N \rbrace,
\]
and
\[
\lbrace (Q_{i},{\Psi_i}_{\bigl|_{Q_{i}}}); \ i=1,...,N \rbrace
\]
where, for each $i=1,...,n$, $\Psi_i:=\Psi_{x_i}$, $Q_{i}:=Q_{A}(x_i,\rho)$ and $Q_{i}^\ast:=Q_{A}^\ast(x_i,\rho)$.

Moreover, since the cover $\lbrace Q_{i} \rbrace_{i=1}^N$ is finite, 
it is possible to introduce 
two partitions of unity $\lbrace \psi_i \rbrace_{i=1}^N$,
$\lbrace \phi_i \rbrace_{i=1}^N$ subordinate to the covers
$\lbrace Q_{i}^\ast \rbrace_{i=1}^N$ and 
$\lbrace Q_{i} \rbrace_{i=1}^N$ respectively.

We need to
construct a non degenerate
$1$-density $d\mu \in \Gamma(M,\Omega_1(M))$ (see \cite{Z}) that will be used to study the convergence of the diffusion.
Define, for all $k=1,...,N$, 
\[
\omega_k:=\psi_k \cdot \Psi_k^\ast \, \abs*{dw^{(k)}}. 
\]
Since $\omega_k  \in \Gamma(M,\Omega_1(M))$, with $\mathrm{supp} \, \omega_k \subseteq Q^\ast_k$ we can define 
\[
d\mu:=\sum_{k=1}^N \omega_k \in \Gamma(M,\Omega_1(M)).
\]
Furthermore, due to the compactness of $M$, by rescaling with a fixed constant,
we may assume that $\int_M d\mu=1$ so that the measure $\mu$, defined through the $1$-density $d\mu$ as
\[
\mu(A):=\int_A d\mu, \quad A \subseteq M \ \text{ Borel measurable},
\]
is a probabibility measure $\mu \in \mathscr{P}(M)$.
Note also that the integrals are taken in the 
sense of the $1$-density.
The main goal of this section is to study the convergence
properties of the subunit diffusion to the measure.

We next define a rescaled version of the vector fields
defined before, i.e. we define, for $i=1,...,N, \, j=1,...,n$, the vector fields constructed on a fixed dilation $Q_\rho^\ast$
of the universal cube $Q_\rho$ on $\mathbb{R}^n$ as 
\[
\tilde{X}_{i,j}:= \frac{c^{(i)}_j\rho^{\kappa_j^{(i)}}}{\tilde{\gamma_i}} \frac{\partial}{\partial w_j^{(i)}}\in \Gamma(Q_\rho^\ast,TQ_\rho^\ast)
\]
and the vector fields
on the corresponding open set on $M$ as  
\[
X_{i,j}:=(\Psi_i^{-1})_\ast \, \tilde{X}_{i,j} \in \Gamma(Q_{i}^\ast ,TQ_{i}^\ast).
\]
\begin{remark}
Note that the vector fields $X_{i,j}$ actually depend on $\rho$
(since $\tilde{X}_{i,j}$ depends on $\rho$).
Thus we should write $X_{i,j,\rho}$,
but in what follows, in order to simplify
the notation, we will omit the subindex $\rho$.
\end{remark}

The aim is now to prove that these local vector fields
satisfy the so-called divergence-free condition (see Definition \ref{def.divfree})
with respect the localization of our strictly positive $1$-density $d\mu$,
that will allow us to obtain a self-adjoint property
for the operators that we are going to construct.

We have already seen that in our case
the $1$-density $d\mu$, for each $i=1,...,N$, 
corresponds locally to a smooth $1$-density $d\tilde{\mu}_i$
on the open sets $Q^\ast_\rho=\Psi_i(Q_i^\ast)$, 
written as $d\tilde{\mu}_i=\tilde{\gamma_i}\abs*{dw^{(i)}}$. 
The latter then determines the smooth $n$-form
$d\tilde{\nu}_i:=\tilde{\gamma_i} \, dw^{(i)}$ on $Q^\ast_\rho$.
In the following, by an abuse of notation,
we will sometimes identify $d\tilde{\nu}_i=d\tilde{\mu}_i$. 

\begin{proposition}\label{propdivfree}
For each $i=1,...,N$ the vector fields
$\tilde{X}_{i,j} \in \Gamma(Q_\rho^\ast,TQ_\rho^\ast)$ are
divergence-free with respect to the local measure $d\tilde{\nu}_i$, 
for all $j=1,...,n$. 
\end{proposition}
\begin{proof}
It suffices to use the local condition \eqref{localdivergence} of Remark \ref{rmk.divfree}. 
In fact, $d\tilde{\nu}_i=\tilde{\gamma_i} \, dw^{(i)}$ 
and $\tilde{X}_{i,j}=a_j \, \partial / \partial w^{(i)}_j$,
with $a_j(w):=(c_j^{(i)}\rho^{\kappa_j^{(i)}})/{\tilde{\gamma_i}}(w)$. Thus,
\[
\mathrm{div}_{\tilde{\nu}_i} \tilde{X}_{i,j}=\frac{\partial}{\partial w_j^{(i)}}\Bigl(\frac{c_j^{(i)}\rho^{\kappa_j^{(i)}}}{\tilde{\gamma_i}}\tilde{\gamma}_i\Bigr)= \frac{\partial}{\partial w_j^{(i)}}c=0,
\]
with $c>0$ a constant. 
\end{proof}
Finally, by possibly rescaling the constants $c_j>0$, one has that
the vector fields $X_{i,j}$ induce subunit curves of $A$ (with speed less than $C\rho$, where $C$ is the universal constant of Theorem \ref{fundtheorem}) say on $W_\rho$, where $Q_\rho \subseteq W_\rho \subseteq Q_\rho^\ast$.

Henceforth, we consider $\rho \in ]0,\rho_0]$ 
where $\rho_0 \in \mathbb{R}_+$ is fixed and small enough 
to reach the construction of the previous subsections 
and we choose $h \in ]0,h_0]$, with $h_0 \in \mathbb{R}_+$
being small enough to be fixed later. Essentially,
the parameters $h$ and $\rho$ represent respectively
the size of the step of our random walk and the
velocity at which it is carried out. 

To achieve our diffusion,
for $i=1,...,N$, $j=1,...,n$ we define the flow
\[
\begin{split}
\Phi_{i,j,h,\rho}:  ]-h,h[ & \times  M \rightarrow M \\
&(t,x) \mapsto 
\begin{cases}
e^{tX_{i,j}}x, \ \mathrm{if} \ (t,x) \in ]-h,h[ \times Q_{i} \\
x, \ \mathrm{otherwise},
\end{cases}
\end{split}
\]
and then construct the \textit{$(h,\rho)$-subelliptic diffusion} as follows. \\
$\bullet$ \textit{Step 0}. Fix $x_0 \in M$;\\
$\bullet$ \textit{Step k}. At step $k \in \mathbb{N}$ we choose 
at random $i \in \lbrace 1,...,N \rbrace$, $j \in \lbrace 1,...,n \rbrace$, $t \in ]-h,h[$ (with respect to the uniform probability) and set 
\[
x_{k+1}=
\begin{cases}
\Phi_{i,j,h,\rho}(t,x_k), & \text{if $x_k\in Q_{i}$ and $\Phi_{i,j,h,\rho}(t,x_k)\in Q_{i}$} \\
x_k, & \text{otherwise}.
\end{cases}
\]
Roughly speaking, the $(k+1)$-step of our $(h,\rho)$-subelliptic random walk
is the arrival point of the integral curve $e^{tX_{i,j}}x_k$ starting from $x_k$,
if both the initial and the final points
belong to the subelliptic cube $Q_{i}$,
or $x_{k+1}=x_k$ otherwise. 

We now want to associate a Markov operator with
this $(h,\rho)$-subelliptic diffusion.
To do that we firstly introduce for $i\in \lbrace 1,...,N \rbrace$, 
$j \in \lbrace 1,...,n \rbrace$, the auxiliary functions 
\[
\begin{split}
\chi_{i,j,h,\rho} \, : \,  ]-h,h[&\times M  \rightarrow \lbrace 0,1 \rbrace \\
&(t,x) \mapsto \chi_{i,j,h,\rho}(t,x)
\end{split}
\]   
defined as
\[
\chi_{i,j,h,\rho}(t,x):=\mathbbm{1}_{\lbrace(t,x) \in ]-h,h[\times M; \ (t,x) \in ]-h,h[ \times Q_{i}, \ \Phi_{i,j,h,\rho}(t,x) \in Q_{i}\rbrace}(t,x).
\]
Then, for $i \in \lbrace 1,...,N \rbrace$ and $j \in \lbrace 1,...,n \rbrace$ fixed we define
\[
S_{i,j,h,\rho}:\mathscr{M}_b(M,\mathbb{R}) \rightarrow \mathscr{M}_b(M,\mathbb{R})
\]
(recall that $\mathscr{M}_b(M,\mathbb{R})$ is the space of all bounded measurable functions on $M$) as
\[
S_{i,j,h,\rho}f(x):=M_{i,j,h,\rho}f(x)+T_{i,j,h,\rho}f(x)
\] 
where 
\[
M_{i,j,h,\rho}f(x):=m_{i,j,h,\rho}(x)f(x),
\]
with
\[
m_{i,j,h,\rho}(x):=1-\frac{1}{2h}\int_{-h}^h\chi_{i,j,h,\rho}(t,x)dt
\]
and
\[
T_{i,j,h,\rho}f(x):=\frac{1}{2h}\int_{-h}^h \chi_{i,j,h,\rho}(t,x)f(\Phi_{i,j,h,\rho}(t,x))dt.
\]
Hence, the $(h,\rho)$-Markov operator of our $(h,\rho)$-Markov process is an operator 
\[
S_{h,\rho}: \mathscr{M}_b(M,\mathbb{R}) \rightarrow \mathscr{M}_b(M, \mathbb{R})
\]
defined as
\begin{equation}\label{defShrho}
S_{h,\rho}:=\frac{1}{nN}\sum_{i=1}^N \sum_{j=1}^n S_{i,j,h,\rho}.
\end{equation}
Therefore, for $\rho \in ]0,\rho_0]$ sufficiently small fixed we have constructed a family of operators $(S_{h,\rho})_{h \in ]0,h_0]}$ that depend on the parameter $h$. In almost all situations we will work with $\rho>0$ fixed and we look at the behaviour for varying $h$. Nevertheless, in some situations we will investigate what happens to our results when $\rho$ goes to $0$.  
\begin{remark}
The operators $S_{i,j,h,\rho}$, and consequently $S_{h,\rho}$, are initially defined for real-valued functions, but they can naturally be extended as operators acting on complex-valued functions $f:M \rightarrow \mathbb{C}$, by setting, for $x \in M$,
\[
S_{i,j,h,\rho}f(x):=(S_{i,j,h,\rho} \mathrm{Re} \, f) (x)+i( S_{i,j,h,\rho} \mathrm{Im} \, f )(x).
\]
\end{remark} 
\subsection{\texorpdfstring{Basic Properties of the {$(h,\rho)$}-subelliptic Markov Operator}{ Basic Properties of the (h,rho)-subelliptic Markov Operator}}
In this subsection we work with $h \in ]0,h_0]$
and $\rho \in ]0,\rho_0]$ fixed and investigate the basic
fundamental properties of the
$(h,\rho)$-Markov operator $S_{h,\rho}$ 
constructed before, stated in the following proposition (cf. \cite{LM}, \cite{DLM1}, \cite{DLM2}). 

\begin{proposition}\label{1propshrho}
The operator $S_{h,\rho}$ satisfies $S_{h,\rho}1=1$ and $\norma*{S_{h,\rho}}_{L^\infty \rightarrow L^\infty}=1$. Furthermore, $S_{h,\rho}$ can be extended to a bounded self-adjoint
operator $S_{h,\rho}:L^2(M,d\mu)\rightarrow L^2(M,d\mu)$ such that $\norma*{S_{h,\rho}}_{L^2 \rightarrow L^2}=1$. In particular, $\mathrm{Spec}(S_{h,\rho})\subseteq [-1,1]$. 
\end{proposition}
\begin{proof}
The fact that $S_{h,\rho}1=1$ follows immediately 
by definition. We then focus to prove 
$\norma*{S_{h,\rho}}_{L^\infty \rightarrow L^\infty}=1$.
If we take $f \in C^\infty(\mathbb{R}^n)$, since
\[
S_{h,\rho}f= \frac{1}{n\cdot N}\sum_{j=1}^n \sum_{i=1}^N S_{i,j,h,\rho}f,
\]
it is sufficient to examine every terms in the sum of the RHS. 
For such terms we have to distinguish the following three cases. \\
$(1)$ Let $x \notin Q_{i}$. In this case $S_{i,j,h,\rho}f(x)=f(x)$, thus 
\[
\sup_{x \in Q_{i}^c}\abs*{S_{i,j,h,\rho}f(x)}\leq \norma*{f}_{L^\infty}.
\]
$(2)$ Let $x \in Q_{i}$ such that $\Phi_{i,j,h,\rho}(t,x) \in Q_{i}$
for all $t \in ]-h,h[$.
In this case
\[
S_{i,j,h,\rho}f(x)=\frac{1}{2h}\int_{-h}^hf(e^{tX_{i,j}}x)dt.
\]
Hence, 
\[
\sup_{x \in Q_{i} \atop \Phi_{i,j,h,\rho}(t,x) \in Q_{i} \; \forall t\in ]-h,h[}\abs*{S_{i,j,h,\rho}f(x)}\leq \norma*{f}_{L^\infty}.
\]
$(3)$ Let $x \in Q_{i}$ such that $\Phi_{i,j,h,\rho}(t,x) \in Q_{i}$
for all $t \in I_{x,h} \subsetneq ]-h,h[$.
Denoting $L_h :=\abs*{I_{x,h}} \in ]0,2h[$ one has
\[
S_{i,j,h,\rho}f(x)=(1-\frac{1}{2h}L_h)f(x)+\frac{1}{2h}\int_{I_{x,h}}f(e^{tX{i,j}}x)dt.
\]
Thus,
\[
\sup_{x \in Q_{i} \atop \Phi_{i,j,h,\rho}(t,x) \in Q_{i} \ \forall t \in I_{x,h}} \abs*{S_{i,j,h,\rho}f(x)}\leq (1-\frac{1}{2h}L_h)\norma*{f}_{L^\infty}+\frac{L_h}{2h}\norma*{f}_{L^\infty}=\norma*{f}_{L^\infty}. 
\]
In conclusion
\begin{equation}\label{eq.normLinfty}
\norma*{S_{i,j,h,\rho}f}_{L^\infty}\leq \norma*{f}_{L^\infty}
\end{equation}
and since $S_{h,\rho}1=1$ we may conclude that 
$\norma*{S_{h,\rho}}_{L^\infty \rightarrow L^\infty}=1$.

\noindent We now prove that the operator $S_{h,\rho}$ 
can be extended to a bounded self-adjoint operator 
$S_{h,\rho}:L^2(M,d\mu)\rightarrow L^2(M,d\mu)$.
To do that, since the global operator $S_{h,\rho}$ is a linear
combination of the operators $S_{i,j,h,\rho}$, in order to prove the self-adjointness property, it is sufficient to study, for $f,g \in C^\infty(M)$,
\[
\int_M (S_{i,j,h,\rho}f)\bar{g} d\mu=\sum_{k=1}^N \int_{Q_{k}}\phi_k(x)S_{i,j,h,\rho}f(x)\overline{g(x)}d\mu(x), 
\]
and in turn, for $k=1,...,N$ fixed, 
\[
\int_{Q_{k}}\phi_k(x)S_{i,j,h,\rho}f(x)\overline{g(x)}d\mu(x). 
\]
By definition of $S_{i,j,h,\rho}$, 
\[
\begin{split}
&\int_{Q_{k}}\phi_k(x)S_{i,j,h,\rho}f(x)\overline{g(x)}d\mu(x) \\
&=\int_{Q_{k}}\phi_k(x)m_{i,j,h,\rho}(x)f(x)\bar{g}(x)d\mu(x)+\int_{Q_{k}}\phi_k(x)T_{i,j,h,\rho}f(x)\overline{g(x)}d\mu(x).
\end{split}
\]
Hence, it is sufficient to examine the last term. Denoting $K=\mathrm{supp} \, \phi_k$ 
one has
\[
\begin{split}
&\int_{Q_{k}}\phi_k(x)T_{i,j,h,\rho}f(x)\overline{g(x)}d\mu(x)\\
&=\int_{Q_{k}}\phi_k(x)\Bigl(\frac{1}{2h}\int_{-h}^h\chi_{i,j,h,\rho}(t,x)f(e^{tX_{i,j}}x)dt\Bigr)\overline{g(x)}d\mu(x)\\
&=\int_{K}\phi_k(x)\Bigl(\frac{1}{2h}\int_{-h}^h\chi_{i,j,h,\rho}(t,x)f(e^{tX_{i,j}}x)dt\Bigr)\overline{g(x)}d\mu(x)
\end{split}
\]
Thus, by Fubini's Theorem and the change of variables
$y=\Phi_t(x)=e^{tX_{i,j}}x$, 
using the fact that by Proposition \ref{propdivfree}
the vector fields $X_{i,j}$ are divergence-free 
with respect to the measure $\mu$, we have 
\[
\begin{split}
&\int_{Q_{k}}\phi_k(x)T_{i,j,h,\rho}f(x)\bar{g}(x)d\mu(x)\\
&=\frac{1}{2h}\int_{-h}^h\Bigl(\int_{K\cap Q_{i}\cap \Phi^{-1}_{t}(Q_i)}\phi_k(x)\chi_{i,j,h,\rho}(t,x)f(e^{tX_{i,j}}x)\overline{g(x)}d\mu(x)\Bigr)dt \\
&=\frac{1}{2h}\int_{-h}^h\Bigl(\int_{\Phi_t(K\cap Q_{i})\cap Q_i}\phi_k(e^{-t{X_{i,j}}}y)\chi_{i,j,h,\rho}(t,e^{-tX_{i,j}}y)f(y)\overline{g(e^{-tX_{i,j}}y)}d\mu(y)\Bigr)dt\\
&= \frac{1}{2h}\int_{-h}^h\Bigl(\int_{Q_k}\phi_k(e^{-t{X_{i,j}}}y)\chi_{i,j,h,\rho}(t,e^{-tX_{i,j}}y)f(y)\overline{g(e^{-tX_{i,j}}y)}d\mu(y)\Bigr)dt    \\
&=\int_{Q_{k}}\overline{\Bigl(\frac{1}{2h}\int_{-h}^h\phi_k(e^{tX_{i,j}}y)g(e^{tX_{i,j}}y)\chi_{i,j,h,\rho}(t,y)dt\Bigr)}f(y)d\mu(y),
\end{split}
\]
where in the last equality we have used
that $\chi_{i,j,h,\rho}(t,e^{-tX_{i,j}}y)=\chi_{i,j,h,\rho}(-t,y)$ 
for every $(t,y) \in ]-h,h[\times M$. 

\noindent Therefore, summing over $k$, we obtain 
\[
(T_{i,j,h,\rho}f,g)_0=(f,T_{i,j,h,\rho}g)_0
\]
and so 
\[
(S_{h,\rho}f,g)_0=(f,S_{h,\rho}g)_0.
\]
Moreover, using again the divergence-free property as before, we have 
\[
\begin{split}
\norma*{S_{h,\rho}f}^2_{L^1}&=\int_M\Bigl|(1-\frac{1}{2h}\int_{-h}^h\chi_{i,j,h,\rho}(t,x)dt)f(x)\\
&\quad \quad +\frac{1}{2h}\int_{-h}^h \chi_{i,j,h,\rho}(t,x)f(\Phi_{i,j,h,\rho}(t,x))dt\Bigr|d\mu(x)\\
&\leq \norma*{f}_{L^1}^2, \quad f \in C^\infty(M).
\end{split}
\]
Therefore, by the Riesz-Thorin interpolation theorem (recall also \eqref{eq.normLinfty}), $S_{h,\rho}$ can be extended to a bounded self-adjoint operator $S_{h,\rho}:L^2(M,d\mu) \rightarrow L^2(M,d\mu)$ such that
\[
\norma*{S_{h,\rho}f}_0\leq \norma*{f}_0, \quad f \in L^2(M,d\mu),
\]
and since $S_{h,\rho}1=1$ we have once again $\norma*{S_{h,\rho}}_{L^2\rightarrow L^2}=1$.

\noindent Finally, since $S_{h,\rho}$ is self-adjoint
and $\norma*{S_{h,\rho}}_{L^2 \rightarrow L^2}=1$, 
by Proposition \ref{prop.specball} and Proposition \ref{prop.SA} we get
\[
\mathrm{Spec}(S_{h,\rho})\subseteq [-1,1].
\]
\end{proof}

We next note that the operator $S_{h,\rho}$
is positive on the space of continuous functions on $M$, i.e. 
\[
f  \in C^0(M), \; f\geq 0 \Rightarrow S_{h,\rho}f(x) \geq 0, \quad \forall x \in M.
\]
Therefore, for every $x \in M$, the functional 
\[
\begin{split}
F_x:C^0(M & ) \rightarrow \mathbb{R} \\
&f \mapsto S_{h,\rho}f(x)
\end{split}
\]
is a positive linear functional.
Hence, by the Markov-Riesz Theorem (see Theorem 5.41 in \cite{T}),  
there exists a unique Radon measure $s_{h,\rho}(x,\cdot)$ on $M$, such that  
\[
F_x(f)=\int_M f(y)s_{h,\rho}(x,dy),
\]
and then for every $x\in M$
\begin{equation}\label{eq.MR}
S_{h,\rho}f(x)=\int_M f(y)s_{h,\rho}(x,dy).
\end{equation}
Thus, the function $ M \ni x\mapsto s_{h,\rho}(x,\cdot) \in \mathscr{P}(M)$ 
is the Markov-Riesz kernel (see \cite{T})
associated with the operator $S_{h,\rho}$. 
From now on, we consider $S_{h,\rho}$ acting on $L^2(M,d\mu)$.
In general, for $f:M \rightarrow \mathbb{R}$ a bounded measurable function, we may write,
\[
S_{h,\rho}f=M_{h,\rho}f+T_{h,\rho}f,
\]
where, for $x \in M$,
\[
M_{h,\rho}f(x):=\frac{1}{nN}\sum_{i=1}^N\sum_{j=1}^n m_{i,j,h,\rho}(x)f(x),
\]
and
\[
T_{h,\rho}f(x):=\frac{1}{n N}\sum_{i=1}^N \sum_{j=1}^n T_{i,j,h,\rho}f(x).
\]
Hence, in terms of their distributional Kernels, one has 
\[
S_{h,\rho}f(x)=\int_M f(y)s_{h,\rho}(x,dy)=\int_M f(y)m_{h,\rho}(x,dy)+\int_M f(y)t_h(x,dy),
\]
where $m_{h,\rho}(x,dy)$ and $t_{h,\rho}(x,dy)$ are the Markov Kernels
associated respectively with the operators $M_{h,\rho}$ and $T_{h,\rho}$.
Our goal now is to study the $(h,\rho)$-Markov Kernel $s_{h,\rho}(x,dy)$
of the operator $S_{h,\rho}$. 
In the following we will also denote by $s_{h,\rho}^k(x,dy)$
the Markov Kernel of the operator $S_{h,\rho}^k=\underbrace{S_{h,\rho}\circ \dots \circ S_{h,\rho}}_{\textit{k-times}}$. 

Finally, to carry out our study
it will be necessary to introduce 
for all $h \in ]0,h_0]$ and $\rho \in ]0,\rho_0]$ 
the operator $\mathsf{A}_{h,\rho}$
acting on $L^2(M,d\mu)$, defined as 
\[
\mathsf{A}_{h,\rho}:=\frac{I-S_{h,\rho}}{h^2}
\]
and the associated bilinear form 
\[
\mathsf{B}_{h,\rho}(u,v):=(\mathsf{A}_{h,\rho}u,v)_0=\Bigl(\frac{I-S_{h,\rho}}{h^2}u,v\Bigr)_0, \quad u,v \in L^2(M,d\mu).
\]
We denote the corresponding \textit{Dirichlet form} by 
\begin{equation}\label{Dirichletform}
\mathsf{E}_{h,\rho}(u):=\mathsf{B}_{h,\rho}(u,u), \quad u \in L^2(M,d\mu).
\end{equation}

\begin{proposition}\label{proppositivity}
The operator $\mathsf{A}_{h,\rho}:L^2(M,d\mu) \rightarrow L^2(M,d\mu)$ 
is a positive self-adjoint operator. 
\end{proposition}
\begin{proof}
The operator $\mathsf{A}_{h,\rho}$ is self-adjoint 
because so is $S_{h,\rho}$.
Regarding the positivity we may reduce the proof
to showing that, with $u \in C^\infty(M)$,
\begin{equation}\label{usefulpositivity}
2((I-S_{i,j,h,\rho})u,u)_0=\int_M \frac{1}{2h}\int_{-h}^h\abs*{u(x)-u(e^{tX_{i,j}}x)}^2\chi_{i,j,h,\rho}(t,x)dtd\mu(x),
\end{equation}
since, by construction of $S_{h,\rho}$ it yields 
\[
((I-S_{h,\rho})u,u)_0 \geq 0.
\]
Note at first that, up to considering the real and the imaginary parts separately,
we may assume that $u$ is real-valued. We have  
\[
\begin{split}
&2((I-S_{i,j,h,\rho})u,u)_0\\
&=2((1-M_{i,j,h\rho}-T_{i,j,h,\rho})u,u)_0 \\
&=2\Bigl(\frac{1}{2h}\int_{-h}^h \chi_{i,j,h,\rho}(t,x)(u(x)-u(e^{tX_{i,j}}x))dt,u\Bigr)_0 \\
&=2\int_{Q_{i}}\Bigl(\frac{1}{2h}\int_{-h}^h\chi_{i,j,h,\rho}(t,x)u(x)^2dt\Bigr)d\mu(x)-2\Bigl(T_{i,j,h,\rho}u,u\Bigr)_0.
\end{split}
\]
As regards the RHS of \eqref{usefulpositivity}, we have
\[
\begin{split}
&\int_{M} \Bigl( \frac{1}{2h}\int^h_{-h}\abs*{u(x)-u(e^{tX_{i,j}}x)}^2\chi_{i,j,h,\rho}(t,x) dt \Bigr) d\mu(x) \\
&=\int_{Q_{i}} \Bigl( \frac{1}{2h}\int_{-h}^h\Bigl(u(x)^2-2u(x)u(e^{tX_{i,j}}x) +u(e^{tX_{i,j}}x)^2\Bigr)\chi_{i,j,h,\rho}(t,x)dt \Bigr) d\mu(x)\\
&=\int_{Q_{i}}\Bigl(\frac{1}{2h}\int_{-h}^hu(x)^2\chi_{i,j,h,\rho}(t,x)dt +\frac{1}{2h}\int_{-h}^h u(e^{tX_{i,j}}x)^2\chi_{i,j,h,\rho}(t,x)dt\Bigr)d\mu(x) \\
& \quad -2\Bigl(T_{i,j,h,\rho}u,u\Bigr)_0.
\end{split}
\]
Then \eqref{usefulpositivity} is satisfied if and only if 
\[
\begin{split}
\int_{Q_{i}}\Bigl(\int_{-h}^hu(x)^2\chi_{i,j,h,\rho}(t,x)dt\Bigr)d\mu(x)= \int_{Q_{i}}\Bigl(\int_{-h}^h u(e^{tX_{i,j}}x)^2\chi_{i,j,h,\rho}(t,x)dt\Bigr)d\mu(x)
\end{split}
\]
but this is true by applying Fubini's Theorem
and the change of variables. 
\end{proof}
\section{Auxiliary Results in Preparation to the Diffusion Theorems}\label{sec.aux}

In this section we will provide 
the necessary auxiliary results
to prove Theorem \ref{Theorem1} and Theorem \ref{Theorem2}.
In the first subsection we begin by investigating some useful
properties of the Markov Kernel (see Theorem \ref{lowestimatekernel} 
and the related corollaries).
Subsequently, in Theorem \ref{theohighlow} 
we will obtain a decomposition of the spectrum of $S_{h,\rho}$
on the \textit{high and low frequency parts} 
with respect to $h$ and through it,
we will also obtain a Weyl-type estimate.
Finally, in the last subsection we will examine the properties
of the \textit{infinitesimal generator} of our Markov process. 
 
\subsection{A Lower Bound for the Markov Kernel}

The goal of this subsection is to obtain
a lower bound for the $n$-th power
of the $(h,\rho)$-Markov Kernel 
and from this result it will be possible
to obtain an $L^\infty$ bound
for the eigenfunctions of the operator $S_{h,\rho}$.
The lower bound is to be understood in the sense of measures.
More precisely, we have the following statement.
 
\begin{theorem}\label{lowestimatekernel}
There exist $h_0>0$, and constants $c,\delta>0$,
such that, for all $x\in M$, for all $\rho \in ]0,\rho_0]$ and $h \in ]0,h_0]$, 
\[
s_{h,\rho}^n(x,dy)=r_h(x,dy)+ch^{-n}\mathbbm{1}_{\lbrace d_g(x,y)<\delta h \rbrace}d_gy,
\]
where $r_h(x,dy)$ is a positive Borel measure for all $x \in M$. 
\end{theorem}
\begin{proof}
We have to prove that there exist $h_0>0$
and constants $c,\delta>0$, independent of $h$, 
such that, for all $h \in ]0,h_0]$ 
and for every $f:M\rightarrow \mathbb{R}$ nonnegative and continuous function,
one has
\begin{equation}\label{globallowbound}
S_{h,\rho}^nf(x)\geq ch^{-n} \int_{\lbrace d_g(x,y)<\delta h \rbrace}f(y)d_gy.
\end{equation}
We may reduce the proof to showing the weaker assertion that
for all $x_0 \in M$ there exist constants
$h_0=h_0(x_0), \alpha=\alpha(x_0), c=c(x_0),\delta=\delta(x_0)>0$,
such that, for all $h \in ]0,h_0]$
and for every $f:M\rightarrow \mathbb{R}$ nonnegative and continuous function, 
\begin{equation}\label{locallowbound}
d_g(x,x_0)<\alpha \ \Rightarrow \ S_{h,\rho}^nf(x)\geq ch^{-n} \int_{\lbrace d_g(x,y)<\delta h \rbrace}f(y)d_gy,
\end{equation}
where $\alpha<r_g(M)$ (with $r_g(M)$ the injectivity radius).
In fact, by compactness $M=\cup_{\ell=1}^L B_g(x_\ell,\alpha(x_\ell))$,
for some $L \in \mathbb{N}$ and choosing $\delta=\min_\ell \, \delta(x_\ell)$,
$h_0=\min_\ell \, h_0(x_\ell)$ and $c=\min_\ell c(x_\ell)$, 
one easily proves that \eqref{locallowbound} implies \eqref{globallowbound}. 

Hence, we fix $x_0 \in M$.
Then $x_0 \in Q_{i_1}$ for some $i_1=1,...,N$, 
and, to simplify notation (up reordering of cubes)
we may assume $x_0 \in Q_{1}$.
We consider the set of $x \in M$ such that $d_g(x,x_0)<\alpha(x_0)$,
for some suitable $\alpha(x_0)>0$,
such that the image of the integral curves 
$\gamma_{x}:]-h_0,h_0[\rightarrow M$ of the vector fields $X_{1,j}$
starting from such $x$ are all contained in $Q_{1}$, for all $j=1,\dots,n$.
We need to examine the following two cases. \\
$(1)$ Suppose first $x \in Q_{1}$, $x \notin \cup_{i=2}^N Q_{i}$.
In this case
\[
S_{i,j,h,\rho}f(x)=f(x), \quad \text{for all $i=2,...,N$, $j=1,...,n$}.
\]
Thus (recall \eqref{defShrho}),
\[
S_{h,\rho}^nf(x)\geq\Bigl(\frac{1}{n N} \sum_{j=1}^n S_{1,j,h,\rho}\Bigr)^nf(x).
\]
Hence, since $f \geq 0$, there exists constants $c,\delta>0$, such that 
\[
\begin{split}
S_{h,\rho}^nf(x)&\geq \frac{1}{(nN)^n} \prod_{j=1}^n S_{1,j,h,\rho}f(x) \\
&=\frac{h^{-n}}{(2nN)^n}\int_{-h}^h \dots \int_{-h}^h f(e^{t_nX_{1,n}} \dots  e^{t_1X_{1,1}}x)dt_1 \dots dt_n \\
& \geq ch^{-n} \int_{\lbrace d_g(x,y)<\delta h \rbrace}f(y)d_gy,
\end{split}
\]
where the last inequality follows
from the fact that $X_{1,1}(x),...,X_{1,n}(x)$
is a basis for $T_xM$, $x \in Q_{i}$. 
 
$(2)$ Suppose now $x \in \cap_{\ell=1}^{N_1}Q_{i_\ell}$
for some $N_1 \in \mathbb{N}$, $2\leq N_1\leq N$. 
As before, w.l.o.g, we may assume $x \in \cap_{i=1}^{N_1} Q_{i}$
and since we deal with a finite of many cubes it is sufficient
to consider the case $x \in Q_{1} \cap Q_{2}$.
Denoting $S_{i,j}=S_{i,j,h,\rho}$, one has
\[
S_{h,\rho}^nf(x) \geq  \frac{1}{(nN)^n} S_{1,1}S_{2,1}S_{1,2}S_{2,1} \dots S_{1,n}S_{2,1}f(x).
\]
For simplicity, we examine the case $n=2$. In this case,
\[
S_{h,\rho}^2f(x) \geq  \frac{1}{(2N)^2} S_{1,1}S_{2,1}S_{1,2}S_{2,1}f(x)
\]
so that, we have to consider two cases. \\
$(i)$ We suppose first that the image of the integral curve
$\gamma_{x}:]-h_0,h_0[\rightarrow M$ of the vector field $X_{2,1}$ 
starting from $x$ is all contained in $Q_{2}$.
In this case 
\begin{equation}\label{shrinkingh}
\begin{split}
&S_{1,1,h,\rho}S_{2,1,h,\rho}S_{1,2,h,\rho}S_{2,1,h,\rho}f(x) \\
&\geq \frac{h^{-4}}{2^4(nN)^2}\int_{-h}^h \int_{-h}^h\int_{-h}^h\int_{-h}^h f(e^{t_4X_{2,1}} e^{t_3X_{1,2}}e^{t_2X_{2,1}}e^{t_1X_{1,1}}x)dt_1 dt_2dt_3dt_4 \\
& \geq \frac{h^{-2}}{2^4(nN)^2}\int_{-\delta_1}^{\delta_1} \int_{-h}^h\int_{-\delta_2}^{\delta_2}\int_{-h}^h f(e^{hs_4X_{2,1}} e^{t_3X_{1,2}}e^{hs_2X_{2,1}}e^{t_1X_{1,1}}x)dt_1 ds_2dt_3ds_4 \\
& \geq c'h^{-2} \int_{\lbrace d_g(x,y)<\delta h \rbrace}f(y)d_gy,
\end{split}
\end{equation}
for a possibly smaller constants $c,\delta>0$,
where the constants $\delta_1,\delta_2>0$
are choosen sufficiently small. \\
$(ii)$ If the integral curve associated with 
$X_{2,1}$ starting from $x$ 
is not all contained on $Q_{2}$, 
since $m_{i,j,h,\rho}\geq 0$, for all $i=1,\dots,N$, $j=1,\dots,n$,
the result follows in the same way.

The case $n>2$ is analogous
by repeating the same argument a finite number of times. 

In conclusion, we have obtained that
there exist $c>0$ and $\delta>0$ sufficiently small,
such that 
\[
S_{h,\rho}^nf(x)\geq ch^{-n} \int_{\lbrace d_g(x,y)<\delta h \rbrace}f(y)d_gy
\]
and thus the statement follows. 
\end{proof}

\begin{remark}\label{rmkpowerk}
Note that, by reasoning as in \eqref{shrinkingh} we get that the statement 
holds for all power $k\geq n$ of the operator. 
\end{remark}
As anticipated, this result allows us to obtain an estimate
for the eigenfunctions of our operator.
To do that, for $x \in M$ and $f:M \rightarrow \mathbb{R}$
bounded measurable function we write 
\[
R_{h,\rho}f(x)=\int_M f(y)r_{h,\rho}(x,dy),
\]
that is the operator with Markov Kernel $r_{h,\rho}(x,dy)$ defined by Theoreom \ref{lowestimatekernel} (note that such a kernel actually depends on $\delta$ given by Theorem \ref{lowestimatekernel}). With this notation we state a first corollary. 

\begin{corollary}\label{corlowest1}
There exist $h_0>0$, $\delta>0$ and $\tau<1$ such that
for all $h \in ]0,h_0]$ the operator $R_{h,\rho}$ 
satisfies
\[
\norma*{R_{h,\rho}}_{L^\infty\rightarrow L^\infty} \leq \tau <1.
\]
\end{corollary}
\begin{proof}
By definition, for $x \in M$,
\[ 
r_{h,\rho}(x,dy)=s_{h,\rho}^n(x,dy)-ch^{-n}\mathbbm{1}_{\lbrace d_g(x,y)<\delta h \rbrace}d_gy,
\]
for $h \in ]0,h_0]$ and for $c,\delta>0$ fixed constants. 
Hence, since $S_{h,\rho}1=1$, for $f \in L^\infty(M,d\mu)$, we have 
\[
\abs*{\int_Mf(x)r_{h,\rho}(x,dy)}\leq \norma*{f}_{L^\infty}\int_Mr_{h,\rho}(x,dy)\leq \norma*{f}_{L^\infty}\Bigl(1-ch^{-n}\inf_{x \in M}\int_{d_g(x,y) \leq \delta h}d_gy\Bigr),
\]
and so, by possibly shrinking $\delta>0$,
\[
\norma*{R_{h,\rho}f}_{L^\infty} \leq \tau \norma*{f}_{L^\infty}
\]
for some fixed constant $0<\tau<1$.
\end{proof}
As a consequence of this result the following estimate easily follows. 
\begin{corollary}\label{corlowest2}
Let $a \in ]\tau^{1/n},1]$ be a fixed constant.
Then, there exists $C=C(a)>0$ such that,
for any $\lambda \in [a,1]$ and any $f \in L^2(M,d\mu)$
eigenfunction $S_{h,\rho}f=\lambda f$ belonging to $\lambda$,
one has 
\[
\norma*{f}_{L^\infty} \leq Ch^{-n/2}\norma*{f}_0. 
\]
\end{corollary}
\begin{proof}
Let $\lambda \in [a,1]$ and $f \in L^2(M,d\mu)$
satisfy $S_{h,\rho}f=\lambda f$.
Then $S_{h,\rho}^nf=\lambda^nf$, whence
\begin{equation}\label{coreigen1}
\norma*{ch^{-n}\int_{\lbrace d_g(\cdot,y)<\delta h\rbrace}f(y)d_gy}_{L^\infty}\geq \lambda^n \norma*{f}_{L^\infty}-\tau \norma*{f}_{L^\infty}\geq c_a \norma*{f}_{L^\infty},
\end{equation}
with $c_a=a^n-\tau >0$. 

\noindent On the other hand, by Cauchy-Schwartz inequality, 
\begin{equation}\label{coreigen2}
\begin{split}
&\abs*{ch^{-n}\int_{\lbrace d_g(x,y)<\delta h\rbrace}f(y)d_gy}\\
&\leq ch^{-n}\mathrm{Vol}_g(\lbrace d_g(x,y) <\delta h\rbrace)^{1/2} \Bigl(\int_{\lbrace d_g(x,y)<\delta h\rbrace}\abs*{f(y)}^2d_gy\Bigr)^{1/2}\\
&\leq \tilde{c}h^{-n/2}\norma*{f}_0,
\end{split}
\end{equation}
for some $\tilde{c}>0$, 
where we have also used the fact that locally
(in coordinates on $\mathbb{R}^n$)
the measures $\mu$ and $\mathrm{Vol}_g$ are both 
absolutely continuous with respect to the Lebesgue measure. 

\noindent Therefore, putting together the estimates \eqref{coreigen1} and \eqref{coreigen2},
we conclude that there exists $C=C(a)>0$, such that 
\[
\norma*{f}_{L^\infty} \leq Ch^{-n/2}\norma*{f}_0. 
\]
\end{proof}
\subsection{Decomposition in High and Low Frequency Part}
For $u\in L^2(M,d\mu)$ that satisfies
the energy estimate \eqref{condhighlow}, in terms of its $L^2$
norm and of the Dirichlet form $\mathsf{E}_{h,\rho}$
(see \eqref{Dirichletform}),
we want to obtain an upper bound for
the so-called \textit{high-frequency part}
and \textit{low-frequency part} of $u$. 
One has the following result.
 
\begin{theorem}\label{theohighlow}
There exists $h_0>0$ such that for any $\rho \in ]0,\rho_0]$ and
any $h \in ]0,h_0]$ the following statement is true.

\noindent If $u \in L^2(M,d\mu)$ satisfies 
\begin{equation}\label{condhighlow}
\norma*{u}_0^2+\mathsf{E}_{h,\rho}(u) \leq 1,
\end{equation}
then there exist $(u^L_{h},u^H_{h}) \in H^1(M)\times L^2(M,d\mu)$ such that 
\[
u=u^L_{h}+u^H_{h}
\]
and
\begin{equation}\label{globalhighlow}
\norma*{u^L_{h}}_1 \leq C_0\rho, \quad \norma*{u^H_{h}}_0 \leq C_0\rho h,
\end{equation}
for some $C_0>0$.
\end{theorem}
\begin{proof}
Let us take $u \in L^2(M,d\mu)$
that satifies \eqref{condhighlow} and write 
\[
u=\sum_{i=1}^N \phi_i u,
\]
where, recall  $\lbrace \phi_i \rbrace_{i=1}^N$ is 
a partition of unity subordinate to the cover 
$\lbrace Q_{i} \rbrace_{i=1}^N$. 
Note that, if $\varphi \in C^\infty(M)$, 
there exists a constant $C=C_\varphi>0$, 
independent of $h \in ]0,h_0]$, such that 
\begin{equation}\label{Ehrhovarphiu}
\mathsf{E}_{h,\rho}(\varphi u)\leq C_\varphi(\norma*{u}^2_0 +\mathsf{E}_{h,\rho}(u)).
\end{equation}
Thus, for $i=1,...,N$ we have that $u_i:=\phi_iu$,
satisfies the estimate \eqref{condhighlow}. 
Fix now $i \in \lbrace 1,\dots,n \rbrace$
and work with $u_i:=\phi_iu$.
Moreover, 
write $c>0$ for a fixed absolute constant
possibly changing from line to line. 

\noindent First of all, by applying \eqref{usefulpositivity}
to $u_i$, we get that, for all $j=1,...,n$,  
\[
2((I-S_{i,j,h,\rho})u_i,u_i)_0=\int_M \frac{1}{2h}\int_{-h}^h\abs*{u_i(x)-u_i(e^{tX_{i,j}}x)}^2\chi_{i,j,h,\rho}(t,x)dtd\mu(x).
\]
Then, by \eqref{Ehrhovarphiu} one has 
\[
\mathsf{E}_{h,\rho}(u_i)=\Bigl(\frac{I-S_{h,\rho}}{h^2}u_i,u_i\Bigr)_0\leq c,
\]
for a fixed constant $c>0$, so that
\begin{equation}\label{estimatehl1}
\begin{split}
&\int_M \Bigl(\frac{1}{2h}\int_{-h}^h\abs*{u_i(x)-u_i(e^{tX_{i,j}}x)}^2\chi_{i,j,h,\rho}(t,x)dt\Bigr)d\mu(x)\\
&=\int_{Q_{i}} \Bigl(\frac{1}{2h}\int_{-h}^h\abs*{u_i(x)-u_i(e^{tX_{i,j}}x)}^2\chi_{i,j,h,\rho}(t,x)dt \Bigr)d\mu(x) \leq ch^2.
\end{split}
\end{equation}
Since $\mathrm{supp} \, u_i=K$ 
is compact and $K\subseteq Q_{i}$,
we may choose $h_0>0$ sufficiently small,
such that 
\begin{equation}\label{estimateshl2}
\begin{split}
&\int_{Q_{i}}\Bigl(\frac{1}{2h}\int_{-h}^h\abs*{u_i(x)-u_i(e^{tX_{i,j}}x)}^2\chi_{i,j,h,\rho}(t,x)dt\Bigr)d\mu(x)\\
&=\int_{Q_{i}}\Bigl(\frac{1}{2h}\int_{-h}^h\abs*{u_i(x)-u_i(e^{tX_{i,j}}x)}^2dt\Bigr)d\mu(x).
\end{split}
\end{equation}
In fact, we may choose $h_0>0$ sufficiently small such that
for every $x \in K$, $e^{tX_{i,j}}x \in Q_{i}$, 
for all $t \in ]-h_0,h_0[$.  
Hence in this case, assuming also that $u_i$ is real valued 
(by possibly considering separately the real and the imaginary part),
we have that 
\eqref{estimateshl2} is satisfied if and only if 
\begin{equation}
\label{estimateshl3}
\begin{split}
& \int_{Q_{i}} \Bigl(\frac{1}{2h}\int_{-h}^hu_i(e^{tX_{i,j}}x)^2dt\Bigr)d\mu(x) \\
&=\int_{Q_{i}} \Bigl(\frac{1}{2h}\int_{-h}^hu_i(e^{tX_{i,j}}x)^2\chi_{i,j,h,\rho}(t,x)dt\Bigr)d\mu(x).
\end{split}
\end{equation}
Thus, we have to prove \eqref{estimateshl3}. 
By applying Fubini's Theorem and 
the fact that $u_i(e^{tX_{i,j}}x)=0$ if 
$x \notin \lbrace y \in Q_{i}; \; e^{tX_{i,j}}y \, \in K \rbrace$ we have 
\[
\begin{split}
&\int_{Q_{i}}\Bigl(\frac{1}{2h}\int_{-h}^h u_i(e^{tX_{i,j}}x)^2\chi_{i,j,h,\rho}(t,x)dt\Bigr)d\mu(x) \\
&=\frac{1}{2h}\int_{-h}^h\Bigl(\int_{Q_{i}} u_i(e^{tX_{i,j}}x)^2\chi_{i,j,h,\rho}(t,x)d\mu(x)\Bigr)dt 
\\
&=\frac{1}{2h}\int_{-h}^h\Bigl(\int_{\lbrace x \in Q_{i}; \; e^{tX_{i,j}}x \, \in K \rbrace} u_i(e^{tX_{i,j}}x)^2\chi_{i,j,h,\rho}(t,x)d\mu(x)\Bigr)dt.
\end{split}
\]
At this point, since for all fixed
$t \in ]-h,h[$ the fact that
$x \in \lbrace y\in Q_{i}; \; e^{tX_{i,j}}y \, \in K \rbrace $, 
yields $\chi_{i,j,h,\rho}(t,x)=1$, we get
\[
\begin{split}
&\frac{1}{2h}\int_{-h}^h\Bigl(\int_{\lbrace x \in Q_{i}; \; e^{tX_{i,j}}x \, \in K \rbrace} u_i(e^{tX_{i,j}}x)^2\chi_{i,j,h,\rho}(t,x)d\mu(x)\Bigr)dt \\
&=\frac{1}{2h}\int_{-h}^h\Bigl(\int_{\lbrace x \in Q_{i}; \; e^{tX_{i,j}}x \, \in K \rbrace} u_i(e^{tX_{i,j}}x)^2d\mu(x)\Bigr)dt \\
&=\frac{1}{2h}\int_{-h}^h\Bigl(\int_{Q_{i}} u_i(e^{tX_{i,j}}x)^2 d\mu(x)\Bigr)dt,
\end{split}
\]
where in the last step we used 
again that $u_i(e^{tX_{i,j}}x)=0$ if
$x \notin \lbrace y \in Q_{i}; \; e^{tX_{i,j}}y \, \in K \rbrace$.

Thus, we have proved \eqref{estimateshl3} 
and consequently \eqref{estimateshl2}. 
In conclusion we have obtained 
\begin{equation}\label{prevest}
\int_{Q_{i}}\Bigl(\frac{1}{2h}\int_{-h}^h\abs*{u_i(x)-u_i(e^{tX_{i,j}}x)}^2dt\Bigr)d\mu(x)\leq ch^2,
\end{equation}
for all $j=1,\dots, n$. 
Moreover, in coordinates on $Q_\rho=\Psi_i(Q_{i})$,
we may rewrite estimate \eqref{prevest} as 
\[
\int_{Q_\rho}\Bigl( \frac{1}{2h}\int_{-h}^h\abs*{\tilde{u}_i(w)-\tilde{u}_i(e^{t\tilde{X}_{i,j}}w)}^2 dt\Bigr)dw\leq ch^2,
\]
where, $\tilde{X}_{i,j}=(\Psi_i)_\ast(X_{i,j})$ 
and $\tilde{u}_i:=u \circ \Psi_i^{-1}$
are defined on $Q_\rho \subseteq \mathbb{R}^n$.

\noindent Next, we extend $\tilde{u}_i$ as a function to $\mathbb{R}^n$
requiring that $\tilde{u}_i$ is vanishing outside $Q_\rho$.
We thus have that $\tilde{u}_i \in L^2(\mathbb{R}^n,d\tilde{\mu}_i)$ 
with $\mathrm{supp} \ \tilde{u} \subseteq Q_\rho$ 
and where $d\tilde{\mu}_i=(\Psi_i^{-1})^\ast d\mu_{\bigl|_{Q_{i}}}$
(extended by $0$ outside the cube $Q_\rho$).
Moreover, we may assume, after a rescaling,
that $\tilde{X}_{i,j}=c_j\rho^{\kappa_j}\frac{\partial}{\partial w_j}$ 
and $d\tilde{\mu}_i$ is the Lebesgue measure in coordinates.
 
\noindent Then, by applying the Fourier transform to $\tilde{u}_i$
(and dropping the index $i$)
in the variable $w_j$, and noting also
\begin{equation}
\begin{split}
\frac{1}{2h}\int_{-h}^h\abs*{1-e^{itc_j\rho^{\kappa_j}\xi_j}}^2dt&=\frac{1}{2h}\int_{-h}^h\Bigl((1-\cos(tc_j\rho^{\kappa_j} \xi_j))^2+(\sin (tc_j\rho^{\kappa_j}\xi_j))^2\Bigr)dt\\
&=\frac{1}{2h}\int_{-h}^h(2-2\cos (tc_j\rho^{\kappa_j}\xi_j))dt\\
&=2\Bigl(1-\frac{\sin (hc_j\rho^{\kappa_j} \xi_j)}{h c_j\rho^{\kappa_j} \xi_j}\Bigr),
\end{split}
\end{equation}
we obtain
\begin{equation}
\label{estimateshl4}
\begin{split}
& 2\int_{\mathbb{R}^n} \Bigl(1-\frac{\sin (hc_j\rho^{\kappa_j}\xi_j)}{h\rho^{\kappa_j}\xi_j}\Bigr)\abs*{\mathscr{F}_{w_j\rightarrow \xi_j}\tilde{u}}^2d\xi_jdw'= \\
& =\int_{\mathbb{R}^n} \Bigl(\frac{1}{2h}\int_{-h}^h\abs*{1-e^{itc_j\rho^{\kappa_j}\xi_j}}^2dt\Bigr)\abs*{\mathscr{F}_{w_j\rightarrow \xi_j}\tilde{u}}^2d\xi_jdw'\leq ch^2,
\end{split}
\end{equation}
where $w'=(w_1,...,\hat{w}_j,...,w_n)$.
At this point, note that there exists $c>0$ such that 
\[
\begin{split}
&(1-\frac{\sin(hc_j\rho^{\kappa_j}\xi_j)}{hc_j\rho^{\kappa_j}\xi_j})\geq ch^2c_j^2\rho^{2\kappa_j}\xi_j^2,  \quad \mathrm{for} \ hc_j\rho^{\kappa_j} \abs*{\xi_j} \leq b , \\
&(1-\frac{\sin(hc_j\rho^{\kappa_j}\xi_j)}{hc_j\rho^{\kappa_j}\xi_j})\geq c, \quad \mathrm{for} \  hc_j\rho^{\kappa_j} \abs*{\xi_j}>b,
\end{split}
\]
where $b>0$ is a fixed small constant,
and then decompose $\tilde{u}$ as 
\[
\begin{split}
\tilde{u}(w)&=\frac{1}{2\pi}\int_{hc_j\rho^{\kappa_j}\abs*{\xi_j} \leq b}e^{iw_jc_j\rho^{\kappa_j}\xi_j}\mathscr{F}_{w_j\rightarrow \xi_j}\tilde{u} \ d\xi_j \\
&\quad +\frac{1}{2\pi}\int_{hc_j\rho^{\kappa_j}\abs*{\xi_j}>b}e^{iw_jc_j\rho^{\kappa_j}\xi_j}\mathscr{F}_{w_j\rightarrow \xi_j}\tilde{u} \ d\xi_j. \\
\end{split}
\] 
Therefore, we choose a suitable cut-off function $\phi \in C_c^\infty(\mathbb{R}^n)$, $0 \leq \phi \leq 1$, with $\phi\equiv 1$ on $\mathrm{supp} \, \tilde{u}=K$ and $\phi\equiv 0$ on $Q_\rho^c$, and define
\[
\begin{split}
\tilde{u}^L_{j,h}(w)&:= \phi(w)\frac{1}{2\pi}\int_{hc_j\rho^{\kappa_j}\abs*{\xi_j} \leq b}e^{iw_jc_j\rho^{\kappa_j}\xi_j}\mathscr{F}_{w_j\rightarrow \xi_j}\tilde{u} \ d\xi_j \\ 
\tilde{u}^H_{j,h}(w)&:= \phi(w)
\frac{1}{2\pi}\int_{hc_j\rho^{\kappa_j}\abs*{\xi_j} > b}e^{iw_jc_j\rho^{\kappa_j}\xi_j}\mathscr{F}_{w_j\rightarrow \xi_j}\tilde{u} \ d\xi_j.
\end{split}
\]
The functions $\tilde{u}^L_{j,h}$
and $\tilde{u}^H_{j,h}$
are functions defined on $\mathbb{R}^n$ satisfying
\[
\tilde{u}=\tilde{u}^L_{j,h}+\tilde{u}^H_{j,h}
\]
that we call respectively
\textit{$j$-low frequency part} and \textit{$j$-high frequency part} of $\tilde{u}$. 
For them we have the following estimates.
\begin{lemma}
For every $j=1,...,n$ there exists $C_j>0$, independent of $h \in ]0,h_0]$, such that 
\begin{equation}
\begin{split}
& (i) \ \norma*{\partial_{w_j}\tilde{u}^L_{j,h}}_{L^2(\mathbb{R}^n)}\leq C_j;\\
& (ii) \ \norma*{\tilde{u}^H_{j,h}}_{L^2(\mathbb{R}^n)} \leq C_j\rho^{\kappa_j}h.
\end{split}
\end{equation}
\end{lemma} 

\begin{proof}
To prove $(i)$ we note that, by Minkowski's inequality, we have 
\[
\begin{split}
\norma*{\partial_{w_j}\tilde{u}^L_{j,h}}_{L^2(\mathbb{R}^{n})}&=\Bigl(\int_{Q_\rho}\abs*{\frac{1}{2\pi}\int_{hc_j\rho^{\kappa_j}\abs*{\xi_j} \leq b}\!\!\!\phi(w)c_j\rho^{\kappa_j}\xi_je^{iw_jc_j\rho^{\kappa_j}\xi_j}\mathscr{F}_{w_j\rightarrow \xi_j}\tilde{u} \, d\xi_j}^2dw\Bigr)^{\frac{1}{2}}\\
&\quad+\Bigl(\int_{Q_\rho}\abs*{\frac{1}{2\pi}\int_{hc_j\rho^{\kappa_j}\abs*{\xi_j} \leq b}(\partial_{w_j}\phi(w))e^{iw_jc_j\rho^{\kappa_j}\xi_j}\mathscr{F}_{w_j\rightarrow \xi_j}\tilde{u} \, d\xi_j}^2dw\Bigr)^{\frac{1}{2}}\\
&\leq \frac{1}{2\pi} \int_{hc_j\rho^{\kappa_j}\abs*{\xi_j} \leq b}\Bigl(\int_{Q_\rho}c_j^2\rho^{2\kappa_j}\xi_j^2\abs*{\phi(w)\mathscr{F}_{w_j\rightarrow \xi_j}\tilde{u}}^2 dw\Bigr)^{\frac{1}{2}} d\xi_j\\
&\quad +\Bigl(\int_{Q_\rho}\abs*{\frac{1}{2\pi}\int_{hc_j\rho^{\kappa_j}\abs*{\xi_j} \leq b}(\partial_{w_j}\phi(w))e^{iw_jc_j\rho^{\kappa_j}\xi_j}\mathscr{F}_{w_j\rightarrow \xi_j}\tilde{u} \, d\xi_j}^2dw\Bigr)^{\frac{1}{2}}.
\end{split}
\]
Now, since $h^{-2}(1-\frac{\sin(hc_j\rho^{\kappa_j}\xi_j)}{hc_j\rho^{\kappa_j}\xi_j})\geq cc_j^2\rho^{2\kappa_j}\xi_j^2$, 
for $hc_j\rho^{\kappa_j} \abs*{\xi_j} \leq b$,
by \eqref{estimateshl4}, we obtain
\[
\norma*{\partial_{w_j}\tilde{u}^L_{j,h}}_{L^2(\mathbb{R}^{n})}\leq C_j,
\]
for some constant $C_j>0$.

\noindent To get $(ii)$, we write
\[
\begin{split}
\norma*{\tilde{u}^H_{j,h}}_{L^2(\mathbb{R}^n)}&=\Bigl(\int_{Q_\rho}\abs*{\frac{1}{2\pi}\int_{hc_j\rho^{\kappa_j}\abs*{\xi_j} > b}e^{iw_jc_j\rho^{\kappa_j}\xi_j}\phi(w)\mathscr{F}_{w_j\rightarrow \xi_j}\tilde{u} \, d\xi_j}^2dw\Bigr)^{\frac{1}{2}}\\
&\leq \frac{1}{2\pi} \int_{hc_j\rho^{\kappa_j}\abs*{\xi_j} > b}\Bigl(\int_{Q_\rho}\abs*{\phi(w)\mathscr{F}_{w_j\rightarrow \xi_j}\tilde{u}}^2 dw\Bigr)^{\frac{1}{2}} d\xi_j.
\end{split}
\]
Thus, since $(1-\frac{\sin(hc_j\rho^{\kappa_j}\xi_j)}{hc_j\rho^{\kappa_j}\xi_j})\geq c$, 
for $hc_j\rho^{\kappa_j} \abs*{\xi_j}>b$,
possibly enlarging the constant $C_j$, once more by \eqref{estimateshl4}
we obtain
\[
\norma*{\tilde{u}^H_{j,h}}_{L^2(\mathbb{R}^n)}\leq C_j\rho^{\kappa_j}h.
\]
\end{proof}

\noindent Therefore, for all $j=1,\dots,n$ (recall we are working on each single $i$)
\[
\tilde{u}=\tilde{u}^L_{j,h}+\tilde{u}^H_{j,h}
\]
with
\[
\norma*{\rho^{\kappa_j}\partial_{w_j}\tilde{u}^L_{j,h}}_{L^2(\mathbb{R}^n)}\leq \rho^{\kappa_j}C_j, \quad \norma*{\tilde{u}^H_{j,h}}_{L^2(\mathbb{R}^n)}\leq C_j\rho^{\kappa_j}h,
\]
and so, there exists a constant $\bar{C}>0$,
independent of $h$, such that, for all $j=1,\dots,n$,
\[
\norma*{\tilde{X}_{i,j}\tilde{u}_{1,j,h}}_{L^2(\mathbb{R}^n)}\leq \rho^{\kappa_j}\bar{C}, \quad \norma*{\tilde{u}_{2,j,h}}_{L^2(\mathbb{R}^n)}\leq \rho^{\kappa_j}\bar{C}h.
\]
Our goal now is to make this decomposition independent of $j$.
To do that we will construct operators $A_h, B_{k,j,h}, D_{j,h}$, 
defined from $\mathscr{S}'(\mathbb{R}^n)$ to $\mathscr{S}'(\mathbb{R}^n)$,
depending on $h$, such that
\begin{equation}
\label{glueoperators}
\begin{split}
&1-A_h=\sum_{j=1}^n D_{j,h}hc_j\rho^{\kappa_j}\frac{\partial}{\partial w_j}, \\
&c_j\rho^{\kappa_j}\frac{\partial}{\partial w_j}A_h=\sum_{k=1}^nB_{k,j,h}c_j\rho^{\kappa_j}\frac{\partial}{\partial w_k},
\end{split}
\end{equation}
and $A_h, B_{k,j,h},D_{j,h}, D_{j,h} h\frac{\partial}{\partial w_j}, B_{k,j}h\frac{\partial}{\partial w_k}$ 
are bounded on $L^2(\mathbb{R}^n)$, uniformly in $h$.

\noindent If there exist such operators we define 
\begin{equation}\label{def.u1hu2h}
\tilde{u}^L_{h}:=A_h\tilde{u}, \quad \quad \tilde{u}^H_{h}:=(1-A_h)\tilde{u},
\end{equation}
and in analogy to what has been done before we call
them respectively the $h$-low frequency part of $\tilde{u}$
and the $h$-high frequency part of $\tilde{u}$.
 
\noindent By this decomposition of $\tilde{u}$ we get 
\[
\begin{split}
& c_j\rho^{\kappa_j}\frac{\partial}{\partial w_j}\tilde{u}^L_{h}=\sum_{k=1}^n c_j\rho^{\kappa_j}B_{k,j,h}\frac{\partial}{\partial w_k}\Bigl(\tilde{u}^L_{k,h}+h\frac{1}{h}\tilde{u}^L_{k,h}\Bigr), \\
& \tilde{u}^H_{h}=\sum_{j=1}^nD_{j,h}hc_j\rho^{\kappa_j}\frac{\partial}{\partial w_j}(\tilde{u}^L_{j,h}+\tilde{u}^H_{j,h})
\end{split}
\]
and consequently,
\begin{equation}
\label{estimatev}
\norma*{c_j\rho^{\kappa_j}\partial_{w_j}\tilde{u}^L_{h}}_{L^2(\mathbb{R}^n)} \leq C'\rho^{\kappa_j}, \quad \norma*{\tilde{u}^H_{h}}_{L^2(\mathbb{R}^n)}\leq C'h\rho^{\kappa_j},
\end{equation}
for some constant $C'>0$.

\noindent Now, through Lemma \ref{lemma.Fourier} we want
to construct the operators $A_h, B_{k,j}, D_j$
satisfying \eqref{glueoperators}.\\
To do that we need to introduce the $h$-scaling operator $\mathsf{T}_h$,
defined as
\[
\mathsf{T}_h g(\cdot):=h^{-n}g(\cdot/h), \quad g \in \mathscr{S}(\mathbb{R}^n)
\]
and if $F \in \mathscr{S}'(\mathbb{R}^n)$ as
\[
\langle \mathsf{T}_h F \, | \, \varphi \rangle:=\langle F | \varphi(h \; \cdot ) \rangle, \quad \varphi \in \mathscr{S}(\mathbb{R}^n).
\] 
Moreover, we take $\psi \in \mathscr{S}(\mathbb{R}^n)$ be a fixed function
such that $\int_{\mathbb{R}^n}\psi(w)dw=1$.

\noindent We first define 
\begin{equation}\label{def.Ah}
A_hf:=f\ast \mathsf{T}_h \psi, \quad f \in L^2(\mathbb{R}^n).
\end{equation}
Next, since $\int_{\mathbb{R}^n}\frac{\partial}{\partial w_j}\psi=0$, 
by Lemma \ref{lemma.Fourier}
we are able to find $\psi_{k,j} \in \mathscr{S}(\mathbb{R}^n)$ that solve the equation
\begin{equation}
\label{Bkjh}
c_j\rho^{\kappa_j}\frac{\partial}{\partial w_j}\psi=\sum_{k=1}^nc_j\rho^{\kappa_j}\frac{\partial}{\partial w_k}\psi_{k,j}.
\end{equation}
Therefore, we define
\[
B_{k,j,h}f:=f \ast \mathsf{T}_h\psi_{k,j}, \quad f \in L^2(\mathbb{R}^n).
\]
Finally, to define the operators $D_{j,h}$ we want to find solutions
$d_1,...,d_n \in \mathscr{S}(\mathbb{R}^n)$ of the equations
\begin{equation}
\label{Djh}
\sum_{j=1}^nc_j\rho^{\kappa_j}\frac{\partial}{\partial w_j}d_j=\delta_0-\psi.
\end{equation}
To do that, let $E$ be a fundamental solution
of the hypoelliptic equation
\[
\tilde{A}E=\delta_0,
\]
where $\tilde{A}=\sum_{j=1}^n c_j^2\rho^{2\kappa_j}\frac{\partial^2}{\partial w_j^2}$,
and let $\chi \in C_c^\infty(\mathbb{R}^n)$ with $\chi\equiv 1$
on a neighborhood of $0$.
Moreover, let $g_1,\dots ,g_n \in \mathscr{S}(\mathbb{R}^n)$,
that solve the equation
\begin{equation}
\label{lj}
\sum_{j=1}^n c_j\rho^{\kappa_j}\frac{\partial}{\partial w_j} g_j=\psi+\sum_{j=1}^nc_j\rho^{\kappa_j}\Bigl[\frac{\partial}{\partial w_j},\chi\Bigr]c_j\rho^{\kappa_j}\frac{\partial}{\partial w_j}E=:\psi_0.
\end{equation}
The solvability of equation \eqref{lj} follows from the fact that $\int_{\mathbb{R}^n}\psi(w)dw=1$ and
\[
\langle \sum_{j=1}^n\Bigl[c_j\rho^{\kappa_j}\frac{\partial}{\partial w_j},\chi\Bigr]c_j\rho^{\kappa_j}\frac{\partial}{\partial w_j}E|1 \rangle =-\langle \chi \tilde{A}E|1\rangle =-\langle \delta_0|\chi \rangle =-\chi(0)=-1,
\]
(note that we deal with compactly supported distribution thanks to presence of $\chi$) and hence $\psi_0 \in \mathscr{S}(\mathbb{R}^n)$ with $\int_{\mathbb{R}^n} \psi_0(w)dw=0$.

\noindent Therefore, for all $j=1,\dots,n$, we may define  
\[
d_j:=\chi c_j\rho^{\kappa_j}\frac{\partial}{\partial w_j}E-g_j, 
\]
and thanks to the equivalences 
\[
\begin{split}
&\sum_{j=1}^nc_j\rho^{\kappa_j}\frac{\partial}{\partial w_j}d_j=\delta_0-\psi\\
& \Longleftrightarrow \sum_{j=1}^nc_j\rho^{\kappa_j}\frac{\partial}{\partial w_j}(\chi c_j\rho^{\kappa_j}\frac{\partial}{\partial w_j}E-g_j)=\delta_0-\psi \\
& \Longleftrightarrow \sum_{j=1}^n c_j\rho^{\kappa_j}\frac{\partial}{\partial w_j} g_j =\psi+\sum_{j=1}^nc_j\rho^{\kappa_j}\Bigl[\frac{\partial}{\partial w_j},\chi\Bigr]c_j\rho^{\kappa_j}\frac{\partial}{\partial w_j}E=:\psi_0,
\end{split}
\]
we have that $d_1, \dots, d_n$ solve equation \eqref{Djh}.
Therefore, for $j=1,\dots, n$, we define
\[
D_{j,h}f:=f \ast \mathsf{T}_h d_j.
\]
Furthermore, since $\mathsf{T}_h\delta_0=\delta_0$, we have that \eqref{Bkjh} and \eqref{Djh} are equivalent to 
\[
\begin{split}
&1-A_h=\sum_{j=1}^n D_{j,h}hc_j\rho^{\kappa_j}\frac{\partial}{\partial w_j}, \\
&c_j\rho^{\kappa_j}\frac{\partial}{\partial w_j}A_h=\sum_{k=1}^nB_{k,j,h}c_j\rho^{\kappa_j}\frac{\partial}{\partial w_k}.
\end{split}
\]
Summing up, we have constructed the operators that satisfy \eqref{glueoperators}.\\
We next have to prove the boundedness properties for them, i.e. that the operators 
\[
A_h, B_{k,j,h}, 
D_{j,h}, D_{j,h} hc_j\rho^{\kappa_j}\frac{\partial}{\partial w_j}, B_{k,j}h\frac{\partial}{\partial w_k}
\]
are bounded on $L^2(\mathbb{R}^n)$, uniformly in $h$.\\
First of all, if $f\in L^2(\mathbb{R}^n)$, by Young's inequality 
\[
\norma*{A_hf}_0\leq \norma*{\mathsf{T}_h\psi}_{L^1}\norma*{f}_0
\]
and since $\norma*{\mathsf{T}_h\psi}_{L^1}=\norma*{\psi}_{L^1}$,
for all $h \in ]0,h_0]$ we have that the operator $A_h$ is bounded in $L^2$,
uniformly in $h \in ]0,h_0]$, 
and by the same argument the operators $B_{k,j,h}$, $D_{j,h}$ are bounded in $L^2$ as well. 

\noindent Furthermore, 
\[
\begin{split}
B_{k,j,h}h\frac{\partial}{\partial w_k}f(w)&=\Bigl(h\frac{\partial}{\partial w_k}f \ast \mathsf{T}_h\psi_{k,j}\Bigr)(w)\\
&=\Bigl(f \ast h\frac{\partial}{\partial w_k}(\mathsf{T}_h \psi_{k,j})\Bigr)(w)\\
&=f \ast \Bigl(\mathsf{T}_h(\frac{\partial}{\partial w_k}\psi_{k,j})\Bigr)(w),
\end{split}
\]
thus, since $\frac{\partial}{\partial w_k}\psi_{k,j}$ is a Schwartz function,
we conclude that also the operators $B_{k,j,h}h\frac{\partial}{\partial w_k}$ are bounded in $L^2(\mathbb{R}^n)$, uniformly in $h \in ]0,h_0]$.

\noindent It remains to examine the 
operators $D_{j,h}hc_j\rho^{\kappa_j}\frac{\partial}{\partial w_j}$.
Since 
\[
\norma*{\mathsf{T}_hf}_{L^2(\mathbb{R}^n)}=h^{-\frac{n}{2}}\norma*{f}_{L^2(\mathbb{R}^n)},
\] 
it suffices to study
\[
L^2(\mathbb{R}^n) \ni f\mapsto f \ast c_j\rho^{\kappa_j}\frac{\partial}{\partial w_j}d_j \in L^2(\mathbb{R}^n).
\]
By construction one has  
\begin{equation}\label{eqFoll}
c_j\rho^{\kappa_j}\frac{\partial}{\partial w_j}d_j=\chi c_j^2\rho^{2\kappa_j}\frac{\partial^2}{\partial w_j^2}E+r_j, \quad r_j \in \mathscr{S}(\mathbb{R}^n),
\end{equation}
so, the distribution $c_j^2\rho^{2\kappa_j}\frac{\partial^2}{\partial w_j^2}E$ is of type $0$ (see \cite{Fo}, pp. 164-167 for the terminology and the proof) and hence 
\[
c_j^2\rho^{2\kappa_j}\frac{\partial^2}{\partial w_j^2}E=k_j\delta_0+PV(f_j),
\]
where $f_j \in C^\infty(\mathbb{R}^n\setminus \lbrace 0 \rbrace)$ and $k_j>0$ is constant.
Moreover, the operator 
\[
L^2(\mathbb{R}^n) \ni f\mapsto f \ast c_j^2\rho^{2\kappa_j}\frac{\partial^2}{\partial w_j^2}E \in L^2(\mathbb{R}^n)
\]
is bounded and consequently, by \eqref{eqFoll} 
\[
L^2(\mathbb{R}^n) \in f\mapsto f \ast c_j\rho^{\kappa_j}\frac{\partial}{\partial w_j}d_j \in L^2(\mathbb{R}^n)
\]
is bounded as well. 
We thus have constructed the operators as
of \eqref{glueoperators} and then the estimates \eqref{estimatev} are proved.

\noindent Finally, starting from \eqref{estimatev} we wish to obtain the global estimates \eqref{globalhighlow}. For such purpose let $\chi_i \in C_c^\infty(Q_{i})$ be $0\leq \chi_i \leq 1$, such that $\chi_i \equiv 1$ in a neighborhood of $\mathrm{supp} \, \phi_i$.
Then $(\Psi_i^{-1})^\ast \chi_i$ is a smooth function defined on $\mathbb{R}^n$ with support in $Q_\rho$ and $(\Psi_i^{-1})^\ast\chi_i\equiv 1$ in a neighborhood of $\mathrm{supp} \, \tilde{u}_i$.  

\noindent Thus, setting
\[
u^L_{i,h}:=\Psi_i^\ast\bigl(((\Psi_i^{-1})^\ast(\chi_i))\tilde{u}^L_{i,h}\bigr), \quad u^H_{i,h}:=\Psi_i^\ast\bigl(((\Psi_i^{-1})^\ast(\chi_i))\tilde{u}^H_{i,h}\bigr),
\]
one has 
\[
\begin{split}
u^L_{i,h}+u^H_{i,h}&=\Psi_i^\ast\bigl(((\Psi_i^{-1})^\ast(\chi_i))\tilde{u}^L_{i,h}\bigr)+\Psi_i^\ast\bigl(((\Psi_i^{-1})^\ast(\chi_i))\tilde{u}^H_{i,h}\bigr)\\
&=\chi_i\Psi_i^{\ast}\tilde{u}_i\\
&=\chi_iu_i=u_i
\end{split}
\]
in $Q_{i}$ and 
\begin{equation}\label{lowhigestlocal}
\norma*{u^H_{i,h}}_0\leq C\rho^{\kappa_j}h, \quad \sup_{1\leq j \leq n}\norma*{X_{i,j}u^L_{i,h}}_0\leq C\rho^{\kappa_j},
\end{equation}
where $C>0$ is a constant depending on $\chi_i$, $\phi_i$ 
and other auxiliary functions but independent of $h$ and $\rho$.
Until now we get the decomposition for $u_i=\phi_i u$, for all $i=1,...,N$.
To obtain this result for $u \in L^2(M,d\mu)$ it is sufficient to write
\[
u=\sum_{i=1}^N \phi_iu=\sum_{i=1}^N u^L_{i,h}+\sum_{i=1}^N u^H_{i,h}=:u^L_{h}+u^H_{h}.
\] 
Since $u^L_{h}$ and $u^H_{h}$ are constructed
by local pieces satisfying \eqref{lowhigestlocal},
they obviously satisfy
\[
\norma*{u^L_{h}}_1 \leq C_0\rho, \quad \norma*{u^H_{h}}_0 \leq C_0\rho h,
\]
for some $C_0>0$ universal constant. 
\end{proof}
\subsection{A Weyl-type Estimate}\label{sec.weylest}
The purpose of this subsection is to give a first result
(see Theorem \ref{theospec})
concerning the spectrum of the family of operators
$(\mathsf{A}_{h,\rho})_{h \in ]0,h_0]}$, for $\rho \in ]0,\rho_0]$, where, recall,
\[
\mathsf{A}_{h,\rho}:L^2(M,d\mu)\rightarrow L^2(M,d\mu)
\]
is defined as 
\[
\mathsf{A}_{h,\rho}=\frac{I-S_{h,\rho}}{h^2}.
\]
where $S_{h,\rho}$ is the $(h,\rho)$-Markov operator given in \eqref{defShrho}.
This result, in the next section, will lead to a characterization 
of the spectrum of our $S_{h,\rho}$.
Let us start by fixing some notation.

In what follows, we denote by $\Lambda=\Lambda^\ast\geq \mathrm{Id} \in \Psi^2(M)$
a fixed elliptic pseudodifferential operator on $M$,
where $\Lambda^\ast$ is the formal adjoint of $\Lambda$ on $L^2(M,d\mu)$
(for the existence of such operators see for instance Lemma 7.1 in \cite{SH}). 
Since $M$ is compact there exists an orthonormal basis
$(e_j)_{j \geq 1}$ of eigenfunctions of $\Lambda$ in $L^2(M,d\mu)$ 
corresponding to the eigenvalues $\nu_1 \leq \nu_2 \leq \dots$ 
counted by multiplicity.
Moreover, by the classical Weyl law (see for instance \cite{GS}, p. 138) one has
\[
\sharp \lbrace j; \ \nu_j^{1/2}\leq r \rbrace \approx r^{\mathrm{dim}(M)}, \quad r \rightarrow +\infty.
\] 
Note that for $s \in \mathbb{R}$,
the topology on $H^s(M)$ is equivalent to the topology defined by the norm
\[
\norma*{f}'_{s}:=(\sum_{j} \nu_{j}^{2s} \abs*{f_j}^2)^{1/2}=(\Lambda^sf,f)_0^{1/2}, \quad f=\sum_jf_je_j \in H^s(M).
\]
Finally, we define $(P_h)_{h \in ]0,h_0]}$, be a family of self-adjoint bounded operators defined on $L^2(M,d\mu)$, for $u=\sum_{j \geq 1} \hat{u}_je_j \in L^2(M,d\mu)$, as
\[
P_hu:=\sum_{j \geq 1} \min(h^{-1},\nu_j^{s/2})\hat{u}_je_j.
\]
The goal now is to characterize the spectrum of the operators for which Theorem \ref{theohighlow} applies. In Lemma \ref{lemspectrum} we will give such a characterization in a more general framework (cf. Appendix of \cite{LM}) and to do that we first need the following estimate (which yields to the control of the eigenvalues of the operator $I+T_h$ 
obtained by the Theorem \ref{minmax} principle).
\begin{proposition}\label{propweyl2}
Let $(T_h)_{h \in ]0,h_0]}$, for $h_0>0$, be a family of nonnegative self-adjoint bounded operators acting on $L^2(M,d\mu)$ that satisfy the following property.\\
There exists a constant $C_0>0$, independent of $h$ such that
\begin{equation}\label{propTh}
\begin{split}
&\textit{if $u \in L^2(M,d\mu)$ satisfies $((I+T_h)u,u)_0 \leq 1$} \\
&\Rightarrow \ \exists \, s>0 \ \text{and} \ (u_1,u_2) \in H^s(M) \times L^2(M,d\mu), \ \text{s.t.} \\
& \quad \quad u=u_1+u_2 \ \text{and} \ \norma*{u_1}_s \leq C_0, \ \norma*{u_2}_0\leq C_0h.
\end{split}
\end{equation}
Then, for all $h \in ]0,h_0]$, 
\begin{equation}
\label{discretespectrum1}
\norma*{P_hu}_0^2 \leq 4C_0^2((I+T_h)u,u)_0.
\end{equation}
\end{proposition}
\begin{proof}
Note first that if $u \in L^2(M,d\mu)$ and $u=u_1+u_2$, with $u_1 \in H^s(M)$ and $u_2 \in L^2(M,d\mu)$ one has 
\begin{equation}\label{discretespectrum0}
\norma*{P_h u}_0^2 \leq 2\norma*{P_h u_1}_0^2+2\norma*{P_h u_2}_0^2\leq 2(\norma*{u_1}_s^2+h^{-2}\norma*{u_2}_0^2).
\end{equation}
Suppose by contradiction there exists $h \in ]0,h_0]$ and $u_\ast \in L^2(M,d\mu)$ such that 
\begin{equation}\label{lembycontr}
4C_0^2((I+T_h)u_\ast,u_\ast)_0 < \norma*{P_hu_\ast}_0^2.
\end{equation}
Hence, setting 
\[
C_\ast:=\sqrt{((I+T_h)u_\ast,u_\ast)_0}>0
\]
(note $((I+T_h)u_\ast,u_\ast)_0>0$ as $T_h\geq 0$) one has
\[
\Bigl((I+T_h)\frac{u_\ast}{C_\ast},\frac{u_\ast}{C_\ast}\Bigr)_0\leq 1
\]
and so, by hypothesis \eqref{propTh}, there exist $v\in H^s(M)$ and $w \in L^2(M,d\mu)$ such that 
\[
\frac{u_\ast}{C_\ast}=v+w, \quad \text{with} \ \norma*{v}_s \leq C_0, \ \norma*{w}_0 \leq C_0h.
\]
Thus, defining $u_{\ast,1}:=C_\ast v$ and $u_{\ast,2}:=C_\ast w$, we have that $u_{\ast,1} \in H^s(M)$, $u_{\ast,2} \in L^2(M,d\mu)$ 
\[
u_\ast=u_{\ast,1}+u_{\ast,2}, \quad \text{with} \ \norma*{u_1^\ast}_s \leq C_\ast C_0, \ \norma*{u_2^\ast}_0 \leq C_\ast C_0h.
\] 
Therefore, by \eqref{discretespectrum0} and \eqref{lembycontr},
\[
\begin{split}
4C_0^2((I+T_h)u^\ast,u^\ast)_0&<2(\norma*{u_1^\ast}_s^2+h^{-2}\norma*{u_2^\ast}_0^2)\\
& \leq 4C_0^2C_\ast^2 \\
&= 4C_0^2((I+T_h)u^\ast,u^\ast)_0
\end{split}
\] 
and then the contradiction.
\end{proof} 
We next have the following Lemma (cf. Appendix of \cite{LM}). 
\begin{lemma}\label{lemspectrum}
Let $(T_h)_{h \in ]0,h_0]}$, for $h_0>0$, be a family of nonnegative self-adjoint bounded operators acting on $L^2(M,d\mu)$ that satisfy property \eqref{propTh}.

\noindent Then, if $C_1<1/(4C_0^2)$
 we have that $\mathrm{Spec}(T_h)\cap [0,\zeta-1]$ is discrete for all $\zeta \leq C_1h^{-2}$ and  there exists $C_2>0$, independent of $h$, such that
\[
\sharp (\mathrm{Spec}(T_h) \cap [0,\zeta-1]) \leq C_2 \braket{\zeta}^{\mathrm{dim}(M)/2s}, \quad \text{for all} \ \zeta \leq C_1h^{-2}.
\]
\end{lemma}
\begin{proof}
By Proposition \ref{propweyl2}, for all $j \geq 1$, we have 
(with $\mu_j$ defined by the min-max principle see Theorem \ref{minmax}).
\begin{equation}
\label{discretespectrum2}
\mu_j(I+T_h) \geq \frac{1}{4C_0^2}\mu_j(P_h^2).
\end{equation}
Moreover, if $\zeta< h^{-2}$, by the Weyl formula \eqref{weylformula}, 
\begin{equation}\label{discretespectrum3}
\sharp \lbrace j; \ \mu_j(P_h^2) \leq \zeta \rbrace =\sharp \lbrace j; \ \nu_j^s \leq \zeta \rbrace \leq C \langle \zeta \rangle ^{\mathrm{dim}(M)/2s},
\end{equation}
for some $C>0$.
Therefore, by \eqref{discretespectrum2}, we may apply Proposition \ref{propweyl1} to the operator $I+T_h$ to obtain that the spectrum of $I+T_h$ in $[0,h^{-2}/4C_0^2[$ is discrete and all the eigenvalues repeated by multiplicity are found by the min-max principle (see Theorem \ref{minmax}). 

\noindent Thus, again by \eqref{discretespectrum2} and \eqref{discretespectrum3}  
\[
\sharp \lbrace j; \ \mu_j(I+T_h) \leq \zeta \rbrace \leq C_2 \langle \zeta \rangle^{\mathrm{dim}(M)/2s}, \quad \forall \zeta \leq C_1h^{-2},
\]
with $0<C_1<1/(4C_0^2)$ and $C_2>0$ universal constants, independent of $h$.
The statement follows by noting that $\mathrm{Spec}(T_h)=\mathrm{Spec}(I+T_h)-1$.
\end{proof}
This result, for $\mathsf{A}_{h,\rho}$, yields to the next Theorem.
\begin{theorem}
\label{theospec}
There exists $h_0 \in \mathbb{R}$ such that for any $h \in ]0,h_0]$ and any $\rho \in ]0,\rho_0]$ one has that $\mathrm{Spec}(\mathsf{A}_{h,\rho})\cap [0,\zeta]$ is discrete for all $\zeta \leq C'h^{-2}$ and
\begin{equation}
\label{cordiscspec}
\sharp(\mathrm{Spec}(\mathsf{A}_{h,\rho}) \cap [0,\zeta])\leq C'' \langle \zeta \rangle ^{\mathrm{dim}(M)/2}, \quad \forall \zeta \leq C'h^{-2},
\end{equation}
for some constants $C',C''>0$, independent of $h$ and $\rho$. 
\end{theorem}
\begin{proof}
By Proposition \ref{1propshrho} and Proposition \ref{proppositivity} we may apply Lemma \ref{lemspectrum} to the operator $\mathsf{A}_{h,\rho}$.

\noindent By Theorem \ref{theohighlow} there exist $h_0\ll 1$ and a constant $C_0>0$ such that if $u \in L^2(M,d\mu)$ satisfies $((I+\mathsf{A}_{h,\rho})u,u)_0\leq 1$ there exists $(u^L_h,u^H_{h}) \in H^1(M) \times L^2(M,d\mu)$ such that $u=u^L_{h}+u^H_{h}$, with $\norma*{u^L_{h}}_{H^1} \leq C_0\rho$ and $\norma*{u^H_{h}}_0 \leq C_0\rho h$. Then, the fundamental hypothesis of Proposition \ref{propweyl2} is satisfied and hence, taking $C_1'<\frac{1}{4\rho^2C_0^2}$, there exists $C_1''>0$ independent of $h$ such that $\mathrm{Spec}(\mathsf{A}_{h,\rho})\cap [0,\zeta-1]$ is discrete for all $\zeta \leq C_1'h^{-2}$ and 
\[
\sharp (\mathrm{Spec}(\mathsf{A}_{h,\rho}) \cap [0,\zeta-1]) \leq C_1'' \braket{\zeta}^{\mathrm{dim}M/2}, \quad \mathrm{for} \ \mathrm{all} \ \zeta \leq C_1'h^{-2}.
\]
Then, by choosing $h_0$ small enough we get that there exist constants $C',C''>0$, such that  $\mathrm{Spec}(\mathsf{A}_{h,\rho})\cap [0,\zeta]$ is discrete for all $\zeta \leq C'h^{-2}$ and 
\[
\sharp (\mathrm{Spec}(\mathsf{A}_{h,\rho}) \cap [0,\zeta]) \leq C'' \braket{\zeta}^{\mathrm{dim}(M)/2}, \quad \mathrm{for} \ \mathrm{all} \ \zeta \leq C'h^{-2}.
\]
\end{proof}
\begin{remark}
\label{rmkrho}
Note that if we shrink $\rho$, we may choose the constant $C_1'$ and consequently $C'$ in the proof ($C_1'<\frac{1}{4\rho^2C_0^2}$) to be larger. Thus, if we shrink $\rho$, by \eqref{cordiscspec}, we obtain a discrete spectrum in a larger interval. 
\end{remark}
\subsection{The Infinitesimal Generator}
The aim of this subsection is to prove some convergence 
properties of the operator $\mathsf{A}_{h,\rho}$
for $h\rightarrow 0$ and $\rho \in ]0,\rho_0]$ fixed. 
In such a way we will obtain 
\textit{the $\rho$-infinitesimal generator} 
(see Definition \ref{def.infgen}) of our process.
To obtain that, for any $x \in M$ and $f \in C^\infty(M)$,
we are going to calculate the pointwise limit
\[
\lim_{h \rightarrow 0}\mathsf{A}_{h,\rho}f(x):=\lim_{h \rightarrow 0}\frac{I-S_{h,\rho}}{h^2}f(x).
\] 
In what follows, for all $i=1,\dots,n$ we denote by 
\[
\begin{split}
\alpha_{i}: \, &M \rightarrow \lbrace 0,1 \rbrace \\
&x \mapsto \begin{cases}
1, \ \mathrm{if} \ x \in Q_{i} \\
0, \ \mathrm{if} \ x \notin Q_{i}.
\end{cases}
\end{split}
\]
Given $f \in C^\infty(M)$ we consider the map 
\[
f \mapsto 
\begin{pmatrix}
f_{\bigl|_{Q_1}} \\
f_{\bigl|_{Q_2}}\\
\vdots \\
f_{\bigl|_{Q_n}}
\end{pmatrix}
\]
where for each $i=1,\dots,N$, $f_{\bigl|_{Q_i}}=\iota^\ast f$, 
with $ \iota:Q_i\rightarrow M$ the inclusion map. Moreover,
since $X_{i,j} \in \Gamma(Q_i^\ast,TQ_i^\ast)$ and $T_xQ_i \cong T_xQ_i^\ast$ 
for all $x \in Q_i$, we may define ${X_{i,j}}_{\bigl|_{Q_i}} \in \Gamma(Q_i,TQ_i^\ast)$ 
and then the function 
\[
\begin{split}
\alpha_i{X_{i,j}}_{\bigl|_{Q_i}} & f_{\bigl|_{Q_i}}:M \rightarrow \mathbb{R}, \\
&x \mapsto
\begin{cases}
0, \quad \text{if} \ x \notin Q_i \\
{X_{i,j}}_{\bigl|_{Q_i}} f_{\bigl|_{Q_i}}(x), \quad \text{if} \ x \in Q_i.
\end{cases}
\end{split}
\]
Therefore, we define the \textit{truncated vector fields $Z_{i,j}$} as the operators 
\[
 C^\infty(M) \ni f \mapsto Z_{i,j}f:=\alpha_i{X_{i,j}}_{\bigl|_{Q_i}} f_{\bigl|_{Q_i}}.
\]
\begin{remark}
The truncated vector fields are operators from $C^\infty(M)$ to $L^\infty(M)$. 
Indeed, since ${X_{i,j}}_{\bigl|_{Q_i}} f_{\bigl|_{Q_i}}=X_{i,j} f_{\bigl|_{Q_i^\ast}}$ 
for all $x \in Q_i$, we get
\[
\sup_{x \in M} \abs{Z_{i,j}f(x)}=\sup_{x \in Q_i}\abs{X_{i,j} f_{\bigl|_{Q_i^\ast}}(x)}<+\infty.
\] 
\end{remark}
In the same way we define
the operators $P_{i,\rho}:L^2(M,d\mu)\rightarrow L^2(M,d\mu)$ by 
\[
P_{i,\rho}:=-\frac{1}{6n} \sum_{j=1}^n Z_{i,j}^2,
\]
where $Z_{i,j}^2f=\alpha_i{X_{i,j}}^2_{\bigl|_{Q_i}} f_{\bigl|_{Q_i}}$, if $f \in C^\infty(M)$.
Note that $P_{i,\rho}$ depends on $\rho$, since $X_{i,j}$,
$\alpha_i$ and hence $Z_{i,j}$ depends on $\rho$.
\begin{theorem}
\label{semilimit1}
Let $f \in C^\infty(M)$, $\rho \in ]0,\rho_0]$ and $x \in M$. Then
\[
\lim_{h \rightarrow 0}\frac{I-S_{h,\rho}}{h^2}f(x)=\frac{1}{N}\sum_{i=1}^N P_{i,\rho} f(x).
\] 
\end{theorem}
\begin{proof}
Let $x \in M$ and $f \in C^\infty(M)$.
Up to rearranging the cubes,
we may assume that $x \in \bigcap_{i=1}^{N_1} Q_{i}$ 
and $x \notin \bigcup_{i=N_1+1}^NQ_{i}$,
for some $N_1 \in \lbrace 1,...,N \rbrace$.

\noindent Furthermore, note that, choosing $h_0>0$ sufficiently small,
for all $i \in \lbrace 1, \dots N_1 \rbrace$ and for all $h \in ]0,h_0]$ one has  
\[
\Phi_{i,j,h,\rho}(t,x) \in Q_{i}, \quad \forall t \in ]-h,h[.
\]
We write 
\[
(I-S_{h,\rho})f(x)=\frac{1}{nN}(nN-\sum_{i=1}^N\sum_{j=1}^n S_{i,j,h,\rho} )f(x)
\]
and since 
\[
S_{i,j,h,\rho}f(x)=f(x), \quad \forall \ i=N_1+1,\dots,N, \ \forall j=1,...,n,
\]
one has 
\[
\begin{split}
(I-S_{h,\rho})f(x)&=\frac{1}{nN}(nN_1\cdot I-\sum_{i=1}^{N_1}\sum_{j=1}^n S_{i,j,h,\rho})f(x) \\
&=\frac{1}{nN}\sum_{i=1}^{N_1}\sum_{j=1}^n(I- S_{i,j,h,\rho})f(x).\\
\end{split}
\]
Hence, we focus on the term $(I-S_{i,j,h,\rho})f(x)$ 
for some fixed $i=1,...,N_1$, $j=1,...,n$ and we work in coordinates,
with the notation used above.
By Taylor's formula, in coordinates on a neighborhood of $w=\Psi_i(x)$, we have
\[
\begin{split}
\widetilde{S_{i,j,h,\rho}f}(w)=\frac{1}{2h}\int_{-h}^h & \Bigl(\tilde{f}(w)+t\tilde{X}_{i,j}\tilde{f}(w)+\frac{t^2}{2}\tilde{X}_{i,j}^2\tilde{f}(w)+ \\
& +\int_{0}^t \tilde{X}_{i,j}^3\tilde{f}(e^{s\tilde{X}_{i,j}}w)\frac{(t-s)^3}{3!}ds\Bigr)dt.
\end{split}
\]
Thus, noting that 
\[
\int_{-h}^{h}t\tilde{X}_{i,j}\tilde{f}(w)dt=0
\]
and that the last term is $o(t^2)$, for $t\rightarrow 0$, we get
\begin{equation}\label{TaylorSij}
\begin{split}
\widetilde{S_{i,j,h,\rho}f}(w)&=\tilde{f}(w)+\frac{1}{2h}\Bigl[\frac{t^3}{6}\Bigr]_{-h}^h\tilde{X}_{i,j}^2\tilde{f}(w)+o(h^2) \\
&=\tilde{f}(w)+\frac{1}{6}h^2\tilde{X}_{i,j}^2\tilde{f}(w)+o(h^2), \quad h \rightarrow 0.
\end{split}
\end{equation}
Then,
\[
\begin{split}
\lim_{h \rightarrow 0} \frac{I-S_{h,\rho}}{h^2}f(x)&=-\frac{1}{6n N}\sum_{i=1}^{N_1} \sum_{j=1}^n \alpha_i{X_{i,j}}_{\bigl|_{Q_i}}^2 f_{\bigl|_{Q_i}}(x)\\
&=-\frac{1}{6nN}\sum_{i=1}^{N} \sum_{j=1}^n \alpha_{\rho,i} X_{i,j}^2 \alpha_{\rho,i} f(x)\\
&=\frac{1}{N}\sum_{i=1}^N P_{i,\rho}f(x).
\end{split}
\]
\end{proof}
\begin{definition}\label{def.infgen}
We call the operator $\mathsf{A}_{0,\rho}=\frac{1}{N}\sum_{i=1}^N P_{i,\rho}$ just
constructed the \textit{$\rho$-infinitesimal generator} 
of our Markov process.
\end{definition}

Finally, we recall that for $h \in ]0,h_0]$ and $\rho \in ]0,\rho_0]$,
the bilinear form $\mathsf{B}_{h,\rho}$ is defined by
\[
\mathsf{B}_{h,\rho}(u,v)=\Bigl(\frac{I-S_{h,\rho}}{h^2} \,u,v\Bigr)_0, \quad u,v \in L^2(M,d\mu)
\]
and through them, we want to give the following
convergence theorem for $h \rightarrow 0$ and $\rho \in ]0,\rho_0]$ fixed.

\begin{theorem}
\label{TheoconvB}
Let $\varphi \in C^\infty(M)$ and $(f_{h},g_{h})_{h \in ]0,h_0]} \subseteq H^1(M) \times L^2(M,d\mu)$
be a sequence such that 
\[
(i) \ f_{h} \xrightarrow{h \rightarrow 0} f \in H^1(M) \quad \text{weakly}, 
\]
and
\[
(ii) \ \sup_{h \in ]0,h_0]} \norma*{g_{h}}_0 < +\infty.
\]
Then
\[
\lim_{h \rightarrow 0} \mathsf{B}_{h,\rho}(f_h+hg_h,\varphi)=\frac{1}{6n N} \sum_{i=1}^N \sum_{j=1}^n(X_{i,j}f,X_{i,j}\varphi)_{L^2(Q_i)}.
\]
\end{theorem}
\begin{remark}
First to start the proof let us point out that
for $u,v\in L^2(M,d\mu)$ one has
\begin{equation}\label{BsumBij}
\mathsf{B}_{h,\rho}(u,v)_0=\frac{1}{nN} \sum_{i=1}^N\sum_{j=1}^n\mathsf{B}_{i,j,h,\rho}(u,v).
\end{equation}
Moreover, if $u,v \in C^\infty(M)$ (and by density also for $u,v \in L^2(M,d\mu)$)
\begin{equation}\label{uv.dopprod}
\begin{split}
&\mathsf{B}_{i,j,h,\rho}(u,v)\\
&=\Bigl(\frac{I-S_{i,j,h,\rho}}{h^2}u,v\Bigr)_0\\
&=h^{-2}\int_M\bigl((I-(M_{i,j,h,\rho}+T_{i,j,h,\rho}))u\bigr) \, \overline{v} \, d\mu\\
&=\frac{1}{2h^3} \int_{Q_i} \int_{-h}^h \chi_{i,j,h,\rho}(t,x)(u(x)-u(e^{tX_{i,j}}x))\overline{v(x)} \, dtd\mu \\
&=\frac{1}{4h^3}\int_{Q_i} \int_{-h}^h \chi_{i,j,h,\rho}(t,x)(u(x)-u(e^{tX_{i,j}}x))\overline{(v(x)-v(e^{tX_{i,j}}x))} \, dtd\mu.
\end{split}
\end{equation}
\end{remark}
\begin{proof}[Proof of the Theorem \ref{TheoconvB}]
We write
\[
\mathsf{B}_{h,\rho}(f_h+hg_h,\varphi)=\mathsf{B}_{h,\rho}(f_h,\varphi)+\mathsf{B}_{h,\rho}(hg_h,\varphi)
\]
and then we have to prove
\begin{equation}
\label{BihrhotoXij1}
\lim_{h \rightarrow 0} \mathsf{B}_{h,\rho}(f_h,\varphi)=\frac{1}{6n N} \sum_{i=1}^N \sum_{j=1}^n(X_{i,j}f,X_{i,j}\varphi)_{L^2(Q_i)}
\end{equation}
and
\begin{equation}
\label{BihrhotoXij2}
\lim_{h \rightarrow 0}\mathsf{B}_{h,\rho}(hg_h,\varphi)=0.
\end{equation}
We first examine \eqref{BihrhotoXij2}. By 
\eqref{BsumBij} we may reduce to show that for any given $i,j$ one has 
\begin{equation}\label{BSij}
\lim_{h \rightarrow 0}\Bigl(hg_h, \frac{I-S_{i,j,h,\rho}}{h^2}\varphi\Bigl)_0=0.
\end{equation}
Since $\varphi \in C^\infty(M)$, for $x \in Q_i$ we have 
\[
\begin{split}
(I-S_{i,j,h,\rho})\varphi(x)&= \frac{1}{2h}\int_{-h}^h \chi_{i,j,h,\rho}(t,x)(\varphi(x)-\varphi(e^{tX_{i,j}}x))dt\\
&=\frac{X_{i,j}\varphi(x)}{2h} \int_{-h}^h t \, \chi_{i,j,h,\rho}(t,x)dt +O_{\mathrm{unif}}(h^2)\\
&=:r_h(x)+O_{\mathrm{unif}}(h^2).
\end{split}
\]
Hence, noting that $\mathrm{supp} \, r_h \subseteq \lbrace x \in Q_i; d(x,\partial Q_i)<h \rbrace$
and $\|r_h\|_{L^\infty}=O_{\mathrm{unif}}(h)$ we get $\|r_h\|_{0}=O_{\mathrm{unif}}(h^{3/2})$. 
Therefore, since $\|hg_h\|_0=O_{\mathrm{unif}}(h)$ we obtain
\[
\begin{split}
B_{i,j,h,\rho}(hg_h,\varphi)&=\frac{1}{h^2}(hg_h,(I-S_{i,j,h,\rho})\varphi)_0\\
&=(g_h, h^{-1}(I-S_{i,j,h,\rho})\varphi)_0\\
&=(g_h,h^{-1}r_h)_0+O_{\mathrm{unif}}(h) \\
&=O_{\mathrm{unif}}(h^{1/2})
\end{split}
\] 
and thus \eqref{BihrhotoXij2} follows from \eqref{BsumBij}.

\noindent We now prove \eqref{BihrhotoXij1}. Also in this case we may reduce to study
\[
\mathsf{B}_{i,j,h,\rho}(f_h,\varphi)=\Bigl(\frac{I-S_{i,j,h,\rho}}{h^2}f_h,\varphi\Bigr)_0.
\]
In the first place, by \eqref{uv.dopprod},
\[
\begin{split}
&\Bigl(\frac{I-S_{i,j,h,\rho}}{h^2}f_h,\varphi\Bigr)_0= \\
&=\frac{1}{4h^3}\int_{Q_i}\Bigl(\int_{-h}^h \chi_{i,j,h,\rho}(t,x)(f_h(x)-f_h(e^{tX_{i,j}}x))\overline{(\varphi(x)-\varphi(e^{tX_{i,j}}x))}dt\Bigr)d\mu(x).
\end{split}
\]
At this point, writing
\[
\begin{split}
&\varphi(x)-\varphi(e^{tX_{i,j}}x)=t\psi(t,x)\\
&f_h(x)-f_h(e^{tX_{i,j}}x)=t \int_0^1X_{i,j}f_h(e^{tsX_{i,j}}x)ds,
\end{split}
\]
with $\psi(t,x)$ smooth function such that $\psi(0,x)=X_{i,j}\varphi(x)$,
and denoting $\Phi_t(x)=e^{tX_{i,j}}x$, one has
\[
\begin{split}
&\mathsf{B}_{i,j,h,\rho}(f,\varphi)\\
&=\frac{1}{4h^3}\int_{Q_i}\Bigl(\int_{-h}^h \chi_{i,j,h,\rho}(t,x)t^2\bigl(\int_0^1X_{i,j}f_h(e^{tsX_{i,j}}x)ds\bigr)\overline{\psi(t,x)}dt\Bigr)d\mu(x)\\
&=\frac{1}{4}\int_{Q_i}\Bigl(\int_{-1}^1 \chi_{i,j,h,\rho}(h\tau,x)\tau^2\bigl(\int_0^1X_{i,j}f_h(e^{h\tau sX_{i,j}}x)ds\bigr)\overline{\psi(h\tau,x)}d\tau\Bigr)d\mu(x)\\
&=\frac{1}{4}\int_{-1}^1 \int_0^1\int_{\Phi_{h\tau s}(Q_i)}\Bigl( \chi_{i,j,h,\rho}(h\tau,e^{-h\tau sX_{i,j}}y)\tau^2 \times \\
&\hspace{5cm}\times X_{i,j}f_h(y)\overline{\psi(h\tau,e^{-h\tau sX_{i,j}}y)}\Bigr)d\mu(y)dsd\tau.
\end{split}
\]
Next, since $\psi(hr,e^{-h \tau sX_{i,j}}y)=X_{i,j}\varphi(y)+O_{\mathrm{unif}}(h)$, 
for any $\delta>0$ 
we get 
\[
\begin{split}
&\mathsf{B}_{i,j,h,\rho}(f_h,\varphi)=
\\
&=\frac{1}{4}\int_{-1}^1 \int_0^1\int_{\Phi_{h\tau s}(Q_i)\cap Q_i} \hspace{-1.5cm}\chi_{i,j,h,\rho}(h\tau,e^{-h\tau sX_{i,j}}y)\tau^2  X_{i,j}f_h(y)\overline{X_{i,j}\varphi(y)}d\mu(y)dsd\tau + O_{\mathrm{unif}}(h^2) \\
&=\frac{1}{4}\int_{-1}^1 \int_0^1\int_{\Phi_{h\tau s}(Q_i)\cap\lbrace y \in Q_i; \; d(y,\partial Q_i)\geq\delta \rbrace} \hspace{-3.5cm}\chi_{i,j,h,\rho}(h\tau,e^{-h\tau sX_{i,j}}y)\tau^2  X_{i,j}f_h(y)\overline{X_{i,j}\varphi(y)}d\mu(y)dsd\tau\\
&\quad + \frac{1}{4}\int_{-1}^1 \int_0^1\int_{\Phi_{h\tau s}(Q_i)\cap\lbrace y\in Q_i; \; d(y,\partial Q_i)<\delta \rbrace} \hspace{-3.5cm}\chi_{i,j,h,\rho}(h\tau,e^{-h\tau sX_{i,j}}y)\tau^2  X_{i,j}f_h(y)\overline{X_{i,j}\varphi(y)}d\mu(y)dsd\tau \\
&\quad +O_{\mathrm{unif}}(h^2) \\
&=:I_\delta(h)+J_\delta(h)+O_{\mathrm{unif}}(h^2).
\end{split}
\]
Now, by Cauchy-Schwarz inequality 
\[
\abs{J_\delta(h)} \leq C(\varphi)\delta^{1/2}\| f_h \|_{H^1(Q_i)}.
\] 
and, furthermore, choosing $h_0=h_0(\delta)$ sufficiently small, we have that for all $h \in ]0,h_0]$
\[
\begin{split}
I_\delta(h)&=\frac{1}{6}\int_{Q_i \cap \lbrace x; \; d(x,\partial Q_i)>\delta\rbrace}X_{i,j}f_h(x)\overline{X_{i,j}\varphi(x)}d\mu(x) \\
&=\frac{1}{6} \int_{Q_i}X_{i,j}f_h(x)\overline{X_{i,j}\varphi(x)}d\mu(x)+O_{\mathrm{unif}}(\delta_{1/2}\|f_h\|_{H^1(Q_i)}).
\end{split}
\]
Therefore, given $\varepsilon>0$ we may find $\delta>0$ 
sufficiently small such that for any $h \in ]0,h_0(\delta)]$
\[
\abs{J_\delta(h)}<\varepsilon
\]
and 
\[
\abs{I_\delta(h)-\frac{1}{6} \int_{Q_i}X_{i,j}f_h(x)\overline{X_{i,j}\varphi(x)}d\mu(x)}<\varepsilon.
\] 
Hence, since $f_h \rightarrow f$ weakly in $H^1(Q_i)$,
for $h \rightarrow 0$, we obtain  
\[
\lim_{h \rightarrow 0}\mathsf{B}_{i,j,h,\rho}(f_h,\varphi)=\frac{1}{6} \int_{Q_i}X_{i,j}f(x)\overline{X_{i,j}\varphi(x)}d\mu(x).
\]
\end{proof}
\section{The Spectral Result and Convergence to Equilibrium}\label{sec.maintheo}
The aim of this section is to prove Theorem \ref{Theorem1} and Theorem \ref{Theorem2}.
The first one concerns the spectrum
of the operator $S_{h,\rho}$
for $h \in ]0,h_0]$ and $\rho \in ]0,\rho_0]$, while
Theorem \ref{Theorem2}
establishes the convergence of our $(h,\rho)$-subelliptic random walk
to the equilibrium when the number of steps goes to infinity. 

\subsection{\texorpdfstring{Spectrum Properties of the {$(h,\rho)$}-Markov Operator}{Spectrum Properties of the {(h,rho)}-Markov Operator}}
We start by proving the spectral results of Theorem \ref{Theorem1}.
More precisely, we will see some spectral properties
near the eigenvalue $1$ and we will investigate what happens
if we shrink the parameter $\rho$. Finally, we will show that
the spectral gap $g_{h,\rho}$ of the operator $I-S_{h,\rho}$ (see Definition \ref{def.specgap}),
that is the distance between the
eigenvalue $1$ for $S_{h,\rho}$ (cf. Remark \ref{rmkspecgap}), 
satisfies $g_{h,\rho}\approx h^2$.
\begin{proof}[Proof of Theorem \ref{Theorem1}]
In the first place note that for all $k \in \mathbb{N}$
\[
2((I-S_{h,\rho}^k)u,u)_0=\int_{M}\Bigl(\int_M\abs{u(x)-u(y)}^2s_{h,\rho}^k(x,dy)\Bigr)d\mu(x).
\]
In fact, to prove this equality we can reduce to prove
\[
\int_M\Bigl(\int_M \abs{u(y)}^2s_{h,\rho}^k(x,dy)\Bigr)d\mu(x)=\int_M \abs{u(x)}^2 d\mu(x),
\]
but this follows by reasoning as in \eqref{usefulpositivity},
using the change of variables and the \textit{divergence-free} property
of the vector fields $X_{i,j}$ with respect the local expression
of the measure $\mu$. 

\noindent We now prove $(i)$. We already saw in 
Proposition \ref{1propshrho} that $S_{h,\rho}1=1$,
thus it remains to show that $1$ is a simple eigenvalue of $S_{h,\rho}$. 
Let $u \in L^2(M,d\mu)$ be a function that satisfies $S_{h,\rho}u=u$.
By Theorem \ref{lowestimatekernel} there exists a universal $\delta>0$ such that 
\begin{equation}\label{1simpleeigenvalue}
\begin{split}
& \int_{M}\Bigl(\int_M\abs{u(x)-u(y)}^2s_{h,\rho}^n(x,dy)\Bigr)d\mu(x) \\
&\geq ch^{-n}\int_{M}\Bigl(\int_{\lbrace y; \; d_g(x,y) <\delta h \rbrace}\abs{u(x)-u(y)}^2d_gy\Bigr)d\mu(x).
\end{split}
\end{equation}
Hence, since 
\[
\int_{M}\Bigl(\int_M\abs{u(x)-u(y)}^2s_{h,\rho}^n(x,dy)\Bigr)d\mu(x)=2((I-S_{h,\rho}^n)u,u)_0=0,
\]
by \eqref{1simpleeigenvalue} we get that
$u(x)-u(y)=0$ for almost all $(x,y)$ such that $d_g(x,y)<\delta h$.
Therefore, since $M$ is connected, we obtain that $u$ is constant. 

\noindent We proceed to prove $(ii)$, i.e. we prove that the spectrum
of $S_{h,\rho}$ is a subset of $[-1+\delta_1,1]$.
To do that it is sufficent to prove the statement for an odd power of the operator,
hence we can reduce to prove that there exist $h_0,c>0$
such that for all $h \in ]0,h_0]$,
\[
((I+S_{h,\rho}^{2k+1})u,u)_0 \geq c \norma*{u}_0^2 
\]
and thus 
\[
\int_{M}\Bigl(\int_M\abs{u(x)+u(y)}^2s_{h,\rho}^{2k+1}(x,dy)\Bigr)d\mu(x)\geq c \norma*{u}_0^2.
\]
Again by Theorem \ref{lowestimatekernel} (see also Remark \ref{rmkpowerk})
it is sufficient to prove that for all $h \in ]0,h_0]$, 
\[
h^{-n}\int_{M}\Bigl(\int_{\lbrace y; \; d_g(x,y) <\delta h \rbrace}\abs{u(x)+u(y)}^2d_gy\Bigr)d\mu(x)\geq c \norma*{u}_0^2
\]
for some constant $c>0$ independent of $h$. 

\noindent By compactness we can cover $M$ by a finite number
of open set $(U_j)_j$, such that $\mathrm{diam} (U_j) <h$, $\mu(U_j)>c_1h^n$ 
and for all $j$ the number of $k$ such that $U_k \cap U_j \neq 0$ is bounded from $c_2$,
with $c_1$ and $c_2$ universal constants (that is they do not depend on $h$).   
In fact, to do so it is sufficient to use Lemma 1.4.9 of \cite{HoV1} on the universal block 
$Q_\rho \subseteq \mathbb{R}^n$ and use that $M$ is a \emph{finite} 
union of $Q_i=\Psi_i^{-1}(Q_\rho)$.  

\noindent In this way, by reasoning as in \cite{DLM1} (p. 246), 
\[
\begin{split}
&c_2h^{-n}\int_M\Bigl(\int_{\lbrace y; \; d_g(x,y) <\delta h \rbrace}\abs{u(x)+u(y)}^2d_gy\Bigr)d\mu(x) \\
&\geq h^{-n} \sum_{j=1} \int_{U_j}\Bigl(\int_{U_j} \mathbbm{1}_{\lbrace y; \; d_g(x,y) <\delta h \rbrace}(y)\abs{u(x)+u(y)}^2d_gy\Bigr)d\mu(x) \\
&\geq ch^{-n} \sum_{j=1} \int_{U_j}\Bigl(\int_{U_j}\abs{u(x)+u(y)}^2d\mu(y)\Bigr)d\mu(x) \\
&\geq ch^{-n} \sum_{j=1} \mu(U_j)\int_{U_j}\abs{u(x)}^2d\mu(x)\\
&\geq c \norma*{u}_0^2,
\end{split}
\]
for some constant $c>0$. In conclusion 
\[
((1+S_{h,\rho}^{2n+1})u,u)_0 \geq c \norma*{u}_0^2, 
\] 
and then there exists $\delta_1>0$ constant, independent of $h$,
such that 
\[
\mathrm{Spec}(S_{h,\rho}) \subseteq [-1+\delta_1,1].
\] 
To obtain $(iii)$ it is sufficient to note that if $\zeta>0$, 
\[
\begin{split}
&\mathsf{A}_{h,\rho}u=\frac{I-S_{h,\rho}}{h^2}u=\lambda' u, \ \text{with} \ \lambda' \in [0,\zeta] \\
& \iff S_{h,\rho}u=\lambda u, \text{with} \ \lambda \in [1-h^2\zeta,1], \ \lambda=1-h^2\lambda'.
\end{split}
\]
Thus, by Corollary \ref{theospec} we conclude
that there exists $\delta_2>0$ such that
$\mathrm{Spec}(S_{h,\rho}) \cap [1-\delta_2,1]$ is discrete and
\begin{equation}
\label{estspec}
\sharp\Bigl(\mathrm{Spec}(S_{h,\rho})\cap [1-h^2\zeta,1]\Bigr)\leq C_1(1+\zeta)^m, \quad\text{for any} \ 0 \leq \zeta \leq \delta_2h^{-2},
\end{equation}
for some $m>0$.

\noindent It remains to show $(iv)$, i.e. that there exist $C_2,C_3>0$, 
independent of $h$, such that the spectral gap $g_{h,\rho}$ satisfies
\begin{equation}\label{boundsgh}
C_3h^2\leq g_{h,\rho} \leq C_2h^2, \quad \forall h \in ]0,h_0].
\end{equation}
Let us recall that if $f\in C^\infty(M)$ and $x \in M$ by Theorem \ref{semilimit1} one has 
\[
\lim_{h \rightarrow 0} \frac{I-S_{h,\rho}}{h^2}f(x)=\mathsf{A}_{0,\rho} f(x)
\]
and then the upper bound follows by the min-max principle. 

\noindent Regarding the lower bound, assume by contradiction that
there exist two sequences $(\varepsilon_k)_k$, $(h_k)_k$,
such that $\varepsilon_k,h_k \rightarrow 0$, $k\rightarrow +\infty$,
and a sequence $(u_k)_k\subseteq L^2(M,d\mu)$,
with $\norma*{u_k}_0=1$ and $(u_k,1)_{0}=\int_{M}u_k d\mu=0$ such that 
\[
\frac{I-S_{h_k,\rho}}{h_k^2}u_k=\varepsilon_k u_k.
\]
This implies, by definition of the Dirichlet form $\mathsf{E}_{h,\rho}$
(see \eqref{Dirichletform}), that $\mathsf{E}_{h_k,\rho}(u_k)=\varepsilon_k$ 
and so, by Theorem \ref{theohighlow}, one has
\[
u_k=u^L_{k}+u^H_{k}=u^L_{k}+h_k\frac{u^H_{k}}{h_k}=:f_{k}+h_kg_{k}, 
\]
with
\[
\sup_k \, \norma*{g_k}_0<C_0, \quad \sup_k \, \norma*{f_{k}}_{1}\leq C_0.
\]
At this point, since the sequence $(f_{k})_k$ is bounded in $H^1(M)$,
we have that there exists $f \in H^1(M)$ such that $f_{k}\rightarrow f$
weakly on $H^1(M)$ as $k \rightarrow + \infty$. Moreover,
up to subsequences we have that $f_k \rightarrow f$ strongly in $L^2(M,d\mu)$ 
and hence $\|f\|_0=1$ and $\int_M f d\mu=0$.
Therefore, by Theorem \ref{TheoconvB}, we have that for all $\varphi \in C^\infty(M)$ 
\[
0=\lim_{k +\infty} \mathsf{B}_{i,j,h,\rho} (u_k,\varphi)=\frac{1}{6n N} \sum_{i=1}^N \sum_{j=1}^n(X_{i,j}f,X_{i,j}\varphi)_{L^2(Q_i)}.
\] 
Finally, since $C^\infty(M)\subseteq H^1(M)$ is dense,
we may choose a sequence
$(\varphi_{k'})_{k'} \in C^\infty(M)$ 
such that $\varphi_{k'} \rightarrow f$, strongly in $H^1(M)$, 
for ${k} \rightarrow +\infty$.

\noindent Hence $\|X_{i,j}f \|_{L^2(Q_i)}=0$ for all $i=1,\dots,N$,
and since for all $i$ the collection $\lbrace X_{i,j} \rbrace_j$
forms a non-degenerate system of vector fields, we may conclude that 
$f$ is constant on $M$. This contradicts the fact $\|f\|_0=1$ and $\int_M f d\mu=0$
and proves the lower bound \eqref{boundsgh}.
This concludes the proof of the Theorem.
\end{proof}
\begin{remark}
    Note that by Remark \ref{rmkrho} the constant $\delta_2$ increases by shrinking $\rho$ and therefore we have a discrete spectrum on a larger interval near $1$. 
\end{remark}

\subsection{Convergence to Equilibrium}
We now prove the convergence to equilibrium stated in Theorem \ref{Theorem2}.
The convergence will be measured by the so-called \textit{total variation distance} defined as follows.
\begin{definition}\label{def.TV}
Let $(X,\mathscr{F})$ be a measurable space, and $\mu, \nu$ 
be two measures on $(X,\mathscr{F})$. \textit{The total variation distance} between $\mu$ and $\nu$
is defined as
\[
\norma*{\mu-\nu}_{TV}=\sup_{A \in \mathscr{F}} \abs*{\mu(A)-\nu(A)}.
\] 
\end{definition} 
Note that the velocity of such convergence 
will depend on the spectral gap $g_{h,\rho}$ of the $(h,\rho)$-Markov operator. 

The estimate \eqref{convequi} encodes a convergence to equilibrium.
In fact, since the local vector fields $X_{i,j}$
through which we have constructed our diffusion
are divergence free with respect to $d\mu$,
by the self-adjoint property of the operator,
it is possible to obtain the reversibility 
of the Markov Chain and 
then the stationarity for the measure associated
with respect the constructed diffusion (see \cite{No}).
\begin{proof}[Proof of Theorem \ref{Theorem2}]
    We first note that for $k \in \mathbb{N}$,
\[
\sup_{x \in M} \norma*{s_{h,\rho}^k(x,dy)-\mu}_{TV}=\frac{1}{2}\norma*{S_{h,\rho}^k-\Pi_0}_{L^\infty\rightarrow L^\infty},
\]
where $\Pi_0$ is the orthogonal projection of $L^2(M,d\mu)$ onto the space of 
constant functions on $M$
\[
\begin{split}
\Pi_0:L^2(M, & d\mu)\rightarrow V_1=\mathrm{Ker}(S_{h,\rho}-1\cdot I) \\
&u \mapsto \int_M u \, d\mu.
\end{split}
\]
We thus have to prove that there exists $h_0>0$ such that,
for all $h\in ]0,h_0]$ and for all $k \in \mathbb{N}$ 
\begin{equation}\label{tve}
\norma*{S_{h,\rho}^k-\Pi_0}_{L^\infty\rightarrow L^\infty} \leq C e^{-kg_{h,\rho}}
\end{equation}
for some $C>0$ constant independent of $h$.
Since by Theorem \ref{Theorem1}-$(iv)$ we have $g_{h,\rho}\approx h^2$ and 
\[
\norma*{S_{h,\rho}^k-\Pi_0}_{L^\infty \rightarrow L^\infty} \leq 2,
\]
we may assume that $k \geq C_0h^{-2}$ for some $C_0>0$ sufficiently large.
To prove \eqref{tve} we denote by
\[
0<\lambda_{1,h}<\dots <\lambda_{j,h}< \dots \leq \delta_2h^{-2}
\]
the eigenvalues of the operator $\mathsf{A}_{h,\rho}$ smaller
than $\delta_2h^{-2}$, with $\delta_2>0$ sufficiently small to be chosen. 
Moreover, we denote by $J_h$ the set of indexes 
\[
J_h:= \lbrace j \in \mathbb{N}; \; \lambda_{j,h} \leq \delta_2h^{-2} \rbrace, 
\]
and by
\[
V_h:=\mathrm{Span} \Bigl(\Bigl\lbrace e_{j,h}\in L^2(M,d\mu); \; \mathsf{A}_{h,\rho}e_{j,h}=\lambda_{j,h}e_{j,h}, \; j \in J_h \Bigr\rbrace\Bigr)\subseteq L^2(M,d\mu). 
\]
Recall that, by Theorem \ref{Theorem1}-$(iii)$, $\mathrm{dim} (V_h) \leq C h^{-\mathrm{dim}(M)}$, for some universal constant $C>0$. 
Moreover, we have
\[
S_{h,\rho}e_{j,h}=(1-\lambda_{j,h} h^{-2})e_{j,h}, 
\]
for all $j \in J_h$ and $(e_{j,h} , e_{k,h})_0=\delta_{jk}$ for all $j,k \in J_h$.

\noindent Write now  
\[
S_{h,\rho}-\Pi_0=K_{1,h}+K_{2,h}
\]
where $K_{1,h}:L^2(M,d\mu)\rightarrow L^2(M,d\mu)$ is the operator associated 
with the Schwartz kernel
\[
k_{1,h}(x,y):=\sum_{j \in J_h}(1-h^2\lambda_{j,h})e_{j,h}(x)\overline{e_{j,h}(y)}.
\]
Since $S_{h,\rho}^k-\Pi_0=K_{1,h}^k+K_{2,h}^k$, for all $k\in N$,
in order to prove \eqref{tve}, we may reduce to give
the $L^\infty$ estimates for the operators $K^k_{1,h}$ and $K^k_{2,h}$.
We begin with the following rough bounds. 
\begin{lemma}
Let $u \in V_h$. Then, there exists a constant $C'>0$ such that
\[
\norma*{u}_{L^\infty}\leq C'h^{-n/2-\mathrm{dim}(M)/2}\norma*{u}_0. 
\]
In addition, for all $k \in \mathbb{N}$,
\[
\norma*{K_{1,h}^k}_{L^\infty \rightarrow L^\infty}\leq C'h^{-n/2-\mathrm{dim}(M)/2}
\]
and
\[
\norma*{K_{2,h}^k}_{L^\infty \rightarrow L^\infty} \leq C'h^{-n/2-\mathrm{dim}(M)/2}.
\]
\end{lemma}
\begin{proof}
As usual, in what follows, by $C'>0$ we mean a universal constant
possibly changing from line to line. By
definition $u= \sum_{j \in J_h} \hat{u}_{j,h}e_{j,h}$.
From Corollary \ref{corlowest2}, possibly shrinking $\delta_2>0$,
for all $j \in J_h$ one has 
\[
\norma*{e_{j,h}}_{L^\infty} \leq C'h^{-n/2}
\]
Therefore, by the Cauchy-Schwarz inequality and Theorem \ref{theospec} 
\begin{equation}\label{eq.CS}
\norma*{u}_{L^\infty} \leq C'h^{-n/2} \Bigl(\sum_{j \in J_h}\abs*{\hat{u}_{j,h}}^2\Bigr)^{1/2}(\mathrm{dim}(V_h))^{1/2} \leq C' h^{-n/2-\mathrm{dim}(M)/2}\|u\|_0.
\end{equation}
Take now $u \in L^2(M,d\mu)$. By definition of $K_{1,h}$ we have
\[
K_{1,h}u(x)=\int_M u(y)k_{1,h}(x,y)d\mu(y)=\sum_{j \in J_h}(1-h^2\lambda_{j,h})\hat{u}_{j,h}e_{j,h}(x).
\]
Hence, since $\abs*{(1-h^2\lambda_{j,h})}\leq 1$, by reasoning as above, we obtain 
\[
\norma*{K_{1,h}^k}_{L^\infty \rightarrow L^\infty} \leq \norma*{K_{1,h}^k}_{L^2 \rightarrow L^\infty} \leq C'h^{-n/2-\mathrm{dim}(M)/2},
\]
Finally, from $S_{h,\rho}^k-\Pi_0=K_{1,h}^k+K_{2,h}^k$ and the fact
that the norm $S_{h,\rho}^k$ is bounded by $1$ from $L^\infty(M,d\mu)$
to $L^\infty(M,d\mu)$ we conclude that 
\begin{equation}\label{estK21}
\norma*{K_{2,h}^k}_{L^\infty \rightarrow L^\infty} \leq C'h^{-n/2-\mathrm{dim}(M)/2}.
\end{equation}
\end{proof}
\noindent Let now $r_h(x,dy)$ and $ch^{-n}\mathbbm{1}_{\lbrace d_g(x,y)<\varepsilon h \rbrace}d_gy$
be the Markov Kernels introduced in Theorem \ref{lowestimatekernel}.
For simplicity in what follows we denote,
respectively, by $R_h$ and $L_h$ the corresponding operators.
By Theorem \ref{lowestimatekernel} one has that $S_{h,\rho}^n=R_h+L_h$
and by Corollary \ref{corlowest1} 
\[
\norma*{R_h}_{L^\infty \rightarrow L^\infty}\leq \tau,
\]
for some $0<\tau <1$. Moreover, by definition of $L_h$,
there exists a constant $c>0$ such that
\[
\norma*{L_h}_{L^2 \rightarrow L^\infty} \leq c h^{-n/2}.
\]
For all $k \geq C_0h^{-2}$ we now write $k=qn+\nu$ for
some $q \in \mathbb{N}$ and $0 \leq \nu < n$.
Thus, for all $q \in \mathbb{N}$, 
\[
S_{h,\rho}^{qn}=F_{q,h}+G_{q,h} 
\]
where $F_{q,h}$ and $G_{q,h}$ are defined by recurrence as  \\
$\bullet$ $F_{1,h}:=R_h$, \quad $G_{1,h}:=L_h$; \\
$\bullet$ $F_{q+1,h}:=R_hF_{q,h}$, \quad $G_{q+1,h}:=R_h G_{q,h}+L_hS_{h,\rho}^{qn}$.
  
\noindent In this way, we obtain 
\begin{equation}\label{estK22}
\norma*{F_{q,h}}_{L^\infty \rightarrow L^\infty} \leq \tau^q
\end{equation}
and, since $S_{h,\rho}^{qn}$ is bounded by $1$ on $L^2$,
\begin{equation}\label{estK23}
\norma*{G_{q,h}}_{L^2 \rightarrow L^\infty} \leq ch^{-n/2}(1+\tau+\dots+\tau^q) \leq c h^{-n/2}/(1-\tau).
\end{equation}
Furthermore, since the operator $K_{2,h}$ is bounded
and self-adjoint from $L^2(M,d\mu)$ onto itself,
by Theorem \ref{theospecradius},
we have $\norma*{K_{2,h}}_{L^2 \rightarrow L^2}=r_{\mathrm{spec}}(K_{h,2})=:\theta <1$.
Therefore, for all $k \in \mathbb{N}$,
\begin{equation}\label{estK24}
\norma*{K_{2,h}^k}_{L^\infty \rightarrow L^2} \leq \norma*{K_{2,h}^k}_{L^2 \rightarrow L^2} \leq \theta^k.
\end{equation}
Thus, for $m \geq 1$, $q \geq 1$ and $0 \leq \nu < n-1$, putting together 
the estimates \eqref{estK21}, \eqref{estK22}, \eqref{estK23}, \eqref{estK24} and
the fact that the norm $S_{h,\rho}$ is bounded by $1$ from $L^\infty(M,d\mu)$ to $L^\infty(M,d\mu)$, we get
\[
\begin{split}
\norma*{K_{2,h}^{qn+\nu+m}}_{L^\infty \rightarrow L^\infty} & \leq \norma*{S_{h,\rho}^\nu S_{h,\rho}^{qn}K_{2,h}^m}_{L^\infty \rightarrow L^\infty} \\
&\leq \norma*{F_{q,h}K_{2,h}^m}_{L^\infty \rightarrow L^\infty}+\norma*{G_{q,h}K_{2,h}^m}_{L^\infty \rightarrow L^\infty} \\
&\leq Ch^{-n/2-\mathrm{dim}(M)/2}\tau^q +ch^{-n/2} \theta^m/(1-\tau).
\end{split}
\]
Thus, there exists $C_1>0$, $b>0$ and $B\gg 1$ such that for all $h \in ]0,h_0]$
\begin{equation}\label{estk2h}
\norma*{K_{2,h}^k}_{L^\infty \rightarrow L^\infty} \leq C_1e^{-bk},
\end{equation}
for all $k \geq B \log(1/h)$ and hence,
to obtain estimate \eqref{tve}, the contribution
of $K_{2,h}$ is far smaller than $e^{-kg_{h,\rho}}$. 
We now focus on the contribution 
of the first operator $K_{1,h}$
and to do that we first recall the following Lemma (cf. \cite{LM}, Lemma 5.3).

\begin{lemma}
There exist $p>2$ and $c'>0$ constant, independent of $h$,
such that for all $u \in V_h$, one has 
\[
\norma*{u}_{L^p}^2 \leq c'(\mathsf{E}_{h,\rho}(u)+\norma*{u}^2_0).
\]
\end{lemma}
\begin{proof}
Once more, we denote by $c'>0$ an absolute constant possibly
changing from line to line. Let $u \in V_h$. 
By reasoning similarly as in the proof of Proposition \ref{propTh} we assume
\begin{equation}\label{eq.lemmaLp}
\mathsf{E}_{h,\rho}(u)+\norma*{u}^2_0 \leq 1.
\end{equation}
Hence, from Theorem \ref{theohighlow}, we have $u=u^L_{h}+u^H_{h}$
with $\norma*{u^L_{h}}_{H^1(M,d\mu)}\leq C$ and $\norma*{u^H_{h}}_0 \leq c'h$.
Thus, by the continuous embedding $H^s(M)\subseteq L^q(M)$ with $s>0$,
$s=n(1/2-1/q)$ we get 
\[
\norma*{u^L_{h}}_{L^q} \leq c', 
\]
with $q=\frac{2n}{n-1}$. 
Moreover, if $u=\sum_{\lambda_{j,h}\leq \delta_2h^{-2}}\hat{u}_{j,h}e_{j,h}$
by our hypothesis \eqref{eq.lemmaLp}
we have 
\[
\sum_{\lambda_{j,h}\leq \delta_2h^{-2}}\abs*{\hat{u}_{j,h}}^2\leq1
\]
and by Corollary \ref{coreigen2},
if $\delta_2>0$ is small enough $\norma*{e_{j,h}}_{L^\infty}\leq c'h^{-n/2}$ for all $j$.
Therefore, by repeating the same argument of \eqref{eq.CS}, we obtain 
\[
\norma*{u}_{L^\infty}\leq c'h^{-n/2-\mathrm{dim}(M)/2}.
\]
Next, since $u^L_{h}$ is defined by \eqref{def.u1hu2h} as sum
of local terms of the form \eqref{def.Ah}, by Young's inequality we have 
\[
\norma*{u^L_h}_{L^\infty}\leq c'\norma*{u}_{L^\infty}.
\]
Thus,
\[
\norma*{u^H_h}_{L^\infty}\leq \norma*{u^L_h}_{L^\infty}+\norma*{u}_{L^\infty}\leq C h^{-n/2-\mathrm{dim}(M)/2}.
\]
Furthermore, since $\norma*{u^H_{h}}_0 \leq c'h$ we obtain 
by interpolation that there exists $q'>2$ such that 
\[
\norma*{u^H_{h}}_{L^{q'}}\leq c'.
\] 
Therefore, the result follows by choosing $p=\min \lbrace q,q' \rbrace$.
\end{proof}

\noindent As a consequence of this lemma, using the interpolation inequality 
\[
\norma*{u}_0^2 \leq \norma*{u}_{L^p}^{\frac{p}{p-1}}\norma*{u}_{L^1}^{\frac{p-2}{p-1}},
\]
we deduce that, for all $u \in V_{h}$,
\begin{equation}\label{Nashinequality}
\norma*{u}_0^{2+1/d} \leq C (\mathsf{E}_h(u)+\norma*{u}_0^2)\norma*{u}^{1/d}_{L^1},
\end{equation}
where $d:=1/(2-4/p)>0$.  Moreover, for $\lambda_{j,h} \leq \delta_2h^{-2}$ one has $h^2\lambda_{j,h} \leq 1$ and hence, for $u \in V_h$, $u=\sum_{j \in J_h}\hat{u}_je_j$ we get
\[
\begin{split}
\mathsf{E}_{h,\rho}(u)&=(\mathsf{A}_{h,\rho}u,u)_0 \\
&=\sum_{j \in J_h} \lambda_{j,h}\abs{\hat{u}_j}^2\\
&=\sum_{j \in J_h} \abs{\hat{u}_j}^2-\sum_{j \in J_h} (1-\lambda_{j,h})\abs{\hat{u}_j}^2\\
&\leq \norma*{u}_0^2-\norma*{S_{h,\rho}u}_0^2.
\end{split}
\]
Therefore, by \eqref{Nashinequality}, for all $u \in V_h$, 
\begin{equation}\label{mainesttve}
\norma*{u}_0^{2+1/d} \leq Ch^{-2}(\norma*{u}_0^2-\norma*{S_{h,\rho}u}_0^2+h^2\norma*{u}^2_0)\norma*{u}_{L^1}^{1/d}.
\end{equation}
Furthermore, by \eqref{estk2h} and $S_{h,\rho}^k-\Pi_0=K_{1,h}^k+K_{2,h}^k$, we get that there exists $C_2>0$ such that, for all $h \in ]0,h_0]$ and $k \geq B\log(1/h)$ 
\[
\norma*{K_{1,h}^k}_{L^\infty \rightarrow L^\infty} \leq C_2
\]
and so, since $K_{1,h}$ is self-adjoint on $L^2(M,d\mu)$ 
\[
\norma*{K_{1,h}^k}_{L^1 \rightarrow L^1} \leq C_2.
\]
We now fix $s \approx B\log(1/h)$, $g \in L^2(M,d\mu)$
such that $\norma*{g}_0 \leq 1$ and define the sequence $(\ell_k)_{k \in \mathbb{N}_0}$ by 
\[
\ell_k:=\norma*{K_{1,h}^{k+s}g}_0^2.
\] 
Then, for all $k \in \mathbb{N}_0$, $0 \leq \ell_{k+1}\leq \ell_k$, 
and since $K_{1,h}^{k+s}g \in V_h$, using \eqref{mainesttve} we obtain
\begin{equation}
\begin{split}
\ell_k^{2+1/d} &\leq Ch^{-2}(\ell_k-\ell_{k+1}+h^2\ell_k)\norma*{K_{1,h}^{k+s}g}_{L^1}^{1/d}\\
&\leq CC_2^{1/d}h^{-2}(\ell_k-\ell_{k+1}+h^2\ell_k)\norma*{g}_{L^1}^{1/d}.
\end{split}
\end{equation}
By \cite{DLM1} (pp. 258-259), there exists
a constant $D=D(C,C_2,d)>0$ such that,
for all $0 \leq k \leq h^{-2}$ and $s \approx B\log(1/h)$ one has 
\[
\norma*{K_{1,h}^{k+s}g}_0 \leq \Bigl(\frac{Dh^{-2}}{1+k}\Bigr)^d\norma*{g}_{L^1},
\]
and since $K_{1,h}$ is self-adjoint on $L^2(M,d\mu)$, by duality, we get 
\[
\norma*{K_{1,h}^{k+s}g}_{L^\infty} \leq \Bigl(\frac{Dh^{-2}}{1+k}\Bigr)^d\norma*{g}_0.
\]
Thus, if $N \approx h^{-2}$, one has
\[
\norma*{K_{1,h}^{N+s}g}_{L^\infty}\leq C''\norma*{g}_0,
\]
for some $C''>0$, and consequently, 
\[
\norma*{K_{1,h}^{N+s+m}g}_{L^\infty}\leq C'' (1-h^2\lambda_{1,h})^m \norma*{g}_0
\]
for $m\geq0$ and $N \approx h^{-2}$.
Finally, since $h^2\lambda_{1,h}=g_{h,\rho}$ 
and $0 \leq (1-y)^m \leq e^{-mt}$ for $y \in [0,1]$ we get 
that, for all $k \geq h^{-2}+N+s$,
\[
\norma*{K_{1,h}^k}_{L^\infty \rightarrow L^\infty}\leq C e^{-(k-(N+s))g_{h,\rho}}=Ce^{(N+s)g_{h,\rho}} e^{-kg_{h,\rho}} \leq C'e^{-kg_{h,\rho}}, 
\]
for some $C,C'>0$ universal constants. The proof of the theorem is complete. 
\end{proof}

\appendix
\section{}
\subsection{Technical Lemmas}
The aim of this section is to provide
the technical Lemmas that we use throughout the work.
We begin with the following equivalent characterization
of subellipticity due to Fefferman and Phong (see \cite{FP2}, Theorem 1). 
\begin{proposition}\label{prop.tec.1}
Let $A$ be a second-order differential operator on $M$ 
(see Definition \ref{diffM}) that satisfies \ref{hp1}, \ref{hp2} and \ref{hp3}. Then 
\begin{equation}
\label{SubOp1}
\norma*{Au}_0+C \norma*{u}_0 \geq c \norma*{u}_{2\varepsilon}, \quad \forall u \in C^\infty (M),
\end{equation}
for some $c,C>0$ and some $0<\varepsilon <1$ if and only if
\begin{equation}
\label{SubOp2}
\mathrm{Re}(Au,u)_0+C' \norma*{u}_0^2 \geq  c' \norma*{u}^2_{\varepsilon}, \quad \forall u \in C^\infty(M),
\end{equation}
for some constants $c',C'>0$.
\end{proposition}
The next result that applies to 
the local expression of our operator is very useful (see \cite{OR}, pp. 64-66).
\begin{proposition}
\label{propmatrixsub}
Let $\mathsf{A}_2(x)=(a_{ij}(x))_{i,j} \geq 0,$ for all $x \in \mathbb{R}^n$.
Then, for all $ 1\leq j\leq n$, one has 
\[
\abs*{\sum_{i=1}^na_{ij}(x) \xi_i}^2 \leq 2a_{jj}(x) \sum_{i,j=1}^n a_{ij}(x) \xi_i \xi_j,  \quad \forall (x, \xi) \in \mathbb{R}^n_x \times \mathbb{R}^n_\xi.
\]
\end{proposition}
We conclude this section with the following result.
\begin{lemma}\label{lemma.Fourier}
Let $\varphi \in \mathscr{S}(\mathbb{R}^{n})$ such that $\int_{\mathbb{R}^n}\varphi(x)dx=0$.
Then there exist $\varphi_1,...,\varphi_n \in \mathscr{S}(\mathbb{R}^{n})$
such that 
\[
\varphi=\sum_{j=1}^n \frac{\partial}{\partial x_j} \varphi_j.
\]  
\end{lemma}
\begin{proof}
We first note that, by applying the Fourier transform,
since $\int_{\mathbb{R}^n}\varphi(x)dx=0$,
we may reduce to finding $\phi_1,...,\phi_n \in \mathscr{S}(\mathbb{R}^n)$
such that  
\[
\sum_{j=1}^n\xi_j \phi_j=\hat{\varphi}, 
\]
with $\hat{\varphi}(0)=0$. 

\noindent Since
\[
\hat{\varphi}(\xi)=\sum_{j=1}^n \xi_j\int_0^1\frac{\partial \hat{\varphi}}{\partial \xi_j}(t\xi)dt,
\]
for $j=1,\dots n$, defining
\[
\tilde{\phi}_j(\xi):=\int_0^1 \frac{\partial \hat{\varphi}}{\partial \xi_j}(t\xi)dt
\]
and
\[
\phi_j(\xi):=\chi(\abs*{\xi}^2) \tilde{\phi}_j(\xi)+(1-\chi(\abs*{\xi}^2))\hat{\varphi}(\xi) \frac{\xi_j}{\abs*{\xi}^2},
\]
where $\chi \in C_c^\infty (\mathbb{R}), \ 0 \leq \chi \leq 1$,  $\chi(t) \equiv 1$ if $\abs*{t} \leq 1$ and $\chi(t) \equiv 0$ if $\abs*{t} \geq 2$.

\noindent In this way $\phi_j \in \mathscr{S}(\mathbb{R}^n)$ for all $j=1,...,n$ and
\begin{equation}
\begin{split}
\sum_{j=1}^n\xi_j\phi_j(\xi)&=\chi(\abs*{\xi}^2)\sum_{j=1}^n \xi_j\tilde{\phi}_j(\xi)+(1-\chi(\abs*{\xi}^2))\sum_{j=1}^n \hat{\varphi}(\xi)\frac{\xi_j^2}
{\abs*{\xi}^2} \\
&=\chi(\abs*{\xi}^2)\hat{\varphi}(\xi)+(1-\chi(\abs*{\xi}^2))\hat{\varphi}(\xi)=\hat{\varphi}(\xi).
\end{split}
\end{equation}
The result follows by setting
$\varphi_j:=\mathscr{F}_{\xi \rightarrow x}^{-1}(\phi_j)$,
for $j=1,...,n$.
\end{proof}
\subsection{Differential Operators on Manifolds}
In this section, we review some useful properties of operators on a smooth manifold
$M$. We begin by defining what we mean by a differential operator.
In the first place, we specify what we mean by a differential operator on a smooth manifold.
\begin{definition}\label{diffM}
Let $M$ be a smooth manifold and let $A:C^\infty(M)\rightarrow C^\infty(M)$ be a linear operator.
We say that \textit{$A$ is a differential operator on $M$} 
of order $m$ with smooth coefficients on $M$,
if, for each $\alpha \in \mathscr{A}$, 
the operator $A_\alpha:C^\infty(\tilde{U}_\alpha)\rightarrow C^\infty(\tilde{U}_\alpha)$
defined by
\[
A_\alpha f=A(f \circ \varphi_\alpha) \circ \varphi_\alpha^{-1}, \quad f \in C^\infty(\tilde{U}_\alpha),
\]
is a differential operator of order $m$ on $\tilde{U}_\alpha\subseteq \mathbb{R}^n$.
The operator $A_\alpha$ is called \textit{a local expression of $A$ on $\tilde{U}_\alpha$}.  
\end{definition}
\begin{remark}\label{rmk.princdiff}
We may define the principal symbol
of a differential operator 
\[
A:C^\infty(M)\rightarrow C^\infty(M)
\]
of order $m$ on $M$ as a well-defined function
$a(x,\xi) \in C^\infty(T^\ast M \setminus 0)$ (where $T^\ast M \setminus 0$ is the cotangent space with the zero section removed, i.e. $T^\ast M \setminus (M \times \lbrace 0 \rbrace)$) that
corresponds to the principal symbol of its local
expression in coordinates (cf. \cite{HoV1}, p. 151 or \cite{GG}, p. 205)
\[
\tilde{a}(\tilde{x},\tilde{\xi})=\sum_{\abs{\gamma}=m}a_{\gamma}(\tilde{x})\tilde{\xi}^\gamma \in C^\infty(\tilde{U}_\alpha \times \mathbb{R}^n).
\]
\end{remark}
\begin{remark}\label{princvecfield}
In a similar fashion, if $X \in \Gamma(M,TM)$
and $X=\sum_{j=1}^n X_j\partial/\partial x_j$ 
in coordinates on $U_\alpha \subseteq M$, then the principal symbol of $X$
can be defined as the function $X(x,\xi) \in C^\infty(T^\ast M \setminus 0)$
that corresponds locally to the function
\[
\tilde{X}(\tilde{x},\tilde{\xi})=\sum_{j=1}^n X_j(\tilde{x})\xi_j \in C^\infty(\tilde{U}_\alpha \times \mathbb{R}^n)
\]
(to be consistent with the previous notations this actually should be the symbol of $i^{-1}X$). 
\end{remark}

At this point we recall the notion of a divergence-free 
of a vector field with respect a measure.
Let $X$ be a vector field on $M$ smooth manifold and let
$d\nu$ be a non-degenerate volume form on $M$.
Since the Lie derivative $\mathscr{L}_Xd\nu$ (see Lee \cite{L1}, p. 321)
of this volume form
is still a smooth $n$-form, 
there exists a unique smooth function $f \in C^\infty(M)$
such that $\mathscr{L}_Xd\nu=f d \nu$.

\begin{definition}
\textit{The divergence of a vector field $X \in \Gamma(M,TM)$
with respect a volume form $d\nu$ 
(or with respect the measure $\mu$ associated with it)}
is defined as the unique function $f \in C^\infty(M)$
that satisfies $\mathscr{L}_Xd\nu=f d \nu$
and we denote it by $\mathrm{div}_{\nu} X:=f$. 
\end{definition}

\begin{definition}\label{def.divfree}
Let $d\nu$ be a volume form on a manifold $M$
and let $X$ be a vector field on $M$. 
We say that \textit{$X$ is divergence-free with respect to $\nu$} 
if $\mathrm{div}_{\nu} X=0$.
\end{definition}

\begin{remark}\label{rmk.divfree}
An easy computation shows that, if locally $d\nu(x)=\gamma (x) dx_1 \wedge \dots \wedge dx_n$, with $\gamma \neq 0$, 
and $X(x)=\sum_{i=1}^n a_i(x) \, \partial/\partial x_i$, 
\[
\mathrm{div}_{\nu} X=\sum_{i=1}^n \frac{1}{\gamma}\frac{\partial}{\partial x_i}(\gamma a_i).
\]
Then, locally 
\begin{equation}
\label{localdivergence}
\mathrm{div}_{\nu} X=0 \iff \sum_{i=1}^n \frac{\partial}{\partial x_i}(\gamma a_i)=0.
\end{equation}
\end{remark}

\subsection{Basic Notions of Spectral Theory}\label{app.spec}
The purpose of this section is to provide the basic notions
of spectral theory on a separable Hilbert space that we use in the work.
To do that we refer to \cite{He} and \cite{RS}.
In what follows $(H, \langle \cdot, \cdot \rangle_H)$
is a Hilbert space with norm $\norma*{\cdot}_H$
induced by the scalar product.
\begin{definition}
An \textit{unbounded operator (or simply an operator) on $H$} is a pair $(T,D(T))$ where $D(T)$
is a subspace of $H$ and $T:D(T) \subseteq H \rightarrow H$
is a linear map.
If $T$ is a continuous map from $D(T)$ (with respect the induced topology)
to $H$ then $T$ is said to be \textit{a bounded operator}.
\end{definition}

In what follows for a general unbounded operator $(T,D(T))$
we will write $T:D(T)\subseteq H \rightarrow H$ and we always suppose
that $D(T)\subseteq H$ is a dense subspace.
It is important to emphasize the case for which $D(T)=H$.
In this case the set of bounded operators from $H$
to itself is denoted by $\mathscr{L}(H)$
and it is possible to endow such a space with the structure of Banach space.
More generally, given $(B_1, \| \cdot \|_{B_1})$ 
and $(B_2, \| \cdot \|_{B_2})$ Banach spaces, we denote by $\mathscr{L}(B_1,B_2)$ the set
of continuous operator from $B_1$ to $B_2$ and for $T:B_1 \rightarrow B_2$ 
we define the norm
\[
\norma*{T}_{B_1\rightarrow B_2}:=\sup_{u \neq 0} \frac{\norma*{Tu}_{B_2}}{\norma*{u}_{B_1}}.
\]
that in what follows we simply denote by $\norma*{T}:=\norma*{T}_{H \rightarrow H}$.

We next give the notion of a graph that is very useful
to study unbounded operators, and consequently
the notion of \textit{closed} and \textit{closable} operator. 

\begin{definition}
Let $T:D(T)\subseteq H \rightarrow H$ be an unbounded operator. The \textit{graph of $T$} is the set 
\[
G(T):=\lbrace (x,y) \in H \times H; \; x \in D(T), \ y=Tx \rbrace 
\]
\end{definition}

\begin{definition}
An unbounded operator $T:D(T)\subseteq H \rightarrow H$
is said to be \textit{closed} if its graph is a closed set in $H\times H$, 
with respect the product topology. 
\end{definition}

\begin{definition}
An unbounded operator $T:D(T)\subseteq H \rightarrow H$
is said to be \textit{closable} if its closure in $H\times H$ is a graph.
\end{definition}

Equivalently an unbounded operator $(T,D(T))$ is closable
if there exists a closed extension,
i.e. if there exists $T_1:D(T_1)\subseteq H \rightarrow H$
such that $T_1$ is closed,
$D(T)\subseteq D(T_1)$ and $T_1u=Tu$ for all $u \in D(T)$. 

Note that bounded operators in general are not excluded from 
the class of unbounded operators. However, 
it is a common use to study first the properties
of bounded operators from $H$ to $H$ and then 
the properties for general unbounded operators.
Let us investigate now the notion of the \textit{adjoint}
that will be fundamental to study the spectrum of an operator.

If $T \in \mathscr{L}(H)$, \textit{the Hilbertian adjoint}
$T^\ast$ is defined by the relation
\begin{equation}\label{bound.adj}
\langle Tu, v \rangle_H=\langle u, T^\ast v \rangle_H
\end{equation}
More precisely, for $v \in H$,
the map $H \ni u \mapsto \langle Tu,v \rangle_H$
defines a continuous linear map in $H$,
thus, by Riesz theorem (see \cite{He} Theorem 3.1),
there exists $w=:T^\ast v \in H$ such that \eqref{bound.adj} is satisfied.
In general we have the following definition (that includes the case of unbounded operators).

\begin{definition}
Let $T:D(T)\subseteq H \rightarrow H$ be an unbounded operator.
We define the \textit{domain of the adjoint} as
\[
D(T^\ast):= \lbrace  v \in H; \abs*{\langle Tu,v \rangle_H} \leq C \norma*{u}_H, \ \text{for some $C>0$} \rbrace
\] 
and the \textit{adjoint} $T^\ast:D(T^\ast)\subseteq H \rightarrow H$, defined 
through the Riesz Theorem, by setting $w=:T^\ast v$,
for $v \in D(T^\ast)$, where $w \in H$ is the unique element that satisfies 
\[
\langle Tu,v \rangle_H=\langle u , w \rangle_H, \quad \forall u \in D(T).
\]
\end{definition}

From this definition immediately follows that if $D(T)=H$ and $T$ is bounded,
its adjoint satisfies \eqref{bound.adj}. 
\begin{definition}
Let $T:D(T) \subseteq H \rightarrow H$. The operator $T$ is said to be \textit{symmetric} if 
\[
\langle Tu, v \rangle_H= \langle u, Tv \rangle_H, \quad \forall u,v \in D(T). 
\]
\end{definition}

\begin{definition}
An operator $T:D(T)\subseteq H \rightarrow H$ is said to be \textit{self-adjoint}
if $D(T)=D(T^\ast)$ and $Tu=T^\ast u$,
for all $u \in D(T)$, i.e. if it is symmetric 
and its domain coincide with the domain of the adjoint. 
\end{definition}

Finally, in the case of $H=L^2(M,d\mu)$
it is useful to recall the notion of formal-adjoint. 
\begin{definition}
Let $A:C^\infty(M) \subseteq L^2(M,d\mu) \rightarrow L^2(M,d\mu)$ be a linear operator. \textit{The formal adjoint} of $A$ is the operator $A^\times$ defined (again by the Riesz Theorem) by the relation 
\[
(Au,v)_0=(u,A^\times v)_0, \quad u,v \in C^\infty(M).
\] 
In this context we thus say that \textit{$A$ is formally self-adjoint} if $A=A^\times$. 
\end{definition} 

We now go back to our abstract framework
and continue to state the following properties
for an operator $T:D(T) \subseteq H \rightarrow H$. 

\begin{definition}
Let $T:D(T)\subseteq H \rightarrow H$ be a symmetric operator.
We say that $T$ is \textit{positive}, and write $T \geq 0$, if 
\[
\langle Tu,u \rangle_H \geq 0, \quad \forall u \in D(T).
\]
In particular, if $T$ and $S$ are symmetric,
we write $T \geq S$, if $T-S \geq 0$. 
\end{definition}

\begin{definition}
Let $T:D(T)\subseteq H \rightarrow H$ be a symmetric operator. We say that $T$ is \textit{semibounded}  if there exists a constant $C>0$ such that 
\[
\langle Tu,u \rangle_H \geq -C \norma*{u}^2_H, \quad \forall u \in D(T).
\]
\end{definition}

We are now in a position to deal with 
the spectrum and its characterizations first in the case of a
bounded operator defined on $H$ and then for an unbounded closed operator.

\begin{definition}
Let $T \in \mathscr{L}(H)$ be a bounded operator. The \textit{resolvent set} of $T$ is defined as 
\[
\varrho(T)=\lbrace \lambda \in \mathbb{C}; \ \textit{$(T-\lambda I)$ is bijective from $H$ to $H$} \rbrace.
\]
The \textit{spectrum} of $T$ is defined as 
\[
\mathrm{Spec}(T)=\mathbb{C}\setminus \varrho(T).
\]
\end{definition}

It is easy to verify that the resolvent set
is an open set and then the spectrum is a closed
set of $\mathbb{C}$.
Moreover, it is useful to clarify now
what we mean to be an \textit{eigenvalue}.

\begin{definition}
Let $\lambda \in \mathbb{C}$. If the nullspace $N(T-\lambda \cdot I)\neq \lbrace 0 \rbrace$, then $\lambda$ is said to be an eigenvalue of $T$. The set of all eigenvalues of $T$ is called \textit{the point spectrum of $T$} and it is denoted by $\mathrm{Spec}_p(T) \subseteq \mathrm{Spec}(T)$.  
\end{definition}

The spectrum of a bounded operator satisfies the following properties
(see \cite{He}, Proposition 6.15).   
\begin{proposition}\label{prop.specball}
Let $T \in \mathscr{L}(H)$.
Then $\mathrm{Spec}(T)$ is a compact set and
\[
\mathrm{Spec}(T) \subseteq \overline{B(0,\norma*{T})}.
\]
\end{proposition}
In the case of a self-adjoint operator we have a more precise characterization
(see \cite{He}, p. 59).
\begin{proposition}\label{prop.SA}
Let $T \in \mathscr{L}(H)$ be a self-adjoint operator.
Then the spectrum of $T$ is real and it is contained in $[r,R]$,
where $r:=\inf_{0 \neq u \in H} \langle Tu,u \rangle_H/\norma*{u}^2_H$
and $R:=\sup_{0 \neq u \in H} \langle Tu,u \rangle_H /\norma*{u}^2_H$.
Moreover, $r$ and $R$ belong to the spectrum of $T$.
Finally, if $T$ is also positive then $R=\norma*{T}$.
\end{proposition}
We may read these properties also in terms of the spectral radius (see \cite{RS}, p. 192). 
\begin{definition}
\textit{The spectral radius} of $T \in \mathscr{L}(H)$ is defined as
\[
r_{\mathrm{spec}}(T):=\sup_{\lambda \in \mathrm{Spec}(T)} \abs*{\lambda}.
\]
\end{definition}
\begin{theorem}\label{theospecradius}
Let $T \in \mathscr{L}(H)$. Then $\lim_{n \rightarrow +\infty}\norma*{T^n}^{1/n}$ exists and is equal to $r_{\mathrm{spec}}(T)$. Moreover, if $T$ is self-adjoint, then $r_{\mathrm{spec}}(T)=\norma*{T}$. 
\end{theorem}
We now define the spectrum for a
general closed operator (not necessarily bounded),
We will characterize it 
in terms of \textit{discrete spectrum} 
and \textit{essential spectrum} and
recall the main properties of 
the spectrum of unbounded self-adjoint operators.
To do that we follow Helffer \cite{He}.
 
\begin{definition}
Let $T:D(T)\subseteq H \rightarrow H$ be a closed operator.
The \textit{resolvent set} of $T$,
denoted again by $\varrho(T)$ is the set of all $\lambda \in \mathbb{C}$ 
such that the range of $(T-\lambda I)$ is equal to $H$
and such that $(T-\lambda I)$ admits a continuous inverse,
denoted by $R(\lambda)$, whose range is therefore included in $D(T)$,
such that $R(\lambda)(T-\lambda I)=I_{D(T)}$ and $(T-\lambda I)R(\lambda)=I_H$.
The \textit{spectrum} of $T$, denoted by $\mathrm{Spec}(T)$, is then again defined by $\mathrm{Spec}(T)=\mathbb{C}\setminus \varrho(T)$. Also in this case if the nullspace $N(T-\lambda \cdot I)\neq \lbrace 0 \rbrace$, with $\lambda \in \mathbb{C}$, then $\lambda$ is said to be an eigenvalue of $T$.
\end{definition}

We have the following property (see \cite{He}, Proposition 11.2).

\begin{proposition}\label{posselfadj}
Let $T:D(T):H\rightarrow H$ be an unbounded self-adjoint operator. Then $T\geq 0$ if and only if $\mathrm{Spec}(T)\subseteq [0,+\infty[$.
\end{proposition}

\begin{definition}
Let $T:D(T)\subseteq H \rightarrow H$ be a self-adjoint operator. We call \textit{the discrete spectrum} denoted by $\mathrm{Spec}_{\mathrm{disc}}(T)$ the set of $\lambda \in \mathbb{C}$
such that $\lambda$ is an isolated point in $\mathrm{Spec}(T)$
and an eigenvalue of finite multiplicity. 
\textit{The essential spectrum} is the complementary set of the discrete spectrum, i.e. 
\[
\mathrm{Spec}_\mathrm{ess}(T)=\sigma(T)\setminus \mathrm{Spec}_{\mathrm{disc}}(T).
\]
\end{definition}

We next recall the following variational description of the spectrum,
the so-called \textit{max-min principle} (see \cite{He}, Theorem 11.7).
\begin{theorem}\label{minmax}
Let $T:D(T)\subseteq H \rightarrow H$ be a self-adjoint semibounded operator. For each $n \geq 1$ we set 
\[
\mu_n(T)=\sup_{\psi_1,\psi_2, \dots, \psi_{n-1}} \; \inf_{	\begin{aligned} \phi & \in [\mathrm{span}(\psi_1,\dots,\psi_{n-1})]^\perp\\
&\phi \in D(T) \ \mathrm{and} \ \norma*{\phi}_H=1 \end{aligned}}\langle T\phi,\phi \rangle_H. 
\] 
Then either, \\
$(a)$ $\mu_n(T)$ is the $n$-th eigenvalue when the eigenvalues are ordered in increasing order (counting the multiplicity) and $T$ has a discrete spectrum in $]-\infty,\mu_n(T)[$;\\
or\\
$(b)$ $\mu_n(T)$ corresponds to the bottom of the essential spectrum. In this case, we have $\mu_j(T)=\mu_n(T)$ for all $j \geq n$. 
\end{theorem}
This theorem leads to the following proposition.
\begin{proposition}\label{propweyl1}
Let $T:D(T)\subseteq H \rightarrow H$ be a semibounded self-adjoint operator on a Hilbert space $H$
and let $a>0$ be a positive constant such that
\begin{equation}\label{weylformula}
\sharp \lbrace j, \mu_j(T) \in [0,a[\rbrace <+\infty.
\end{equation}
Then $\mathrm{Spec}({T})\cap [0,a[$ is discrete and all the eigenvalues of $T$ in $[0,a[$, repeated
according to multiplicity, are all found by the max-min principle.
\end{proposition}

Finally, the last tool that we need, 
in order to study the velocity of the convergence
to the equilibrium is the notion of \textit{spectral gap}.

\begin{definition}\label{def.specgap}
Let $T$ be a bounded operator of norm $1$ such that $T1=0$.
\textit{The spectral gap} $g \in \mathbb{R}$ associated with some bounded operator $T$,
is the best constant $1/g>0$ that satisfies the inequality
\[
\norma*{u}^2_{L^2}- (u,1)_0^2 \leq \frac{1}{g}(T u, u)_0, \quad \forall u \in L^2(M,d\mu).
\]
\end{definition}

\begin{remark}
\label{rmkspecgap}
The spectral gap $g$ of an operator $T$
acting on $L^2(M,d\mu)$ as above,
can be interpreted as the minimum eigenvalue. 
In fact, observe first
\[
\norma*{u}_0^2-(u,1)_0^2=\int_M \abs*{u}^2 d\mu - \Bigl(\int_M u d\mu\Bigr)^2=\mathrm{Var}_\mu(u),
\]
from which it follows 
\[
g=\inf \Set{ \frac{(T u,u)_0}{\mathrm{Var}_\mu(u)}, \quad \mathrm{Var}_\mu(u)\neq 0}.
\]
Then 
\begin{equation}
\begin{split}
&\inf \Set{\frac{(T u,u)_0}{\mathrm{Var}_\mu(u)}; \quad \mathrm{Var}_\mu(u)\neq 0} \\
&=\inf \Set{\frac{(T u,u)_0}{\int_M \abs*{u}^2 d\mu - (\int_M u d\mu)^2 }; \quad u \neq \mathrm{cost}} \\
&=\inf \Set{(T u,u)_0; \quad \norma*{u}_0^2=1, \; \int_M u d\mu=0}.
\end{split}
\end{equation}
where the last equality follows by writing $u=u_1+u_2$, where $u_2$ is the orthogonal projection of $u$
on the space of constant functions, and using the fact that $Tu_2=0$.
Finally, in view of the min-max characterization 
(see for instance \cite{P}, Theorem 4.1.1) we have 
\[
\lambda_1=\inf \Set{(T u,u)_0; \quad \norma*{u}_0=1},
\]
so that, since $T1=0$,
\begin{equation}
\lambda_1=\inf \Set{(T u,u)_0; \quad \norma*{u}_0^2=1, \; \int_M u d\mu=0}=g.
\end{equation}
\end{remark} 

\subsection{Exponential Map and Some Related Topics}\label{app.exp}
Let $(M,g)$ be a smooth Riemannian manifold with the Levi-Civita connection.
The aim of this section is to construct an atlas through the \textit{exponential map}.
For such arguments we refer to the book of Lee \cite{L1}.
In the first place it is useful to recall that \textit{the Riemannian distance} $d_g(x,x')$ 
between two points $x,x' \in M$ is defined as 
\[
\begin{split}
d_g(x,x'):= \inf \Bigl\lbrace L_g(\gamma); \quad \gamma:[a,b]\subseteq \mathbb{R} \rightarrow M \ & \text{piecewise regular}, \\
&\gamma(0)=x, \quad \gamma(1)=x' \Bigr\rbrace,
\end{split}
\]
where, if $\gamma:[a,b] \rightarrow M$ is a piecewise regular curve
(i.e. if there exists a subdivision $a=t_0<t_1<\dots t_n=b$ such that,
for each $i$, the restriction $\gamma_{\bigl|_{[t_{i-1},t_i]}}$ 
is $C^\infty$ and $\gamma'_{\bigl|_{[t_{i-1},t_i]}}$ never vanish) 
\[
L_g(\gamma):= \int_a^b \abs{\dot{\gamma}(t)}^g_{{\gamma(t)}}dt. 
\]
It is well known that each initial point $x \in M$
and each initial velocity vector $v_x \in T_xM$ determine a unique maximal geodesic $\gamma_v$, with $v=(x,v_x)$.
Then we define the domain of the exponential map as 
\begin{center}
$\mathscr{E}=\lbrace v \in TM;$  $\gamma_v$ is defined on an interval containing $[0,1] \rbrace$. 
\end{center}
The exponential map is thus defined as the map 
\[
\begin{split}
\mathrm{exp}: \ & \mathscr{E} \rightarrow M, \\
& v \mapsto \gamma_v(1).
\end{split}
\]
From it we also define the exponential map at point $x$
as the restriction of the exponential map
to the set $\mathscr{E}_x:=\mathscr{E}\cap T_xM$,
i.e.
\[
\begin{split}
\mathrm{exp}_x: \  \mathscr{E}_x \subseteq  T_xM & \ \rightarrow M \\
& v_x \mapsto \gamma_v(1).
\end{split}
\]
For every $x \in M$, $\mathrm{exp}_x$ is smooth and has the property that  
\[
\mathrm{d}(\mathrm{exp}_x)_{0}:T_{0}(T_xM)\cong T_xM \rightarrow T_xM
\]
(with $0$ the origin in $T_xM$) is the identity map.
Then, by the inverse function theorem, there exists
a neighborhood of the origin $U'_{0}$ in $T_xM$
and a neighborhood $U_x$ of $x$ in $M$ such that $\mathrm{exp}_x:U'_{0}\rightarrow U_x$
is a diffeomorphism.
Note also that every choice of basis $(b_i)_i$ for $T_xM$
determines an isomorphism $B:\mathbb{R}^n \rightarrow T_xM$.
By combined these two diffeomorphisms 
it is possible to get a smooth coordinate map 
\begin{equation}\label{varphiexp}
\varphi_x:B^{-1} \circ \mathrm{exp}_{\bigl|U'_{0}}^{-1}:U_x \rightarrow \mathbb{R}^n.
\end{equation}
and such $U_x\subseteq M$ is called \textit{a normal neighborhood of $x$}. 
  
Let us denote now by $B(0,a)\subseteq T_xM$
the ball centered at $0 \in T_xM$ of radius $a>0$
with respect to the metric $g_x$ and 
\begin{center}
$r_g(x):=\sup\lbrace a>0;$  $\mathrm{exp}_x$ is a diffeomorphism from $B(0,a) \subseteq T_xM$ onto its image $\rbrace$.
\end{center}
The injective radius of $M$ is then defined as
\begin{equation}\label{rgM}
r_g(M):=\inf_{x \in M} r_g(x).
\end{equation}

\begin{theorem}
If $(M,g)$ is a compact Riemannian manifold
then $r_g(M)$ is strictly positive. 
\end{theorem}
Thanks to this result, setting $r:=r_g(M)>0$
we may consider for every $x \in M$ the map 
\[
\mathrm{exp}_x:B(0,r) \rightarrow B(x,r)
\]
where $B(x,r):=\mathrm{exp}_x(B(0,r))$ is a normal
neighborhood of $x$ called \textit{geodesic ball} centered at $x$ in $M$. 
In this way we have an open cover of $M$ 
\[
\mathscr{U}'=\lbrace (B(x,r),\varphi_x); \ \ x \in M \rbrace,
\]
where $\varphi_x:=B^{-1}\circ \mathrm{exp}_{\bigl|_{B(0,r)}}^{-1}$
is defined as in \eqref{varphiexp}.

Moreover, if $B(x,r)\cap B(x',r)\neq \emptyset$, 
\[
\varphi_x \circ \varphi_{x'}^{-1}:\varphi_{x'}(B(x,r)\cap B(x',r))\rightarrow \varphi_x(B(x,r)\cap B(x',r))
\]
is a diffeomorphism by definition.
So the open cover $\mathscr{U}'$ is an atlas for $M$,
contained in the maximal atlas $\mathscr{U}$. 
Furthermore, we also have that the geodesic balls
we have just defined are equal to the metric balls on $M$. 
\begin{theorem}\label{BandB_g}
If $(M,g)$ is a Riemannian manifold, $x \in M$ one has
\[
B(x,a)=B_g(x,a):=\lbrace y \in M; \ d_g(x,x')<a \rbrace, 
\]
for every $0<a<\mathrm{inj}(x)$. 
\end{theorem}
The last tool we need is the Lebesgue number given as follows.
\begin{lemma}\label{Lebesguelemma}
Let $(X,d)$ be a compact metric space and
$\mathscr{U}=\lbrace U_\alpha\rbrace_{\alpha \in \mathscr{A}}$ an open covering of $X$.
Then there exists $\delta>0$ such that every set
of diameter less than $\delta$ is contained in some $U_\alpha \in \mathscr{U}$. 
\end{lemma}
In our context, since $M$ is compact and $(M,g)$ with the geodesic distance
is a compact metric space, we have a finite cover 
\[
\lbrace B_g(x_\ell,r); \ \ell=1,\dots,m \rbrace,
\]
for some $m \in \mathbb{N}$, and by Lemma \ref{Lebesguelemma}
there exists $\delta>0$, such that for every $x \in M$ the
geodesic ball $B_g(x,\delta)$ is contained
in $B_g(x_\ell,r)$, for some $\ell \in \lbrace 1,...,m \rbrace$.

\vspace{2cm}

\noindent\textbf{Acknoledgments:} The author would like to thank Alberto Parmeggiani for his constant support throughout the preparation of this work, as well as for many helpful discussions.

\end{document}